\documentclass[onecolumn, twoside, a4paper, 11pt]{article}

\usepackage[utf8]{inputenc}
\usepackage[T1]{fontenc}

\usepackage{fullpage}

\usepackage{newpxtext,newpxmath}

\usepackage[detect-weight=true, binary-units=true]{siunitx}

\usepackage{amsmath}
\usepackage{amssymb}
\usepackage{bm}
\usepackage{nicefrac}

\newtheorem{remark}{Remark}

\usepackage{graphicx}
\usepackage{caption}
\usepackage{subcaption}
\usepackage{adjustbox}

\usepackage{tikz}
\usetikzlibrary{positioning,shapes,arrows,calc,intersections}
\usepackage{pgfplots}
\usepgfplotslibrary{dateplot}
\pgfplotsset{compat=1.8}

\usepackage{nicefrac}
\usepackage{authblk}

\newcommand{\dd}{{\;\mathrm{d}}}

\begin{document}

\title{Fast divergence-conforming reduced basis methods for steady Navier-Stokes flow}
\author[1]{E.~Fonn}
\author[2]{E.H.~van~Brummelen}
\author[1,3]{T.~Kvamsdal}
\author[1]{Adil Rasheed}
\affil[1]{SINTEF Digital, Trondheim, Norway}
\affil[2]{Department of Mechanical Engineering, Eindhoven University of Technology}
\affil[3]{Department of Mathematical Sciences, Norwegian Univ.~of Science and Technology}

\maketitle

\begin{abstract}
  Reduced-basis methods (RB methods or RBMs) form one of the most promising techniques to
  deliver numerical solutions of parametrized PDEs in real-time performance
  with reasonable accuracy. For incompressible flow problems, RBMs based on LBB
  stable velocity-pressure spaces do not generally inherit the stability
  of the underlying high-fidelity model and, instead, additional stabilization techniques
  must be introduced. One way of bypassing the loss of LBB stability in the RBM is to
  inflate the velocity space with supremizer modes.
  This however
  deteriorates the performance of the RBM in the performance-critical online stage, as
  additional DOFs must be introduced to retain stability, while these DOFs do not effectively
  contribute to accuracy of the RB approximation. In this work we consider a velocity-only
  RB approximation, exploiting a solenoidal velocity basis. The solenoidal reduced basis
  emerges directly from the high-fidelity velocity solutions in the offline stage. By means
  of Piola transforms, the solenoidality of the velocity space is retained under geometric
  transformations, making the proposed RB method suitable also for the investigation of
  geometric parameters. To ensure exact solenoidality of the high-fidelity velocity solutions
  that constitute the RB, we consider approximations based on divergence-conforming compatible
  B-splines. We show that the velocity-only RB method leads to a significant improvement in
  computational efficiency in the online stage, and that the pressure solution can be recovered
  a posteriori at negligible extra cost. We illustrate the solenoidal RB approach by modeling
  steady two-dimensional Navier-Stokes flow around a NACA0015 airfoil at various angles of attack.
\end{abstract}

\section{Introduction}

Conventional methods for simulating partial differential equations include
well-established techniques such as Finite Volume Methods (FVM),
Finite Difference Methods (FDM) and Finite Element Methods (FEM).
Common to all of these methods is the large number of
degrees of freedom that is typically required to accurately model a physical system, often
numbering in the millions or billions. Given the good and well-established
approximation properties of FVM, FDM and in particular FEM, such models are
usually classified as \emph{high-fidelity} models. Problems of this size are
generally not possible to solve in realistic timeframes except on specialized
high-performance computing facilities, which renders such computations very expensive,
and even then they may require several days of computing time. In addition, the
need for dedicated high-end computing resources prohibits the use of high-fidelity
models in time-critical on-site analyses such as biomedical applications, control strategies
and hybrid analytics, in view of the need for additional special infrastructure and
the incurrence of further communication overheads.

The excessive computational complexity of high-fidelity models is also at odds
with the increasing demand for real-time low-cost models for the repetitive solution of
partial differential equations (PDEs)
in many-query scenarios. This is particularly relevant in optimization, control systems,
inverse and inference problems, and uncertainty quantification. Common for many of these
applications is that the PDE in question is
\emph{parametrized} by a suitably small number of input parameters. It is often necessary
to demand solve times in the sub-second regime, which is entirely infeasible with conventional
methods.

Reduced-order modeling (ROM) provides a paradigm to address the aforementioned challenges.
ROM is rapidly developing field~\cite{Bazaz2012rpm}. The general aim of ROM is to replace
the original model with a \emph{reduced model} of very modest computational complexity. Within this
general framework of ideas, one of the most promising is that of \emph{reduced
basis methods} (RBM). This methods dates back to the 1980s with work from
\cite{Almroth1978acg,Almroth1981gfi,Nagy1979mrg,Noor1980rbt,Noor1981bpb,Noor1981rar,Noor1982mln}
and the first theoretical foundations of the method were given by Fink and
Rheinboldt in~\cite{Fink1983ebr,Fink1984sms}.
It was used first used for flow problems in work by Peterson and Gunzburger
\cite{Gunzburger2012fem,Peterson1989rbm}, and much of the more recent work can
be attributed to the theoretical foundations in
\cite{Prudhomme2002mcf,Veroy2003peb}.
Excellent modern introductions can be found in \cite{Quarteroni2016rbm,Haasdonk2017rbm}.

The fundamental concept of reduced basis methods is to formulate the problem on
a function space with very low dimension that is tailored to the solution of the
PDE in the parameter regime of interest. In comparison, while finite-element spaces
generally have well-established asymptotic approximation properties, e.g. density and
quasi-optimality for smooth functions, and the computational cost of constructing a FEM basis
function is generally negligible, the approximation properties of a FEM basis per-degree-of-freedom
(i.e.~for each of the basis function separately) is clearly limited. RB methods seek
to construct an approximation space that has optimal approximation properties with respect
to dimension for a restricted class of functions, viz. the solutions of the PDE in the parameter space
under consideration, under the premise that the cost of constructing the basis is inconsequential. The
latter premise arises because the construction of the basis is carried out a priori in the so-called
offline stage. To do this, some high-fidelity solutions must be computed in advance, called
an \emph{ensemble} of \emph{snapshots}.

RB methods are strictly divided in two stages, viz.~the \emph{offline} and
\emph{online} stages. The offline stage is run only once, and the online stage
is run once for each problem instance to solve (that is to say, for each
parameter query). The goal is then to off-load as much work as possible from the
online to the offline stage, so that the total cost per online execution is very
small in an amortized sense. In other words, work in the offline stage is
considered ``free'' and all efforts will be focused on creating a cheap online stage.

RB methods have been applied in various application areas,
e.g.~optimal control \cite{Ravindran2000roa,Ito2001rbm},
inverse problems \cite{Liu2005iit,Lieberman2010psm},
shape optimization \cite{Manzoni2012sov,Rozza2011rba},
quantum models \cite{Pau2007rbm}
and solid mechanics \cite{Noor1980rbt,Krysl2001dmr}. Development and
application of RB methods to flow problems have been considered in for
instance~\cite{Peterson1989rbm,Gunzburger2012fem,Iollo2000spp,Stabile2017fvp}.
The principal novelty of the present paper is the development of a solenoidal (divergence-free)
velocity-only RB method for incompressible-fluid flow problems with geometric parametrizations, based on
divergence-conforming high-fidelity approximations. The advantage of the solenoidal-velocity RBM is that
the pressure variable
is redundant and can be ignored in the online stage. The absence of the pressure
variable in the RBM in fact carries a two-fold benefit. First, velocity-pressure RBs do not generally inherit
the inf-sup stability of the underlying high-fidelity models. Instability of the RB
velocity-pressure pair can be resolved by extending the velocity space with additional
supremizer modes~\cite{Ballarin2015ssp,Stabile2017fvp}. This however introduces additional DOFs in the
performance-critical online stage, while these DOFs do not effectively contribute to
accuracy of the RB approximation. By incorporating the solenoidality constraint in the RB,
the approximate formulation is released of its mixed-character, and inf-sup stability
considerations for the RBM are voided. Second, the removal of the pressure variable in the RB reduces
the DOF count in the online stage, accelerating the solution process in the time-critical
phase.

The use of solenoidal ROM for incompressible-flow problems has been propounded before;
see, e.g.~\cite{Ballarin2015ssp,Lovgren2006rbe}. However, it appears that the concept was not followed
through for the important scenario of geometric parametrizations, in view of the complexity pertaining
to the required Piola transforms in the ROM. For physical parameters, solenoidal RB is in fact
straightforward. For parameters
related to essential boundary conditions, solenoidal RB is more complicated by the fact that solenoidal
lifts of the boundary data are required. Geometry parametrizations in addition require that the
solenoidality property of the RB is retained under the considered geometric transformations. A
central component in a solenoidal RB method is therefore the Piola transform, which provides
solenoidality-preserving transformations. The Piola transform however leads to complicated
non-polynomial expressions in the ROM, for which it is not generally possible to construct
the affine representations that are required for an efficient online stage. In this work we consider
truncated series expansions of these expressions that do admit an affine representation, and
we assess the computational efficiency of such an approach.

An important additional aspect of the present work is the use of divergence-conforming B-spline
approximations for the high-fidelity simulations.
Such methods were first introduced by Buffa and co-workers~\cite{Buffa2010iae,Buffa2011ias}
and then by Evans and Hughes~\cite{Evans2012dsa,Evans2013idc1,Evans2013idc2,Evans2013idc3}.
Extension to locally refined B-splines~\cite{Dokken2013psl,Johannessen2014ial}
were done in~\cite{Johannessen2015dcd} and to divergence-conforming multiscale turbulence
models in~\cite{Opstal2017idc}.
They rely on formulating a pressure space that
is exactly equal to the divergence of the velocity space. This is enabled by the use of
B-spline basis functions, recently popularized in the field of isogeometric
analysis~\cite{Cottrell2009iat}. The divergence-conforming approximations provide pointwise
solenoidal snapshots, as opposed to the weakly solenoidal high-fidelity solutions provided by, for
instance, conventional Taylor--Hood elements. The Piola transform in turn ensures that this property
is retained in the solenoidal RB method.

The paper is organized as follows:
\begin{itemize}
\item Section~\ref{sec:setting} formulates the classical stationary
  Navier-Stokes problem in a Finite Element setting, introduces the notation
  needed for handling parametrized versions of the problem, and introduces the
  critically important (to Reduced Basis Methods) affine representations in
  Section~\ref{sec:affrep}.
\item Section~\ref{sec:reduction} describes the reduced basis method as applied
  to the the Navier-Stokes problem in detail. One method for pressure
  stabilization is introduced in Section~\ref{sec:stab}.
\item Section~\ref{sec:divconf} introduces the main novelty of the paper, the
  application of divergence-conforming high-fidelity methods to a reduced basis
  method.
\item Section~\ref{sec:airfoil} develops the necessary mathematics for a
  numerical example: flow around a NACA0015 airfoil, parametrized by
  angle-of-attack and inflow velocity. This development is made in parallel for
  two methods: a divergence-conforming reduced basis method, and a conventional
  one (read: not divergence-conforming).
\item Section~\ref{sec:results} presents the most important results from this
  example. Particular emphasis is made on evaluating the convergence properties
  of the reduced methods (that is, the degree to which their solutions agree
  with the corresponding high-fidelity method) and their speed of execution.
\item Finally, Section~\ref{sec:conc} summarizes the findings.
\end{itemize}

\section{Stationary Navier-Stokes problem}

\subsection{Weak formulation}
\label{sec:setting}
We consider the stationary Navier-Stokes equations,
\begin{alignat}{2}
  \label{eqn:ns-1}
  -\nu \Delta \bm u + (\bm u \cdot \nabla) \bm u + \nabla p &= \bm f && \qquad \text{in } \Omega, \\
  \label{eqn:ns-2}
  \nabla \cdot \bm u &= 0 && \qquad \text{in } \Omega, \\
  \label{eqn:ns-3}
  \bm u &= \bm g && \qquad \text{on } \Gamma_\text{D}, |\Gamma_D| > 0 \\
  \label{eqn:ns-4}
  -p \bm n + \nu (\nabla \bm u) \bm n &= \bm h && \qquad \text{on } \Gamma_\text{N}.
\end{alignat}
where $\nu$ is the viscosity, $\bm u, p$ are the unknown velocity and pressure,
$\bm{f},\bm{g},\bm{h}$ correspond to exogenous data,
$\Omega \subset \mathbb R^d$ is the domain of interest with boundary
$\partial \Omega = \Gamma_\text{D} \cup \Gamma_\text{N}$, $\Gamma_\text{D} \cap \Gamma_\text{N} = \emptyset$
and~$\bm{n}$ denotes the external unit normal vector.

We define the function spaces
\begin{align}
  U_{\bm{\zeta}} &= \left\{ \bm u \in [H^1(\Omega)]^d \;|\; \bm u = \bm \zeta \text{ on } \Gamma_\text{D} \right\}, \\
  P &= L^2(\Omega).
\end{align}
and note that the weak Galerkin formulation of the problem is to find
$(\bm u, p) \in (U_{\bm g}, P)$ such that, for all $(\bm w, q) \in (U_0, P)$ it holds
that
\begin{align}
  a(\bm u, \bm w) + c(\bm u, \bm u, \bm w) + b(p, \bm w) &= d(\bm w), \label{eqn:var-1} \\
  b(q, \bm u) &= 0, \label{eqn:var-2}
\end{align}
where the linear, bilinear and trilinear forms $a,b,c,d$ are defined as
\begin{subequations}
\begin{align}
  a(\bm u, \bm w) &= \nu \int_\Omega \nabla \bm u : \nabla \bm w, \label{eqn:form-3} \\
  b(p, \bm w) &= -\int_\Omega p \nabla \cdot \bm w, \label{eqn:form-2} \\
  c(\bm u, \bm v, \bm w) &= \int_\Omega (\bm u \cdot \nabla) \bm v \cdot \bm w. \label{eqn:form-4}\\
  d(\bm w) &= \int_{\Gamma_\text{N}} \bm h \cdot \bm w + \int_{\Omega} \bm f \cdot \bm w, \label{eqn:form-1}
\end{align}
\end{subequations}

For the purposes of solving \eqref{eqn:var-1}--\eqref{eqn:var-2}, it is
customary to introduce a \emph{lift function} $\bm \ell \in U$ satisfying
$\bm{\ell} = \bm{g}$ on $\Gamma_\text{D}$ and to solve for the difference $\bm u - \bm
\ell$, which satifisies homogeneous boundary conditions and thus lives in a
linear space, as opposed to $\bm u$, which resides in an affine space. To this
end, write $\bm u = \bm u_0 + \bm \ell$. The modified problem then reads: find
$(\bm u_0, p) \in (U_0, P)$ such that, for all $(\bm w, q) \in (U_0, P)$ it
holds that
\begin{align}
a(\bm u_0, \bm w) + c(\bm u_0, \bm u_0, \bm w) + c_0(\bm u_0, \bm w) +
b(p, \bm w) &= d_0(\bm w)
\label{eqn:vvar-1}\\
  b(q, \bm u_0) &= d_1(q).
\label{eqn:vvar-2}
\end{align}
where the new forms read
\begin{align}
c_0(\bm u, \bm w)&=
c(\bm u, \bm{\ell}, \bm w)+c(\bm{\ell}, \bm u, \bm w)
\\
d_0(\bm w) &=d(\bm w)-c(\bm{\ell},\bm{\ell},  \bm w)-a(\bm{\ell}, \bm w),
\\
d_1(q)&=-b(q,\bm{\ell})
\end{align}
One may then recover the original through $\bm u = \bm u_0 + \bm \ell$. The lift-based procedure for
treating Dirichlet boundary conditions is well-known, and most FEM implementations handle
boundary conditions automatically by categorizing degrees of freedom as internal
(free) and boundary (constrained), requiring no cognitive overhead for the user.
However, it is worth alluding to this point here since the reduced-basis method
presupposes homogeneous Dirichlet boundary conditions, in order to form a linear space.

\subsection{Velocity-only formulation}

For what follows it is worth considering the effects of choosing divergence-free
velocity function spaces. In particular, we will assume in this subsection that
$\nabla \cdot \bm{u}_0 = 0$ for all $\bm{u}_0 \in U_0$, and also that
$\nabla \cdot \bm \ell = 0$. This ensures that $\nabla \cdot \bm u = 0$, so the
continuity equation~\eqref{eqn:ns-2} and its variational counterpart
\eqref{eqn:vvar-2} are both satisfied by construction. The momentum equation
\eqref{eqn:vvar-1} then reduces to
\begin{equation}
a(\bm u_0, \bm w) + c(\bm u_0, \bm u_0, \bm w) + c_0(\bm u_0, \bm w)
= d_0(\bm w)
\label{eqn:velonly}
\end{equation}
In other words, we have a \emph{velocity-only} formulation.

Considering numerical approximations of~\eqref{eqn:ns-1}--\eqref{eqn:ns-4},
formulation~\eqref{eqn:velonly} has some practical benefits: in view of the
absence of the pressure variable, it voids all concerns about
pressure stability, and it reduces the size of the discretization.
On the other hand, two problems arise. The first problem is that it is
practically infeasible to create divergence-free function spaces, which is
why this formulation does not see much use in conventional approximation methods.
The second problem is that the pressure
is often of interest, e.g.~in post-processing operations to determine flow-induced loads by means
of extraction functions~\cite{Melbo2003goe,Brummelen2012fep}, depreciating a solver
based on~\eqref{eqn:velonly} in many applications.

The proposed solenoidal RB method forms an instance of the velocity-only formulation~\eqref{eqn:velonly}.
However, the problems mentioned above for conventional approximation methods are mitigated in the
RBM setting. First, the inherent solenoidality of the snapshots facilitates the construction of a
divergence-free~\emph{reduced basis}. Each snapshot corresponds to a high-fidelity solution
of~\eqref{eqn:vvar-1}--\eqref{eqn:vvar-2} and, hence, the computational cost of constructing the
solenoidal basis is indeed significant. However, the computational cost of constructing snapshots is
inherent to RBM, and it is ignored because this expense is incurred in the offline stage.  Second,
the pressure solution can be recovered \emph{a posteriori}, by solving \eqref{eqn:vvar-1} with a
different test space, to find $p\in P$ such that for all $\bm s \in S$ it holds that
\begin{equation}
  b(p, \bm s) = d_0(\bm s) - a(\bm u_0, \bm s) - c(\bm u_0, \bm u_0, \bm s) - c_0(\bm u_0, \bm s),
  \label{eqn:sup-vvar}
\end{equation}
where $\bm u_0$ is the RB velocity solution. Here, the test space $S$ should
be chosen so that the left-hand-side becomes suitably non-singular.
Section~\ref{sec:stab} elaborates on a procedure to this purpose.

\subsection{Parametric dependence}
\label{sec:pardep}

We consider the case where the problem \eqref{eqn:var-1}--\eqref{eqn:var-2}
depends on a number of \emph{parameters}. Denote by $\mathcal{P}$ the parameter
space (which will for the moment remain abstract) and by $\bm \mu$ any given
element of $\mathcal{P}$. The typical parameters of interest are generally
grouped into classes: \emph{physical parameters}, \emph{data parameters}
and \emph{geometric parameters}. The purpose of the following
discussion is to explicitly resolve the parametric dependence of the multi-linear forms.
We shall denote this with notation such as $a(\bm u, \bm w; \bm\mu)$, and similarly for other forms.

Loosely stated, physical parameters are those parameters impacting the physical
quantities involved in the model, e.g.~the viscosity $\nu$ in
\eqref{eqn:ns-1}. Data parameters impact the source term~$\bm{f}$, and Dirichlet or Neumann
boundary data, i.e.~the functions $\bm g$ and $\bm h$ in
\eqref{eqn:ns-3}--\eqref{eqn:ns-4}, and through $\bm g$ also the lifting
function $\bm \ell$. The aforementioned types of parameters appear directly in
the multi-linear forms, viz.
\begin{align}
  a(\bm u, \bm w; \bm \mu) &= \nu(\bm \mu) \int_\Omega \nabla \bm u : \nabla \bm w, \\
  d_0(\bm w; \bm \mu) &= \int_{\Gamma_\text{N}} \bm h(\bm \mu) \cdot \bm w
                        + \int_{\Omega} \bm{f}(\bm{\mu}) \cdot \bm w, \\
  c_0(\bm u, \bm w; \bm \mu) &= c(\bm u, \bm \ell(\bm \mu), \bm w) + c(\bm \ell(\bm \mu), \bm u, \bm w), \\
  d_1(q; \bm \mu) &= - b(q, \bm \ell(\bm \mu)).
\end{align}
One may note that while the original variational form
\eqref{eqn:var-1}--\eqref{eqn:var-2} is set in a function space that is
possibly parameter-dependent via the Dirichlet data, the corresponding homogeneous variational form
\eqref{eqn:vvar-1}--\eqref{eqn:vvar-2} is in canonical form in the sense that the
ambient spaces are parameter independent, and the parameter dependence occurs only in
the functionals.

Let us next consider geometric parameters, by which we mean that
the computational domain $\Omega = \Omega(\bm \mu)$ varies in relation to the parameters.
To cast this problem in canonical form, i.e. to transfer the parameter dependence
from the ambient spaces to the functionals in the weak formulation, we pull back
the corresponding multi-linear forms to (function spaces on) a fixed \emph{reference domain}.
Let $\hat{\Omega}$ denote a suitable reference domain, (usually
$\hat{\Omega} = \Omega({\bm \mu})$ for some $\bm \mu$, but this is not
necessarily the case), such that there exists a bi-Lipschitz bijection
$\chi_{\bm \mu} : \Omega(\bm \mu) \to\hat{\Omega}$,
which maps any of the parametrized domains to the reference domain. We insist that
the images of $\Gamma_\text{D}$ and $\Gamma_\text{N}$ are fixed in $\hat{\Omega}$.
In other words, it holds that
\begin{equation}
  \chi_{\bm \mu}( \Gamma_\text{D}(\bm \mu) ) =
  \hat{\Gamma}_\text{D} \subset \partial \hat{\Omega}
\end{equation}
and that $\hat{\Gamma}_\text{D}$ as defined is independent of $\bm
\mu$. Based on the domain map $\chi_{\bm \mu}$, we introduce
two homeomorphisms ${\pi}_{\bm\mu}^\textsc{v}:[H^1(\hat{\Omega})]^d\to{}[H^1(\Omega(\bm \mu))]^d$
and ${\pi}_{\bm\mu}^\textsc{p}:L^2(\hat{\Omega})\to{}L^2(\Omega(\bm{\mu}))$. For instance,
$\bm{\pi}_{\bm\mu}:=({\pi}_{\bm\mu}^\textsc{v},{\pi}_{\bm\mu}^\textsc{p})$ can be defined as the canonical
transformation:
\begin{equation}
  \label{eqn:cantrf}
\bm{\pi}_{\bm\mu}:(\hat{\bm{u}},\hat{p})\mapsto\big(\hat{\bm{u}}\circ\chi_{\bm\mu},\hat{p}\circ\chi_{\bm\mu}\big)
\end{equation}
However, alternate maps can be constructed. The pullback of the multi-linear forms
$a,b,c,c_0$ (to the reference domain, by the map~$\bm{\pi}_{\bm\mu}$) is then defined by:
\begin{equation}
  \label{eqn:newform}
\begin{aligned}
({\bm\pi}^*_{\bm\mu}a)\big(\hat{\bm{u}},\hat{\bm{w}};\bm{\mu}\big)
&=
a\big({\pi}_{\bm\mu}^\textsc{v}\hat{\bm{u}},{\pi}_{\bm\mu}^\textsc{v}\hat{\bm{w}};\bm{\mu}\big)
\\
({\bm\pi}^*_{\bm\mu}b)\big(\hat{p},\hat{\bm{w}};\bm{\mu}\big)
&=
b\big({\pi}_{\bm\mu}^\textsc{p}\hat{p},{\pi}_{\bm\mu}^\textsc{v}\hat{\bm{w}};\bm{\mu}\big)
\\
({\bm\pi}^*_{\bm\mu}c)\big(\hat{\bm{u}},\hat{\bm{v}},\hat{\bm{w}};\bm{\mu}\big)
&=
c\big({\pi}_{\bm\mu}^\textsc{v}\hat{\bm{u}},{\pi}_{\bm\mu}^\textsc{v}\hat{\bm{v}},{\pi}_{\bm\mu}^\textsc{v}\hat{\bm{w}};\bm{\mu}\big)
\\
({\bm\pi}^*_{\bm\mu}c_0)\big(\hat{\bm{u}},\hat{\bm{w}};\bm{\mu}\big)
&=
c\big({\pi}_{\bm\mu}^\textsc{v}\hat{\bm{u}},{\pi}_{\bm\mu}^\textsc{v}\hat{\bm{\ell}}(\bm{\mu}),{\pi}_{\bm\mu}^\textsc{v}\hat{\bm{w}};\bm{\mu}\big)
+
c\big({\pi}_{\bm\mu}^\textsc{v}\hat{\bm{\ell}}(\bm{\mu}),{\pi}_{\bm\mu}^\textsc{v}\hat{\bm{u}},{\pi}_{\bm\mu}^\textsc{v}\hat{\bm{w}};\bm{\mu}\big)
\end{aligned}
\end{equation}
and the pullback of the linear forms $d_0$ and $d_1$ is defined analogously. By means of the
pullback operation ${\bm\pi}^*_{\bm\mu}$, the parametric version
of~\eqref{eqn:vvar-1}\nobreakdash--\eqref{eqn:vvar-2} can be cast
in the canonical form: find $(\hat{\bm u}_0, \hat{p}) \in (\hat{U}_0, \hat{P})$ such
that, for all $(\hat{\bm w}, \hat{q}) \in (\hat{U}_0, \hat{P})$ it holds
that
\begin{align}
  \nonumber (\bm\pi_{\bm\mu}^* a)(
    \hat{\bm u}_0,
    \hat{\bm w};
    \bm \mu
  ) + (\bm\pi_{\bm\mu}^* c)(
    \hat{\bm u}_0,
    \hat{\bm u}_0,
    \hat{\bm w};
    \bm \mu
  ) + (\bm\pi_{\bm\mu}^* c_0)(
    \hat{\bm u}_0,
    \hat{\bm w};
    \bm \mu
  )
  \\
  + \; (\bm\pi_{\bm\mu}^* b)(
    \hat{p},
    \hat{\bm w};
    \bm \mu
  ) &= (\bm\pi_{\bm\mu}^* d_0)(
    \hat{\bm w};
    \bm \mu
  )
  \label{eqn:war-1}, \\
  (\bm\pi_{\bm\mu}^* b)(
    \hat{q},
    \hat{\bm u}_0;
    \bm \mu
  ) &= (\bm\pi_{\bm\mu}^* d_1)(
    \hat{q};
    \bm \mu
  ). \label{eqn:war-2}
\end{align}
One may note that~\eqref{eqn:war-2} is precisely \eqref{eqn:vvar-1}--\eqref{eqn:vvar-2} expressed
with pullbacks through $\bm \pi_{\bm \mu}^*$ and parameter-dependent forms.

\begin{remark}
  \label{rem:lift}
  It is to be noted that the lift $\hat{\bm \ell}$ in~\eqref{eqn:newform} is constructed on
  the reference domain, and then transported
  via $\pi_{\bm \mu}^\textsc{v}$ to the parametrized domain. In this construction it is important
  that $\smash[t]{\pi_{\bm \mu}^\textsc{v} \hat{\bm \ell}}$ and~${\bm \ell}$ are compatible, so that
  the trace of $\bm{\ell}=\smash[t]{\pi_{\bm \mu}^\textsc{v}\hat{\bm \ell}}$ on $\Gamma_\text{D}$
  coincides with the prescribed Dirichlet data.
\end{remark}

\subsection{Affine representation}
\label{sec:affrep}

We shall assume, and this is critical to the success of what follows, that the
forms in \eqref{eqn:war-1}--\eqref{eqn:war-2} may be expressed as linear
combinations of forms that are parameter-\emph{in}dependent.
\begin{align}
  (\bm\pi_{\bm\mu}^* a)(
    \hat{\bm u},
    \hat{\bm w};
    \bm \mu
  ) &= \sum_{i=1}^{N_a} \theta^a_i(\bm \mu) \hat{a}_i(\hat{\bm u}, \hat{\bm w}), \label{eqn:split-a} \\
    (\bm\pi_{\bm\mu}^* b)(
    \hat{p},
     \hat{\bm w};
    \bm \mu
  ) &= \sum_{i=1}^{N_b} \theta^b_i(\bm \mu) \hat{b}_i(\hat{p}, \hat{\bm w}), \label{eqn:split-b} \\
  (\bm\pi_{\bm\mu}^* c)(
    \hat{\bm u},
    \hat{\bm v},
    \hat{\bm w};
    \bm \mu
  ) &= \sum_{i=1}^{N_c} \theta^c_i(\bm \mu)
      \hat{c}_i(\hat{\bm u}, \hat{\bm v}, \hat{\bm w}), \label{eqn:split-c} \\
  (\bm\pi_{\bm\mu}^* c_0)(
    \hat{\bm u},
    \hat{\bm w};
    \bm \mu
  ) &= \sum_{i=1}^{N_c} \theta^{c_0}_i(\bm \mu)
      \hat{c}{}^{0}_i(\hat{\bm u}, \hat{\bm w}), \label{eqn:split-c0} \\
  (\bm\pi_{\bm\mu}^* d_0)(
    \hat{\bm w};
    \bm \mu
  ) &= \sum_{i=1}^{N_1} \theta^{d_0}_i(\bm \mu) \hat{d}^0_i(\hat{\bm w}), \label{eqn:split-d0} \\
  (\bm\pi_{\bm\mu}^* d_1)(
    \hat{q};
    \bm \mu
  ) &= \sum_{i=1}^{N_2} \theta^{d_1}_i(\bm \mu) \hat{d}^1_i(\hat{q}) \label{eqn:split-d1}
\end{align}
Most boundary-condition parameters naturally produce such representations so
long as the conditions themselves adhere to the same form, i.e.~if
$\bm g = \sum_i \xi_i(\bm \mu) \bm g_i$ then we can also write
$\hat{\bm \ell} = \sum_i \xi_i(\bm \mu) \hat{\bm \ell_i}$ where
$\pi_{\bm\mu}^\textsc{v} \bm \hat{\bm{\ell}}_i = \bm{g}_i$ on
$\Gamma_\text{D}$. This expression is substituted into
\eqref{eqn:newform}, and the resulting affine
representation follows directly. The same derivation can be done for Neumann
data $\bm h$ and physical parameters such as $\nu$.

Geometric parameters can lead to significantly more complicated representations
unless they are of the trivial sort. In this case one must map the integrals
\eqref{eqn:form-1}--\eqref{eqn:form-4} to the reference configuration, and the
resulting expressions must be brought into forms compatible with
\eqref{eqn:split-a}--\eqref{eqn:split-d1}. Since these expressions usually
involve the Jacobian and its inverse, it is desirable that they are as simple as
possible.

\section{Model order reduction}
\label{sec:reduction}

\subsection{High-fidelity approximation}
The reduced basis is built from underlying high-fidelity finite-element or isogeometric~\cite{Cottrell2009iat}
approximations of~\eqref{eqn:war-1}--\eqref{eqn:war-2}. To provide a setting for the RB, we therefore
first introduce the high-fidelity setting.
We introduce a partition $\mathcal{T}^h$ of the reference domain~$\hat{\Omega}$ into non-overlapping
element domains. The superscript $h$ indicates dependence
on a mesh-resolution parameter. We insist that $\mathcal{T}^h$ satisfies the usual uniformity and
regularity conditions. The mesh $\mathcal{T}^h$ serves as a support structure for conforming
finite-dimensional approximation spaces $\hat{U}_0^h\subset\hat{U}_0$ and $\hat{P}^h\subset\hat{P}$.
We assume that the pairing of the velocity and pressure spaces is stable in the inf-sup sense; see,
for instance, \cite{Taylor1973nsn,Nedelec1986nfm,Raviart1981mfe,Johannessen2015dcd}
and~\cite{Evans2013idc1,Evans2013idc3,Buffa2011ias} for stable velocity-pressure pairs in
the finite-element and isogeometric setting,
respectively. For each $\bm{\mu}\in\mathcal{P}$, the high-fidelity approximation
$(\hat{\bm{u}}_0^h(\bm{\mu}),\hat{p}^h(\bm{\mu}))$ corresponds to the solution
of~\eqref{eqn:war-1}\nobreakdash--\eqref{eqn:war-2} with $\hat{U}_0,\hat{P}$ replaced by their
discrete counterparts $\hat{U}_0^h,\hat{P}^h$.

The considered high-fidelity approximation spaces $\hat{U}_0^h,\hat{P}^h$ must be sufficiently fine to
ensure that the error in the high-fidelity approximations is acceptable in accordance with a prescribed tolerance
for all considered parameter values $\bm{\mu}\in\mathcal{P}$. For instance, one may insist that
\begin{equation}
\label{eqn:hifi_con}
\sup_{\bm{\mu}\in\mathcal{P}}
\big\|\hat{u}_0(\bm{\mu})-\hat{u}_0^h(\bm{\mu})\big\|_{H^1(\hat{\Omega})}\leq\textsc{tol}_0,
\qquad
\sup_{\bm{\mu}\in\mathcal{P}}
\big\|\hat{p}(\bm{\mu})-\hat{p}^h(\bm{\mu})\big\|_{L^2(\hat{\Omega})}\leq\textsc{tol}_1,
\end{equation}
for certain prescribed tolerances $\textsc{tol}_{0,1}$.
The requirements imposed by~\eqref{eqn:hifi_con} generally lead to very high dimensionality of $\hat{U}_0^h\times\hat{P}^h$,
impeding direct use of the high-fidelity model for fast- or many-query applications.

\subsection{Construction of the Reduced Basis by Proper Orthogonal Decomposition}
The main notion underlying RB methods is that in many applications, the dimension of the
high-fidelity discrete solution space  far exceeds the natural dimension of the parametrized
model itself. In other words, it is anticipated that the high-fidelity solution space
\begin{equation}
  \left\{ (\hat{\bm u}^h_0(\bm \mu), \hat{p}^h(\bm \mu)) \;|\; \bm \mu \in \mathcal{P} \right\}
  \subset \hat{U}_0^h \times \hat{P}^h
  \label{eqn:nspace}
\end{equation}
admits an adequately accurate approximation by a finite-dimensional space whose dimension is
significantly lower than that of the high-fidelity approximation space $\smash[t]{\hat{U}_0^h \times \hat{P}^h}$.
This dimension may be
considerably more manageable than the dimension of any discretization of $\smash[t]{\hat{U}_0 \times \hat{P}}$
necessary to achieve sufficient accuracy with conventional methods.

We use standard Proper Orthogonal Decomposition (POD)~\cite{Chatterjee2000ipo,Quarteroni2016rbm} to produce
a low-dimensional approximation of the space~\eqref{eqn:nspace}. POD proceeds by sampling $\mathcal{P}$
and generating an \emph{ensemble} of high-fidelity approximations
\begin{equation}
\label{eqn:ensemble}
  \Phi = \left\{
    \varphi_i = (\hat{\bm u}^h_0(\bm \mu_i), \hat{p}^h(\bm \mu_i))
  \right\}_{i=1}^N
\end{equation}
with conventional finite element or isogeometric methods. A covariance matrix is then
constructed,
\begin{equation}
  \label{eqn:covmx}
  C_{ij} = z(\varphi_i, \varphi_j)
\end{equation}
where the bilinear symmetric positive definite covariance function
$z(\cdot,\cdot)$ is left unspecified for the time being. Given eigenpairs
$(\lambda_j, \bm v^j)$ ($\lambda_1\geq\lambda_2\geq\cdots\geq\lambda_N$)
of $\bm C$, one may then compose a suitably small number $M$
of basis functions~\cite[(6.10)]{Quarteroni2016rbm}:
\begin{equation}
  \label{eqn:spd}
  \psi_j = \frac{1}{\sqrt{\lambda_j}} \sum_{i=1}^N v^j_i \varphi_i,
  \qquad 1 \leq j \leq M.
\end{equation}
The basis $V=\{\psi_j\}_{j=1}^M$ constitutes a reduced basis in the sense
that $\operatorname{span}V\subset\operatorname{span}\Phi$. The basis~$V$ is orthonormal
in the $z(\cdot,\cdot)$ inner product by construction and it
is optimal in the following sense: Let
$W =\{ w_j\}_{j=1}^M$ be any orthonormal basis of size $M$, and
$\bm P_W:\Phi\to{}W$ be the orthogonal projection onto $W$, i.e.
\begin{equation}
  \bm P_W \varphi = \sum_{j=1}^M z(\varphi, w_j) w_j
\end{equation}
Then it holds (see \cite[Proposition 6.2]{Quarteroni2016rbm})
\begin{equation}
  \sum_{i=1}^N \| \varphi_i - \bm P_W \varphi_i \|_z^2 \geq
  \sum_{i=1}^N \| \varphi_i - \bm P_V \varphi_i \|_z^2 =
  \sum_{k=M+1}^N \lambda_k.
\end{equation}
Indeed, since we have by construction that
\begin{equation}
  \sum_{i=1}^N \| \varphi_i \|_z^2 = \sum_{k=1}^N \lambda_k,
\end{equation}
we obtain the relationship
\begin{equation}
  \frac{ \sum_{i=1}^N \| \varphi_i - \bm P_V \varphi_i \|_z^2 }{ \sum_{i=1}^N \| \varphi_i \|_z^2 }
  = \frac{ \sum_{k=M+1}^N \lambda_k }{ \sum_{k=1}^N \lambda_k }.
\end{equation}
From this it follows that if $M = M(\epsilon)$ is chosen according to the inequality
\begin{equation}
  \label{eqn:error}
  \sum_{j=1}^M \lambda_j
  \geq \left(1 - \epsilon^2\right) \sum_{j=1}^N \lambda_j,
\end{equation}
then the reduced basis $V$ satisfies the error bound
\begin{equation}
\label{eqn:PVbound}
\left( \sum_{i=1}^N \big\| \varphi_i - \bm P_V \varphi_i \big\|_z^2 \right)^{1/2}
  \leq \epsilon \left( \sum_{i=1}^N \big\| \varphi_i \big\|_z^2 \right)^{1/2}.
\end{equation}
If the ensemble $\Phi$ is representative of the complete high-fidelity solution space~\eqref{eqn:nspace}
in the sense that the latter is contained in $\operatorname{span}\Phi$, then~\eqref{eqn:PVbound}
implies the following bound for the projection error:
\begin{equation}
\label{eqn:RBerror}
\frac{\left( \int_\mathcal{P}\big\|\varphi(\bm \mu) - \bm P_V \varphi(\bm \mu)\big\|_z^2\dd\bm \mu \right)^{1/2}}
{\left( \int_\mathcal{P} \big\| \varphi(\bm \mu) \big\|_z^2 \dd\bm \mu \right)^{1/2}}
  \leq \epsilon
  .
\end{equation}

\subsection{Field separation}
\label{sec:fieldsep}
In practice, the extraction of a reduced basis from~\eqref{eqn:ensemble} is performed separately
on the velocity and pressure components. In other words, first one selects a velocity covariance
function $z=z_\textsc{v}$. Typically,
\begin{equation}
\label{eqn:z_u}
z_\textsc{v}\big((\hat{\bm u}_0, \hat{p}), (\hat{\bm v}_0, \hat{q})\big) =
  \left( \hat{\bm u}_0, \hat{\bm v}_0 \right)_{H^1(\hat{\Omega})}
  = \int_{\hat{\Omega}} \nabla \hat{\bm u}_0 : \nabla \hat{\bm v}_0.
\end{equation}
This is equivalent to using the $H^1$-seminorm for the velocity space. By virtue of the standing
assumption $|\hat{\Gamma}_\text{D}| > 0$, this covariance function is
symmetric positive definite on~$U_0$. On the product space $U_0 \times P$, it is only semidefinite.
The covariance function~\eqref{eqn:z_u} yields the velocity covariance matrix:
\begin{equation}
  C^{\textsc{v}}_{ij} = \left( \hat{\bm u}^h_0(\bm{\mu}_i), \hat{\bm u}^h_0(\bm{\mu}_j) \right)_{H^1(\hat{\Omega})}
\end{equation}
whose eigenpairs $(\lambda_j, \bm v^{j})$ give the reduced velocity
basis functions
\begin{equation}
\label{eqn:uRj}
  \hat{\bm u}^\textsc{r}_j = \frac{1}{\sqrt{\lambda_j}} \sum_{i=1}^N v^{j}_i \hat{\bm u}^h_0(\bm{\mu}_i),
  \qquad 1 \leq j \leq M_\textsc{v}
\end{equation}
where $M_\textsc{v}$ is the chosen dimension of the reduced velocity space. Even though
$z_\textsc{v}$ is only semidefinite, since the basis is restricted to only those
components on which $z_\textsc{v}$ \emph{is} definite, the basis functions remain
linearly independent.

Independent of the construction of the reduced velocity bases, a pressure basis can
be generated by exactly the same method using, for example, the covariance function
$z=z_\textsc{p}$ according to:
\begin{equation}
z_\textsc{p}\big((\hat{\bm u}_0, \hat{p}), (\hat{\bm v}_0, \hat{q})\big) =
  \left( \hat{p}, \hat{q} \right)_{L^2(\hat{\Omega})}.
\end{equation}
providing a reduced pressure basis:
\begin{equation}
  \label{eqn:pRj}
  \hat{p}^\textsc{r}_j = \frac{1}{\sqrt{\lambda_j}} \sum_{i=1}^N v^{j}_i \hat{p}^h(\bm{\mu}_i),
  \qquad 1 \leq j \leq M_\textsc{p}
\end{equation}
where $M_\textsc{p}$ is the specified dimension of the reduced velocity space, and $(\lambda_j,\bm{v}^j)$ now
refers to the eigenpairs of the pressure covariance matrix.

\subsection{Pressure stabilization}
\label{sec:stab}
The reduced velocity and pressure spaces, $U_0^{\textsc{r}},P^{\textsc{r}}$, do not generally inherit the
inf-sup stability of the high-fidelity approximation spaces. Given the reduced spaces, the solution of
the reduced model generally involves solving a sequence of linear-algebraic problems of the block form
\begin{equation}
  \label{eqn:sys-instab}
  \begin{pmatrix} \bm A_\textsc{vv} & \bm B_\textsc{vp} \\ \bm B_\textsc{pv} & 0\end{pmatrix}
  \begin{pmatrix} \bm u \\ \bm p \end{pmatrix}
  =
  \begin{pmatrix} \bm f_\textsc{v} \\ \bm f_\textsc{p} \end{pmatrix}
\end{equation}
with $(\bm B_\textsc{vp})_{ij}=\bm{\pi}^*_{\bm{\mu}}b(\hat{p}_j^\textsc{r},\hat{\bm{u}}_i^\textsc{r};\bm{\mu})$.
On account of the loss of inf-sup stability of the reduced velocity-pressure pair,
one typically observes that the matrix $\bm B_\textsc{vp}$ is row-rank-deficient, i.e.~it
has a nontrivial kernel, therefore leading to unstable pressure solutions. One procedure to resolve this
instability is to artificially enrich the velocity space with so-called
\emph{supremizers}, viz.~basis functions whose purpose is not to achieve greater approximative power but
rather to keep the system well-posed; see \cite{Ballarin2015ssp}. We will use this technique
both for comparative purposes and for generating pressure solutions for velocity-only
discretizations, see Section \ref{sec:conforming}.

For each pressure snapshot $\hat{p}^h(\bm{\mu}_i)$, the corresponding supremizer is the
velocity function in $\hat{U}^h_0$ that realizes the supremum in the
LBB\footnote{Ladyzhenskaya-Babu\v{s}ka-Brezzi, also known as the inf-sup condition.} condition,
\begin{equation}
  \label{eqn:lbb}
  \inf_{\hat{p}^h \in \hat{P}^h} \sup_{\hat{\bm u}_0^h \in \hat{U}_0^h}
  \frac{\bm{\pi}^*_{\bm{\mu}}b(\hat{p}^h, \hat{\bm u}^h_0; \bm \mu)}
  {\|\hat{p}^h\|_{L^2(\hat{\Omega})} \|\hat{\bm{u}}_0^h\|_{H^1(\hat{\Omega})}}
  \geq \beta_h > 0;
\end{equation}
see \cite{Ballarin2015ssp} for the theoretical background.
%
The supremizers can in fact be conceived of as Riesz representations of the pressure solutions as
viewed through the $\bm{\pi}^*_{\bm{\mu}}b(\cdot,\cdot\,;\bm \mu)$-form. In other words,
given a reduced pressure function $\hat{p}^\textsc{r}$, its corresponding supremizer $\smash[t]{\hat{\bm s}^h_0}$ is
the Riesz representation in $\hat{U}_0^h$ of the linear functional
$\smash[t]{\bm{\pi}^*_{\bm{\mu}}b(\hat{p}^h,\cdot\,; \bm \mu)}$, i.e.~the solution of
\begin{equation}
  \label{eqn:riesz}
\hat{\bm s}^h_0\in  \hat{U}^h_0:
\qquad
  \left( \hat{\bm w}^h_0, \hat{\bm s}^h_0 \right)_{H^1(\hat{\Omega})} = \bm{\pi}^*_{\bm{\mu}}b(\hat{p}^h, \hat{\bm w}^h_0; \bm \mu)
  \qquad\forall\hat{\bm w}^h_0 \in \hat{U}^h_0.
\end{equation}

From~\eqref{eqn:riesz} it follows that the supremizers are formally $\bm{\mu}$\nobreakdash-dependent functions,
via the dependence of the right member of~\eqref{eqn:riesz} on~$\bm{\mu}$, although this dependence
is generally restricted to geometric parameters. They are therefore not directly suitable for use
as basis functions in a reduced-basis method. To bypass the $\bm{\mu}$\nobreakdash-dependence of the supremizers,
we apply the process described as ``approximate supremizer enrichment'' in
\cite[4.2.3]{Ballarin2015ssp}. In this method, supremizer snapshots $\hat{\bm s}^h_0(\bm\mu_i)$ are constructed
offline in tandem with the pressure snapshots $p(\bm \mu_i)$. These supremizer snapshots are then
reduced via POD in exactly the same manner as the velocity snapshots, into a separate and independent reduced space,
which is subsequently added to the reduced velocity space via a direct-sum composition. Formally,
the sum space should be orthogonalized to ensure linear independence. Following~\cite{Ballarin2015ssp}, we have ignored
this on the presumption that the reduced velocity space and the supremizer space are sufficiently disparate to
retain linear independence. This hypothesis is verified in the numerical experiments in Section~\ref{sec:stability}.

Introduction of the supremizer modes yields the following extension of the linear-algebraic
systems~\eqref{eqn:sys-instab}:
\begin{equation}
  \label{eqn:block}
  \begin{pmatrix}
    \bm A_\textsc{vv} & \bm A_\textsc{vs} & \bm B_\textsc{vp} \\
    \bm A_\textsc{sv} & \bm A_\textsc{ss} & \bm B_\textsc{sp} \\
    \bm B_\textsc{pv} & \bm B_\textsc{ps} & 0 \end{pmatrix}
  \begin{pmatrix} \bm u \\ \bm s \\ \bm p \end{pmatrix}
  =
  \begin{pmatrix} \bm f_\textsc{v} \\ \bm f_\textsc{s} \\ \bm f_\textsc{p} \end{pmatrix}
\end{equation}
One may note that with \emph{exact} supremizer enrichment, the practical
consequence of~\eqref{eqn:riesz} is that $\bm B_\textsc{sp}$ coincides
with the Gramian matrix $\smash[t]{(\bm B_\textsc{sp})_{ij}=( \hat{\bm s}^h_i, \hat{\bm s}^h_j )_{H^1(\hat{\Omega})}}$,
which is nonsingular by design as discussed in Section~\ref{sec:reduction}.

The computation of supremizers introduces an additional offline cost which in
our experience is quite negligible. The solution of \eqref{eqn:riesz} for every
parameter sample is cheap next to the more challenging Navier-Stokes problem. A
cost is also paid in the online stage, as the system \eqref{eqn:block} is larger
than \eqref{eqn:sys-instab} by several degrees of freedom that do not
meaningfully contribute to better approximations, as we will see shortly.

It must be noted that the supremizer method is not the only stabilization method
proposed in the literature, see for example
\cite{Weller2008nml,Bergmann2009erp}, where reduced basis methods are stabilized
based on fine-scale Variational Multiscale Methods (VMS).

\subsection{Assembly of the reduced system}
Let $\bm V \in \mathbb{R}^{N \times M}$ be a $\bm \mu$-independent
tall \emph{transformation matrix}, mapping coefficient vectors in the
reduced basis to coefficient vectors in the original high-fidelity
basis, i.e.~the columns of $\bm V$ are the high-fidelity coefficient
vectors of the reduced basis functions.

We aim then to express the solution to a generic finite element system
$\bm A(\bm \mu) \bm u(\bm \mu) = \bm f(\bm \mu)$ by the reduced
coefficient vector $\bm u(\bm \mu) = \bm V \bm u^\textsc{r}(\bm \mu)$, i.e.
\begin{equation}
  \label{eqn:redsys}
  \bm A(\bm \mu) \bm V \bm u^\textsc{r}(\bm \mu) = \bm f(\bm \mu).
\end{equation}
Since \eqref{eqn:redsys} is overdetermined, we choose the Galerkin
approach of enforcing that the residual must be orthogonal to the
reduced space, i.e.
\begin{equation}
  \label{eqn:redorth}
  \bm V^\intercal \left(
  \bm A(\bm \mu) \bm V \bm u^\textsc{r}(\bm \mu) - \bm f(\bm \mu)
  \right) = \bm 0,
\end{equation}
or, more traditionally,
\begin{equation}
  \label{eqn:redorth-final}
  \bm V^\intercal \bm A(\bm \mu) \bm V \bm u^\textsc{r}(\bm \mu) =
  \bm V^\intercal \bm f(\bm \mu).
\end{equation}

The importance of \eqref{eqn:split-a}--\eqref{eqn:split-d1} then
becomes clear.  Although the solution of the reduced problem involves
nonlinear systems of size $M$, which are quick to solve when $M$ is
small (even if they are dense, as they will be), the assembly of these
systems remains nontrivial. However, given \eqref{eqn:redorth-final},
each matrix and vector corresponding to a bilinear or linear form in
\eqref{eqn:split-a}--\eqref{eqn:split-d1} may be reduced to the new
basis independently without regard to parametric dependence. It is
then vital that $\bm V$ is indeed parameter-independent.

For example, let us consider \eqref{eqn:split-a}. Let the matrix
$\bm M^\textsc{h}_i$ be the high-fidelity matrix corresponding to the bilinear
form $\hat{a}_i$. (Here, superscript H refers to \emph{high-fidelity} and
superscript R refers to the reduced model.) Then we have
\begin{equation}
  \bm M^\textsc{r}(\bm \mu)
  = \bm V^\intercal \bm M^\textsc{h}(\bm \mu) \bm V
  = \bm V^\intercal \left( \sum_i \theta^d_i(\bm \mu) \bm M^\textsc{h}_i \right) \bm V
  = \sum_i \theta^d_i(\bm \mu) \underbrace{\bm V^\intercal \bm M^\textsc{h}_i \bm V}_{\bm M^\textsc{r}_i},
\end{equation}
where each of the small matrices $\bm M^\textsc{r}_i$ can be computed in the
offline stage. Assembling the reduced model system matrices then becomes
trivial.

\subsection{Assembly of nonlinear terms}

Typically, the presence of a nonlinear term such as $(\bm\pi_{\bm\mu}^* c)(\hat{\bm u}_0,\hat{\bm u}_0,\hat{\bm w};\bm \mu)$ in~\eqref{eqn:war-1} is treated by introducing a nonlinear solver (e.g.~Newton's method or Picard fixed-point
iteration) that is based on solving a sequence of \emph{linear} systems. These linear systems generally involve
(an approximation to) the Fr\'echet derivative of the nonlinear term, evaluated at a sequence of approximations to
the solution of the nonlinear problem. Because the sequence of approximations is evidently not known in advance,
the assembly of the Jacobian matrices corresponding to the Fr\'echet derivatives is generally done by a standard assembly
procedure that includes an integration loop. However, the computational cost of such an assembly operation is prohibitive
in the online stage of an RBM. Currently, there is no standard methodology for the efficient treatment of general
nonlinear terms in the online stage of an RBM.

The nonlinear term in the Navier--Stokes equations has a particular structure that can be exploited to obtain an efficient assembly in the online stage. Specifically, the nonlinear term in~\eqref{eqn:war-1} is quadratically nonlinear, as indicated by the fact that it can be represented as a trilinear form with a repeated argument. The
quadrature-based assembly of the Jacobian matrix in the online stage can then be avoided by assembling in the
offline stage the RB third-order tensors corresponding to the components $\hat{c}_i$ in the affine representation~\eqref{eqn:split-c} of the trilinear form.
In particular, one can construct the third-order tensors $\smash[t]{\{{\bm C}^i\}_{i=1}^{N_c}}$ with components:
\begin{equation}
  C^i_{jkl} = \hat{c}{}_i\big(\hat{\bm u}^\textsc{r}_j, \hat{\bm u}^\textsc{r}_k, \hat{\bm u}^\textsc{r}_l\big)
  \qquad
  j,k,l\in\{1,2,\ldots,M_{\textsc{v}}\}
  \label{eqn:tensor}
\end{equation}
For any approximation $\underline{\hat{\bm u}}{}_0\in{}U_0^{\textsc{r}}$ the Jacobian matrix associated with the
Fr\'echet derivative of $\bm\pi_{\bm\mu}^* c$ restricted to the reduced velocity space~$U_0^{\textsc{r}}$ can then be assembled in accordance with the ultimate expression in the chain of identities:
\begin{multline}
\frac{d}{ds}(\bm\pi_{\bm\mu}^* c)(\underline{\hat{\bm u}}{}_0+s\hat{\bm u}^\textsc{r}_j,\underline{\hat{\bm u}}{}_0+s\hat{\bm u}^\textsc{r}_j,\hat{\bm u}^\textsc{r}_k;\bm \mu)\bigg|_{s=0}
=
(\bm\pi_{\bm\mu}^* c)(\underline{\hat{\bm u}}{}_0,\hat{\bm u}^\textsc{r}_j,\hat{\bm u}^\textsc{r}_k;\bm \mu)
+
(\bm\pi_{\bm\mu}^* c)(\hat{\bm u}^\textsc{r}_j,\underline{\hat{\bm u}}{}_0,\hat{\bm u}^\textsc{r}_k;\bm \mu)
\\
=
\sum_{i=1}^{N_c}\theta^c_i(\bm\mu)
\sum_{l=1}^{M_{\textsc{v}}}
\big(C^i_{ljk}+C^i_{jlk}\big)\eta_l
\end{multline}
where $\eta_l$ denote the coefficients of $\underline{\hat{\bm u}}{}_0$ relative to the reduced basis, i.e.
$\underline{\hat{\bm u}}{}_0=\sum_{l=1}^{M_{\textsc{v}}}\eta_l\hat{\bm u}^\textsc{r}_l$. Hence, for any
approximation $\underline{\hat{\bm u}}{}_0\in{}U^{\textsc{r}}_0$, the Jacobian matrix can be assembled on the
basis of the pre-computed third-order tensors. For large models, the storage requirements for the third-order
tensor are unmanageable. However, the storage of such a tensor for a reduced space {\em is\/} feasible.
Although these tensors are dense, we find that the modest sizes keep the storage requirements realistic,
typically in the regions of hundreds of megabytes for ${\sim}100$ basis
functions.

\section{Divergence-Conforming methods}
\label{sec:divconf}

\subsection{Separated velocity and pressure bases}
\label{sec:conforming}

In high-fidelity terms, a divergence-conforming method is any method which
ensures that the pressure space is exactly equal to the divergence of the
velocity space. Today this sort of method is often seen in the field of
isogeometric analysis, where the use of function spaces based on B-splines,
together with the Piola transformation, creates sufficient flexibility to easily
work with such spaces.

The practical consequence of such a method is that it eliminates all forms of
pressure instabilities and simultaneously guarantees both a stable method for
convection-dominated flow problems and a velocity solution that is pointwise
solenoidal (divergence-free). The latter point is interesting, since it paves
the way to generate fully divergence-free reduced velocity spaces.

To ensure that solenoidal velocity fields remain solenoidal even in the case of
parameter-dependent geometries, one employs the Piola transform, whereby the
physical manifestation of the reference velocity field $\hat{\bm u}$ is not the
canonical transformation \eqref{eqn:cantrf} (which, indeed, does not
preserve the solenoidal property in general), but rather
\begin{equation}
  \label{eqn:piola}
  \bm u = \frac{\bm J}{|\bm J|} (\hat{\bm u} \circ \chi_{\bm\mu})
\end{equation}
where $\bm J$ is the Jacobian of the mapping $\chi_{\bm \mu}^{-1}$. This modified
mapping must then be used in \eqref{eqn:war-1}--\eqref{eqn:war-2}. It should be
noted that the affine representations \eqref{eqn:split-a}--\eqref{eqn:split-d1}
may be considerably more complicated in this case, since $\bm J$ in general is
parameter-dependent.

We must also point out that to ensure the success of this method, the lift
functions $\hat{\bm \ell}(\bm \mu)$ \emph{must} also be solenoidal. Since the total
solution $\hat{\bm u} = \hat{\bm u}_0 + \hat{\bm \ell}$, and since the reduced spaces are all
concerned with the approximation of $\hat{\bm u}_0$ (and not $\hat{\bm u}$), this is
necessary to ensure the equivalence between the solenoidality of $\hat{\bm u}_0$ and
that of $\hat{\bm u}$. In practice one may solve the Stokes problem with a
divergence-conforming method to generate solenoidal lift functions. Since the lift is constructed in
the offline stage, the cost of this construction can be ignored.

With a fully solenoidal reduced velocity space, we are free to employ the
velocity-only formulation \eqref{eqn:velonly}, so that \eqref{eqn:block} reduces
to
\begin{equation}
  \bm A_\textsc{vv} \bm u = \bm f_\textsc{v}.
  \label{eqn:blocksolve-1}
\end{equation}
which is a system of size $M$ compared to \eqref{eqn:block} of size $3M$.

The supremizers of Section~\ref{sec:stab} function as the test functions in the
pressure recovery equation \eqref{eqn:sup-vvar}. Indeed, one obtains the system
\begin{equation}
  \bm B_\textsc{sp} \bm p = \bm f_\textsc{s} - \bm A_\textsc{sv} \bm u.
  \label{eqn:blocksolve-2}
\end{equation}
It is worth noting that the same expressions can be reached by setting $\bm B_\textsc{vp} = 0$ and
$\bm f_\textsc{p} = 0$ in \eqref{eqn:block}, obtaining
\begin{equation}
  \label{eqn:piolablock}
  \begin{pmatrix}
    \bm A_\textsc{vv} & \bm A_\textsc{vs} & \\
    \bm A_\textsc{sv} & \bm A_\textsc{ss} & \bm B_\textsc{sp} \\
    & \bm B_\textsc{ps} & \end{pmatrix}
  \begin{pmatrix} \bm u \\ \bm s \\ \bm p \end{pmatrix}
  = \begin{pmatrix} \bm r_\textsc{v} \\ \bm r_\textsc{s} \\ \bm 0 \end{pmatrix}
\end{equation}
If all spaces are of equal size $M$, the system \eqref{eqn:piolablock} is
block-triangular, and can be solved as
\eqref{eqn:blocksolve-1}--\eqref{eqn:blocksolve-2} together with $\bm s = 0$.
This illustrates the fact that supremizers serve no approximative purpose and
merely function as a test space for the pressure.

The system \eqref{eqn:blocksolve-1}--\eqref{eqn:blocksolve-2} also facilitates a
much faster reduced method than \eqref{eqn:block}, solving two systems of size
$M$ rather than one of size $3M$. We shall see in the following that this may
lead to speedups of more than an order of magnitude. It is also worth pointing out
that since reduced systems are full rather than dense, the benefits of this
reduction in size is amplified further compared to the usual situation of high-fidelity
sparse systems.

For other methods for pressure stabilization based on pure velocity solutions followed by pressure
recovery in various forms, see \cite{Akhtar2009ser,Caiazzo2014niv,Baiges2013ero}. Other methods for
pressure recovery from pure velocity solutions have been investigated. In \cite{Caiazzo2014niv},
three different methods based on modified equations and spaces are compared. The authors in
\cite{Akhtar2009ser} compute the pressure by taking the divergence of the Navier-Stokes equation,
projecting the result onto the pressure modes.

\subsection{Combined velocity-pressure basis}
\label{sec:combined}

One interesting method of pressure reconstruction is given in \cite{Baiges2013ero}, where all
reduced basis functions are \emph{joint} velocity \emph{and} pressure functions. The pressure is
then reconstructed from the coefficients of the solution vector in the reduced space, exactly as
with the velocity field. In the following, we will refer to this technique as
\emph{pressure reconstruction}, to distinguish it from \emph{pressure recovery} as discussed in
Section~\ref{sec:conforming}.

It is tempting to introduce this approach to divergence-conforming RBMs because the system
\eqref{eqn:blocksolve-1} remains the same, that is, adding pressure data to divergence-free velocity
basis functions does not change the reduced system.

\begin{remark}
  \label{rem:linearity}
  This method implies a linear dependence of pressure on velocity, an assumption that contradicts
  the Navier--Stokes equations. Its application should be restricted to linear models (e.g.~Stokes
  flow) or problems where the parameter space is suitably small for the linearization to be valid.
\end{remark}

There are at least two ways to construct a combined basis. One way is to omit the separation
discussed in Section~\ref{sec:fieldsep}, and instead pick $z$ in \eqref{eqn:covmx} as joint
velocity-pressure inner product, for example
\begin{equation}
  \label{eqn:z_c}
  z\big((\hat{\bm u}_0, \hat{p}), (\hat{\bm v}_0, \hat{q})\big) =
  \left( \hat{\bm u}_0, \hat{\bm v}_0 \right)_{H^1(\hat{\Omega})}
  + \rho \left( \hat{p}, \hat{q} \right)_{L^2(\hat{\Omega})}
\end{equation}
where $\rho > 0$ can be chosen to appropriately scale the relative importance of velocity
and pressure. The POD procedure will then produce joint velocity-pressure basis functions whose
velocity parts are divergence-free, so long as the snapshots are also divergence-free.

Another approach is to construct the velocity basis according to Section~\ref{sec:fieldsep}, and
then to use the same coefficients $\smash[t]{v_i^j}$ as in \eqref{eqn:uRj} to generate the pressure functions
in \eqref{eqn:pRj}. Note that this does not guarantee any kind of linear independence in the
pressure component.

With \eqref{eqn:z_c} one obtains combined velocity-pressure basis functions that are orthonormal in a
joint inner product. The second method produces the same velocity basis as the velocity-only
formulation discussed in Section~\ref{sec:conforming} with attached pressure data that may or may
not be linearly independent. However, because pressure reconstruction is based on direct linear
combination of the pressure modes, linear independence of these modes is not required.

In the following we have used the second approach for pressure reconstruction.

\section{Numerical experiments}
\label{sec:airfoil}

\begin{figure}
  \begin{center}
    \begin{tikzpicture}
      \draw[densely dotted, thick] (0,-3) arc (-90:90:3);
      \draw[thick] (0,3) arc (90:270:3);
      \draw[->] (-3,0) -- (-2.2,0);
      \draw[->] (-2.95,0.6) -- (-2.15,0.6);
      \draw[->] (-2.75,1.2) -- (-1.95,1.2);
      \draw[->] (-2.40,1.8) -- (-1.60,1.8);
      \draw[->] (-1.80,2.4) -- (-1.00,2.4);
      \draw[->] (-2.95,-0.6) -- (-2.15,-0.6);
      \draw[->] (-2.75,-1.2) -- (-1.95,-1.2);
      \draw[->] (-2.40,-1.8) -- (-1.60,-1.8);
      \draw[->] (-1.80,-2.4) -- (-1.00,-2.4);
      \node[anchor=east] at (-3,0) {$\bm u_\infty$};
      \draw[->] (1.6,0) arc (0:30:1.6);
      \draw[->] (1.6,0) arc (0:-30:1.6);
      \node[anchor=west] at (1.6,0) {$\varphi$};
      \begin{scope}[scale=0.3, shift={(-3,-1.4)}]
        \begin{axis}[xmin=-0.01, xmax=1.01, ymin=-0.2, ymax=0.2, unit vector ratio*=1 1, axis lines=none]
          \addplot[black, line width=2.5]
          table[x index={0}, y index={1}]{data/NACApts.dat};
        \end{axis}
      \end{scope}
      \begin{scope}[scale=0.3, rotate=30, shift={(-3,-1.5)}]
        \begin{axis}[xmin=-0.01, xmax=1.01, ymin=-0.2, ymax=0.2, unit vector ratio*=1 1, axis lines=none]
          \addplot[black, line width=1.5, dotted]
          table[x index={0}, y index={1}]{data/NACApts.dat};
        \end{axis}
      \end{scope}
      \begin{scope}[scale=0.3, rotate=-30, shift={(-3,-1.3)}]
        \begin{axis}[xmin=-0.01, xmax=1.01, ymin=-0.2, ymax=0.2, unit vector ratio*=1 1, axis lines=none]
          \addplot[black, line width=1.5, dotted]
          table[x index={0}, y index={1}]{data/NACApts.dat};
        \end{axis}
      \end{scope}
    \end{tikzpicture}
  \end{center}
  \caption{
    Sketch of the airfoil flow problem and its parameters. Solid lines indicate Dirichlet
    boundaries, and dotted lines indicate Neumann boundaries.
  }
  \label{fig:airfoil-sketch}
\end{figure}
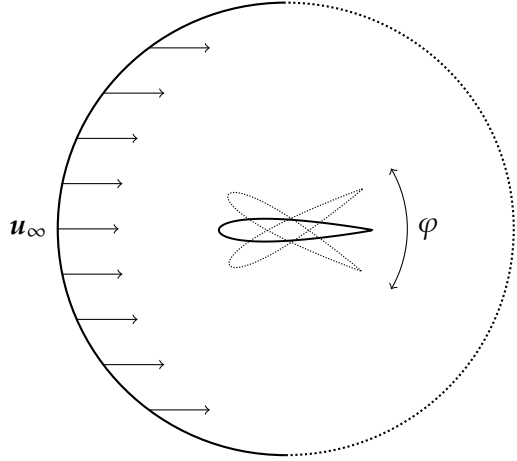

We consider a two-dimensional airfoil suspended in a mesh at the origin. We wish
to perform airflow simulations around this airfoil with the angle of attack
$\varphi$ and the inflow velocity $u_\infty$ as parameters, in other words
$\bm \mu = (\varphi, u_\infty)$; See Figure~\ref{fig:airfoil-sketch} for a
sketch. The purpose of the present section is to develop an affine
representation of this problem, as understood by
\eqref{eqn:split-a}--\eqref{eqn:split-d1}.

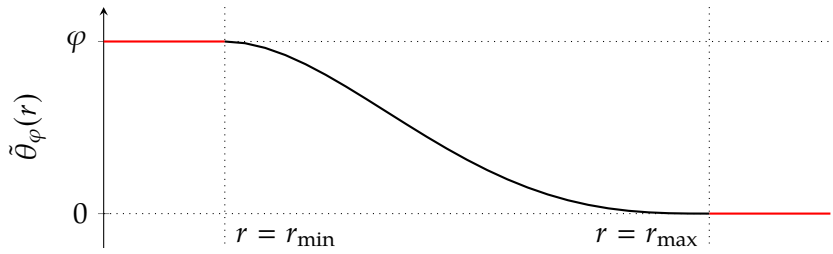
\begin{figure}
  \begin{center}
    \begin{tikzpicture}
      \begin{axis}[
        ylabel={$\tilde{\theta}_\varphi(r)$},
        width=0.7\textwidth,
        height=0.3\textwidth,
        axis lines=left,
        hide x axis,
        ymin=-0.2, ymax=1.2,
        ytick={0,1},
        yticklabels={$0$,$\varphi$},
        ]
        \addplot[thick, red, domain=0:1] {1};
        \addplot[thick, red, domain=5:6] {0};
        \addplot[thick, domain=1:5] {(1-((x-1)/4))^3 * (3*(x-1)/4+1)};
        \draw[dotted] (axis cs:1,1.2) -- (axis cs:1,-0.2);
        \draw[dotted] (axis cs:5,1.2) -- (axis cs:5,-0.2);
        \draw[dotted] (axis cs:0,0) -- (axis cs:5,0);
        \draw[dotted] (axis cs:1,1) -- (axis cs:6,1);
        \node[anchor=south west] at (axis cs:1,-0.25) {$r=r_\text{min}$};
        \node[anchor=south east] at (axis cs:5,-0.25) {$r=r_\text{max}$};
      \end{axis}
    \end{tikzpicture}
  \end{center}
  \caption{$\tilde{\theta}_\varphi(r)$ as a function of $r$.}
  \label{fig:thetatilde}
\end{figure}

While the velocity is trivial, the geometric
parameter can be accomplished by rotating the mesh through an angle that depends
on the distance from the center of the airfoil (hereafter referred to as the
\emph{radius} $r$), i.e.~$\tilde{\theta}_\varphi = \tilde{\theta}_\varphi(r)$
so that
\[
  \tilde{\theta}_\varphi(r) = \begin{cases}
    \varphi, & r \le r_\text{min}, \\
    0, & r \ge r_\text{max},
  \end{cases}
\]
and which smoothly interpolates between the two in the region
$r_\text{min} < r < r_\text{max}$. See Figure~\ref{fig:thetatilde} for an
example. Note that we also require
\[
  \tilde{\theta}_\varphi'(r_\text{min}) =
  \tilde{\theta}_\varphi'(r_\text{max}) = 0
\]
in accordance with Remark~\ref{rem:lift}, i.e.~the mapping of lifts between reference and
transformed domains will not modify the boundary data.
It is desirable that the entire airfoil is contained within the region $r <
r_\text{min}$, and that the boundary of the mesh is entirely inside the region
$r > r_\text{max}$. Given these restrictions, the mesh may otherwise be freely
chosen.

In the following we will denote by $\theta(r)$ a \emph{canonical} angle
function, satisfying
\[
  \theta(r) = \begin{cases}
    1, & r \le r_\text{min}, \\
    0, & r \ge r_\text{max}.
  \end{cases}
\]
With this, we can write $\tilde{\theta}_\varphi(r) = \varphi \theta(r)$.
The mapping from reference coordinates to physical coordinates can then be
expressed as
\[
  \chi_{\varphi}^{-1} \hat{\bm r}
  = \bm R(\varphi \theta(r)) \; \hat{\bm r},
\]
where $\hat{\bm r} = \left( \hat{x}, \hat{y} \right)$,
$r = \| \hat{\bm r} \| = \sqrt{\hat{x}^2 + \hat{y}^2}$
and $\bm R$ is the rotation matrix,
\[
  \bm R(a) = \begin{pmatrix} \cos a & - \sin a \\ \sin a & \cos a \end{pmatrix}.
\]

The Jacobian $\bm J = \bm J(\chi_{\varphi}^{-1})$ can be expressed,
using $a = \varphi \theta(r)$ as shorthand, as
\begin{align}
  \nonumber
  \bm J &= \bm R(a) + \bm R'(a) \hat{\bm r} \nabla a^\intercal
  = \bm R(a) \left( \bm I + \bm P \hat{\bm r} \nabla a^\intercal \right) \\
  &= \bm R(a) \left(
    \bm I + \frac{\varphi \theta'}{r} \bm P \hat{\bm r} \hat{\bm r}^\intercal
  \right)
  = \bm R(a) \left( \bm I + \varphi \bm P \bm Q \right)
\end{align}
where $\bm R'(a) = \partial \bm R(a) / \partial a$ and
where the two utility matrices $\bm P$ and $\bm Q$ are defined as
\begin{equation}
  \bm P = \bm R(\nicefrac{\pi}{2}), \qquad
  \bm Q = \frac{\theta'}{r} \hat{\bm r} \hat{\bm r}^\intercal,
\end{equation}
noting the useful property that $\bm R \bm P = \bm R'$.

The determinant of $\bm R$ is $1$ and the determinant of $\bm J$ follows
from the matrix determinant lemma (as a one-rank update to an invertible
matrix),
\[
  \det \bm J = 1 + \frac{\varphi\theta'}{r} \hat{\bm r}^\intercal \bm P \hat{\bm r} = 1,
\]
since the last term may be recognized as the inner product between two
orthogonal vectors. The inverse of $\bm J$ follows from the Sherman-Morrison theorem,
\[
  \bm J^{-1} = \left( \bm I - \varphi \bm P \bm Q \right)
  \bm R(a)^\intercal.
\]
It can then be seen that the problem of finding affine representations of
$\bm J$ and $\bm J^{-1}$ reduces to the affine representation of $\bm R$.
We therefore proceed as follows.
\[
  \bm R(a)
  = \begin{pmatrix} \cos a & -\sin a \\ \sin a & \cos a \end{pmatrix}
  = \sum_{i=0}^\infty \left( \frac{(-1)^i a^{2i}}{(2i)!} \bm I
    + \frac{(-1)^i a^{2i+1}}{(2i+1)!} \bm P \right).
\]
Noting that, since $\bm P^0 = \bm I$, $\bm P^1 = \bm P$, $\bm P^2 = -\bm I$,
$\bm P^3 = -\bm P$ etc., we obtain
\[
  (-1)^i \bm I = \bm P^{2i}, \qquad (-1)^i \bm P = \bm P^{2i+1},
\]
so the series expansion can be more succinctly written as
\begin{equation}
  \label{eqn:rotsum}
  \bm R(a) = \sum_{i=0}^\infty \frac{a^i}{i!} \bm P^i
  = \sum_{i=0}^\infty \varphi^i \underbrace{\frac{\theta^i}{i!}\bm P^i}_{\bm R_i}
  = \sum_{i=0}^\infty \varphi^i \bm R_i.
\end{equation}

To investigate the effect of the Piola transformation, we consider the affine representations
\eqref{eqn:split-a}--\eqref{eqn:split-c0} for two different methods. First, a conventional method
using the pullback via $\pi_{\bm \mu}^\textsc{v}$ to map between function spaces, and second a
divergence-conforming method based on the Piola mapping \eqref{eqn:piola}.  In doing so we will
truncate the series for $\bm R(a)$ to a number of terms that can achieve the desired precision. It
is important to note that the system \eqref{eqn:piolablock} only obtains its desired form for high
accuracy approximations, and that if the affine representation form for $b(\cdot,\cdot,\bm \mu)$ is
not exact, the matrix $\bm B_\textsc{vp}$ from \eqref{eqn:block} will correspondingly not be exactly
zero, and the solution of \eqref{eqn:blocksolve-1}--\eqref{eqn:blocksolve-2} will not agree with the
solution of \eqref{eqn:block}.

The error made by truncating \eqref{eqn:rotsum} to $n$ terms is approximately
$\varphi_\textsc{max}^n / n!$. For a maximal angle of attack of
$\SI{35}{\degree}$ for example, we can expect about $10$ digits of accuracy
with $n=10$ terms.

In the following, only the parameter $\varphi$ is considered, as $u_\infty$ is
comparatively trivial. In all transformations, the Jacobian determinant $|\bm J|$ is identically
equal to $1$.

\subsection{Non-Piola formulation}

For the Laplacian form $a$ we get
\begin{align}
  ({\bm\pi}^*_{\bm\mu}a)(
    \hat{\bm u},
    \hat{\bm w};
    \varphi
  )
  &= \nu \int_{\hat{\Omega}} (\bm J^{-\intercal} \nabla) \hat{\bm u} : (\bm J^{-\intercal} \nabla)
    \hat{\bm w} |\bm J|
  = \nu \int_{\hat{\Omega}} \nabla \hat{\bm u} : (\bm J^{-1} \bm J^{-\intercal} \nabla) \hat{\bm w}
    |\bm J|,
\end{align}
and we find for the matrix $\bm J^{-1} \bm J^{-\intercal}$ that
\begin{align}
  \nonumber
  \bm J^{-1} \bm J^{-\intercal}
  &= (\bm I - \varphi \bm P \bm Q) \bm R^\intercal
    \bm R (\bm I + \varphi \bm Q \bm P) \\
  \nonumber
  &= \bm I + \varphi \underbrace{(\bm Q \bm P - \bm P \bm Q)}_{\bm D_1}
    - \varphi^2 \underbrace{\bm P \bm Q^2 \bm P}_{\bm D_2} \\
  &= \bm I + \varphi \bm D_1 - \varphi^2 \bm D_2.
\end{align}
Since $|\bm J| = 1$, this gives a three-term affine representation of $a$ as
\begin{equation}
  ({\bm\pi}^*_{\bm\mu}a)(
    \hat{\bm u},
    \hat{\bm w};
    \varphi
  ) =
  \nu \int_{\hat{\Omega}} \nabla \hat{\bm u} : \nabla \hat{\bm w}
  + \nu \varphi \int_{\hat{\Omega}}
  \nabla \hat{\bm u} : (\bm D_1 \nabla) \hat{\bm w}
  - \nu \varphi^2 \int_{\hat{\Omega}}
  \nabla \hat{\bm u} : (\bm D_2 \nabla) \hat{\bm w}.
\end{equation}

For the divergence form $b$ we have
\begin{equation}
  ({\bm\pi}^*_{\bm\mu}b)(
    \hat{p},
    \hat{\bm w};
    \varphi
  ) =
  \int_{\hat{\Omega}} \hat{p} (\bm J^{-\intercal} \nabla) \cdot \hat{\bm w} |\bm J|
  = \int_{\hat{\Omega}} \hat{p} \bm J^{-\intercal} : \nabla \hat{\bm w} |\bm J|,
\end{equation}
meaning we need a series representation of $\bm J^{-\intercal}$:
\begin{align}
  \nonumber
  \bm J^{-\intercal}
  &= \bm R(a) (\bm I - \varphi \bm Q^\intercal \bm P^\intercal)
  = \bm R(a) (\bm I + \varphi \bm Q \bm P)
  = \sum_{i=0}^\infty
    \varphi^i \bm R_i
    (\bm I + \varphi \bm Q \bm P) \\
  &= \sum_{i=0}^\infty
    \varphi^i \bm R_i
    + \varphi^{i+1} \bm R_i \bm Q \bm P
  = \sum_{i=0}^\infty
    \varphi^i \underbrace{\left(
    \bm R_i + \bm R_{i-1} \bm Q \bm P
    \right)}_{\bm B^{(-)}_i}
  = \sum_{i=0}^\infty \varphi^i \bm B^{(-)}_i,
\end{align}
with the understanding that $\bm R_{-1} = 0$. This gives an affine
representation of $b$ in $2n$ terms, where $n$ is a suitable number of terms for
a truncated Taylor series for $\sin$ or $\cos$, given the range of $\varphi$
under consideration.
\begin{equation}
  ({\bm\pi}^*_{\bm\mu}b)(
    \hat{p},
    \hat{\bm w};
    \varphi
  ) \approx \sum_{i=0}^{2n} \varphi^i
  \int_{\hat{\Omega}} \hat{p} \bm B^{(-)}_i : \nabla \hat{\bm w}
\end{equation}
The same expansion will work with the convective term $c$, viz.
\begin{equation}
  ({\bm\pi}^*_{\bm\mu}c)(
    \hat{\bm u},
    \hat{\bm v},
    \hat{\bm w};
    \varphi
  )
  = \int_{\hat{\Omega}} (\hat{\bm u} \cdot \bm J^{-\intercal}\nabla) \hat{\bm v} \cdot \hat{\bm w}
  |\bm J|
  \approx \sum_{i=0}^{2n} \varphi^i \int_{\hat{\Omega}}
    (\hat{\bm u} \cdot \bm B^{(-)}_i \nabla) \hat{\bm v} \cdot \hat{\bm w}.
\end{equation}

\subsection{Piola formulation}

This proceeds as for the non-Piola case, except every vector field is pre-multiplied with the
Jacobian (recall, the determinant is $1$, which simplifies \eqref{eqn:piola} somewhat). For this, we
need an affine representation of $\bm J$. It looks deceptively like that of $\bm J^{-\intercal}$:
\begin{align}
  \nonumber
  \bm J
  &= \bm R(a) (\bm I + \varphi \bm P \bm Q)
  = \sum_{i=0}^\infty \varphi^i \bm R_i (\bm I + \varphi \bm P \bm Q)
  = \sum_{i=0}^\infty \varphi^i \bm R_i + \varphi^{i+1} \bm P \bm Q \\
  &= \sum_{i=0}^\infty
    \varphi^i \underbrace{\left(
    \bm R_i + \bm R_{i-1} \bm P \bm Q
    \right)}_{\bm B^{(+)}_i} = \sum_{i=0}^\infty \varphi^i \bm B^{(+)}_i.
\end{align}
For the form $b$ we then get
\begin{align}
  \nonumber
  ({\bm\pi}^*_{\bm\mu}b)(
    \hat{p},
    \hat{\bm w};
    \varphi
  ) &= \int_{\hat{\Omega}} \hat{p} \bm J^{-\intercal} : \nabla (\bm J \hat{\bm w}) |\bm J|
    \approx \int_{\hat{\Omega}} \hat{p}
      \left( \sum_{i=0}^{2n} \varphi^i \bm B^{(-)}_i \right) : \nabla
      \left( \sum_{j=0}^{2n} \varphi^j \bm B^{(+)}_j \hat{\bm w} \right) \\
    &= \sum_{i,j=0}^{2n} \varphi^{i+j}
      \int_{\hat{\Omega}} \hat{p} \bm B^{(-)}_i :
      \nabla \left( \bm B^{(+)}_j \hat{\bm w} \right).
\end{align}
Given truncated expressions for $\bm J$ and $\bm J^{-\intercal}$ with $2n$
terms, this is an affine representation with $4n$ terms.

The form $c$ has similar complexity.
\begin{align}
  \nonumber
  ({\bm\pi}^*_{\bm\mu}c)(
    \hat{\bm u},
    \hat{\bm v},
    \hat{\bm w};
    \varphi
  )
  &= \int_{\hat{\Omega}} (\bm J \hat{\bm u} \cdot \bm J^{-\intercal}\nabla)
    \bm J \hat{\bm v} \cdot \bm J \hat{\bm w} |\bm J|
  \approx \int_{\hat{\Omega}} (\hat{\bm u} \cdot \nabla)
      \left( \sum_{i=0}^{2n} \varphi^i \bm B^{(+)}_i \hat{\bm v} \right) \cdot
    \left( \sum_{i=0}^{2n} \varphi^j \bm B^{(+)}_j \hat{\bm w} \right) \\
  &= \sum_{i,j=0}^{2n} \varphi^{i+j}
    \int_{\hat{\Omega}} (\hat{\bm u} \cdot \nabla) \bm B^{(+)}_i \hat{\bm v} \cdot \bm B^{(+)}_j \hat{\bm w}.
\end{align}

And finally, $a$ can be represented as
\begin{align}
  \nonumber
  ({\bm\pi}^*_{\bm\mu}a)(
    \hat{\bm u},
    \hat{\bm w};
    \varphi
  ) &= \int_{\hat{\Omega}}
      (\bm J^{-\intercal} \nabla) (\bm J \hat{\bm u}) :
      (\bm J^{-\intercal} \nabla) (\bm J \hat{\bm w}) |\bm J|
    = \int_{\hat{\Omega}}
      \nabla (\bm J \hat{\bm u}) :
      (\bm J^{-1} \bm J^{-\intercal} \nabla) (\bm J \hat{\bm w}) |\bm J| \\
  \nonumber
    &\approx \int_{\hat{\Omega}}
      \nabla \left( \sum_{i=0}^{2n} \varphi^i \bm B^{(+)}_i \hat{\bm u} \right) :
      ((\bm I + \varphi \bm D_1 - \varphi^2 \bm D_2) \nabla)
      \left( \sum_{j=0}^{2n} \varphi^j \bm B^{(+)}_j \hat{\bm w} \right) \\
  \nonumber
    &= \sum_{i,j=0}^{2n} \varphi^{i+j} \int_{\hat{\Omega}}
      \nabla (\bm B^{(+)}_i \hat{\bm u}) : \nabla (\bm B^{(+)}_j \hat{\bm w}) \\
  \nonumber
    &+ \sum_{i,j=0}^{2n} \varphi^{i+j+1} \int_{\hat{\Omega}}
      \nabla (\bm B^{(+)}_i \hat{\bm u}) : (\bm D_1 \nabla) (\bm B^{(+)}_j \hat{\bm w}) \\
    &- \sum_{i,j=0}^{2n} \varphi^{i+j+2} \int_{\hat{\Omega}}
      \nabla (\bm B^{(+)}_i \hat{\bm u}) : (\bm D_2 \nabla) (\bm B^{(+)}_j \hat{\bm w}),
\end{align}
which is an affine representation with $4n+2$ terms.

\section{Results}
\label{sec:results}

For both the Taylor-Hood and the conforming method an ensemble of 225 solutions was generated at
the $15 \times 15$ Gauss points for the parameter set
\[
  \mathcal{P} = \left\{ (\varphi,u_\infty) \;|\;
    \varphi \in [-\SI{35}{\degree},\SI{35}{\degree}],
    u_\infty \in [1, 20]
  \right\}.
\]
The viscosity was set at $\nu = \nicefrac{1}{6}$. The airfoil was chosen as a
NACA0015 profile with a chord length of $1$, thus giving an approximate maximal
Reynold's number of $120$.

The reference domain $\hat{\Omega}$ was chosen as an disk of radius $10$, centered at the
center of the airfoil, excluding the airfoil at $\varphi=\SI{0}{\degree}$. The mesh had
$120$ elements in the circumferential direction and $40$ elements in the
radial direction. Given discretization nodes $\{n^{\textsc{int}}_i\}_{i=1}^{120}$ for the airfoil, and uniformly
spaced nodes $\{n^{\textsc{ext}}_i\}_{i=1}^{120}$ for the external boundary, the interior nodes of the mesh were
positioned as
\begin{equation}
  \label{eqn:meshgen}
  n^{\textsc{dom}}_{ij} = \bigg(\frac{e^{\gamma j / 120} - 1}{e^{\gamma} - 1}\bigg) n^{\textsc{int}}_i
  + \bigg(1 - \frac{e^{\gamma j / 120} - 1}{e^{\gamma} - 1}\bigg) n^{\textsc{ext}}_i
\end{equation}
where $\gamma = 120 \log{\alpha}$ and $\alpha$ is a grading factor. This produces a mesh where the elements
grow in size by a factor of approximately $\alpha$ for each radial layer, starting from the airfoil. We
used $\alpha = 1.2$.

The boundary conditions were enforced strongly at the inflow boundary,
with no-slip conditions on the airfoil and a do-nothing homogeneous Neumann
boundary condition on the outflow. The lift function was chosen as the
solution to the Stokes problem with $\varphi=0$, which is solenoidal.

\begin{figure}
  \begin{center}
    \includegraphics[width=0.5\textwidth]{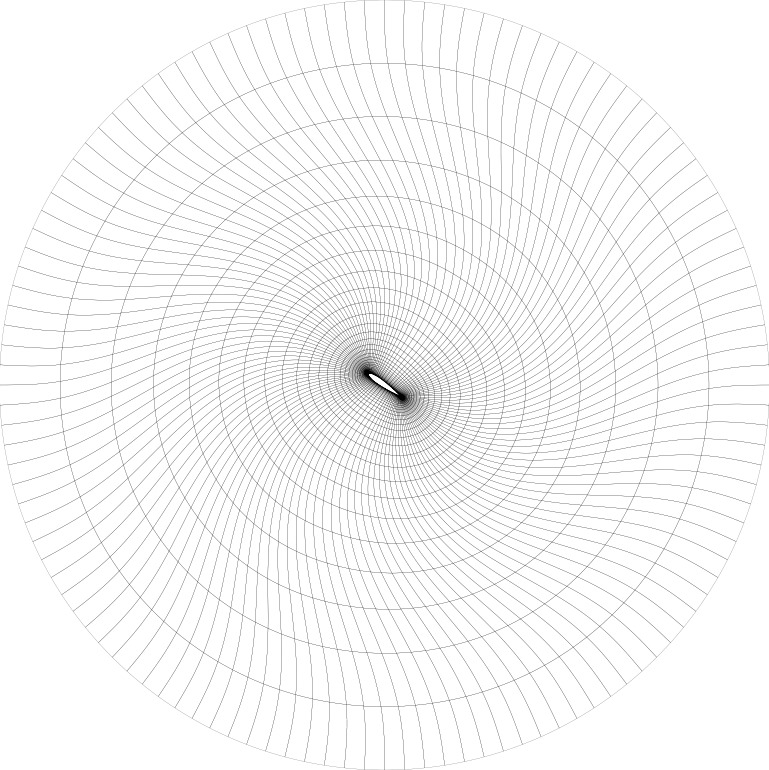}
  \end{center}
  \caption{Sample domain with $\varphi = -\nicefrac{\pi}{4}$.}
  \label{fig:domain}
\end{figure}

For parametrizing the geometry, the radius-dependent rotation angle function
$\theta$ was chosen as
\begin{align}
  \theta(r) = \begin{cases}
    1, & r < r_\text{min}, \\
    0, & r > r_\text{max}, \\
    (1-\overline{r})^3 (3\overline{r}+1), & \text{otherwise},
  \end{cases}
\end{align}
with
\[
  \overline{r} = \frac{(r - r_\text{min})}{(r_\text{max} - r_\text{min})},
  \qquad r_\text{min} = 1, \qquad r_\text{max} = 10,
\]
see Figures~\ref{fig:thetatilde} and \ref{fig:domain}.

For the conventional high-fidelity model, we chose a Lagrangian bi-quadratic basis for velocity,
and Lagrangian bilinear basis for the pressure. This fulfills the Taylor--Hood property, leading to an
inf-sup stable formulation. For the conforming high-fidelity model, we chose a bilinear Lagrangian basis for the
pressure, and a mixed quadratic and linear B-spline basis for the velocity, giving a fully
divergence-conforming method for the reference geometry at $\varphi=0$ (see
\cite{Evans2013idc1}). Note that by design, the divergence-conforming method will retain this
property for other angles.

For solving the nonlinear equation we used Newton iteration, stopping when the
velocity update reached $10^{-10}$ or less, as measured in the $H^1$ seminorm.
The affine representations were derived from a truncated version of
\eqref{eqn:rotsum} with $n=12$ terms, sufficient to represent the rotation
matrix $\bm R(a)$ to $10$ digits accuracy within the range of angles considered.

Results have been generated for reduced models with $M=10,20,\ldots,50$ degrees
of freedom each in the three spaces (velocity, supremizers and pressure), as
well as the corresponding un-stabilized models (only velocity and pressure).
Additionally, for the conforming stabilized method, a distinction is made
between a naive solver, a block solver
(see \eqref{eqn:blocksolve-1}--\eqref{eqn:blocksolve-2})
and a combined solver
(using the pressure reconstruction technique in Section~\ref{sec:combined}).

\subsection{Spectrum}

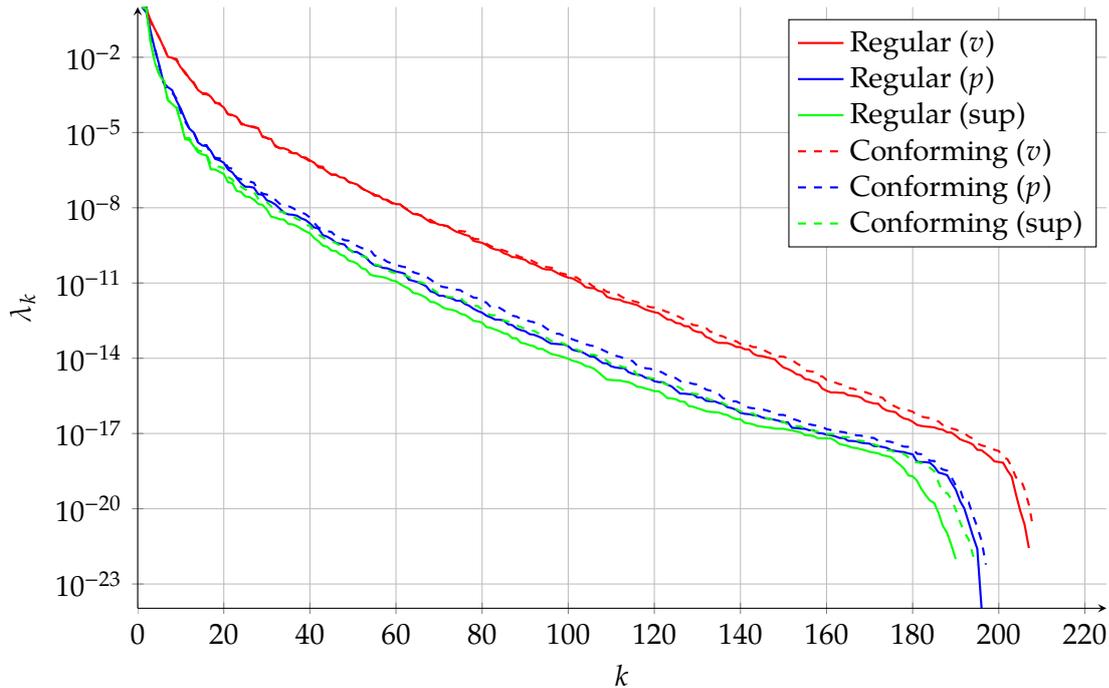
\begin{figure}
  \begin{tikzpicture}
    \begin{axis}[
      xlabel={$k$},
      ylabel={$\lambda_k$},
      ymode=log,
      xmin=0, xmax=225,
      width=0.9\textwidth,
      height=0.6\textwidth,
      grid=both,
      axis lines=left,
      legend style={
        at={(0.99, 0.98)},
        anchor=north east,
      },
      legend cell align=left,
      ]
      \addplot[red, thick]
      table[x index={0}, y index={1}]{data/airfoil-spectrum-no-piola.csv};
      \addplot[blue, thick]
      table[x index={0}, y index={2}]{data/airfoil-spectrum-no-piola.csv};
      \addplot[green, thick]
      table[x index={0}, y index={3}]{data/airfoil-spectrum-no-piola.csv};
      \addplot[red, thick, dashed]
      table[x index={0}, y index={1}]{data/airfoil-spectrum-piola.csv};
      \addplot[blue, thick, dashed]
      table[x index={0}, y index={2}]{data/airfoil-spectrum-piola.csv};
      \addplot[green, thick, dashed]
      table[x index={0}, y index={3}]{data/airfoil-spectrum-piola.csv};
      \legend{
        Regular ($v$),
        Regular ($p$),
        Regular (sup),
        Conforming ($v$),
        Conforming ($p$),
        Conforming (sup),
      }
    \end{axis}
  \end{tikzpicture}
  \caption{
    Ensemble spectra for 225 ensemble solutions for both methods, and all three
    spaces (velocity, pressure and supremizers).
  }
  \label{fig:spectra}
\end{figure}

The decay rate of the eigenvalue spectrum of the solution ensembles indicates to which degree one
might expect a reduced basis to adequately capture the most significant behavioral patterns of the
model. As can be seen from Figure~\ref{fig:spectra}, the decay is rapid for the first $20$ or so
modes, and then flattens out somewhat after that. (There is reason to believe that this might be
improved for different choices of the angle function $\theta$.) The spectra for the supremizers
closely mimic those of the pressure solutions, as should be expected.

\begin{remark}
  The difference in the pressure spectra can be partially explained by quadrature. We have used the
  same quadrature rule (9-point Gaussian quadrature) for both methods. However, the conforming
  method has higher order integrands than the Taylor-Hood method, and thus requires stronger
  quadrature rules. The relative quadrature error for the conforming method is estimated at roughly
  $\SI{e-4}{}$ for both velocity and pressure, while the same errors for the Taylor-Hood method are
  $\SI{e-5}{}$ and $\SI{e-6}{}$, respectively.
\end{remark}

\subsection{Basis functions}

\begin{figure}
  \begin{center}
    \adjustbox{trim={0.13\width} {0.12\height} {0.13\width} {0.12\height}, clip}{
      \includegraphics[width=0.4\textwidth]{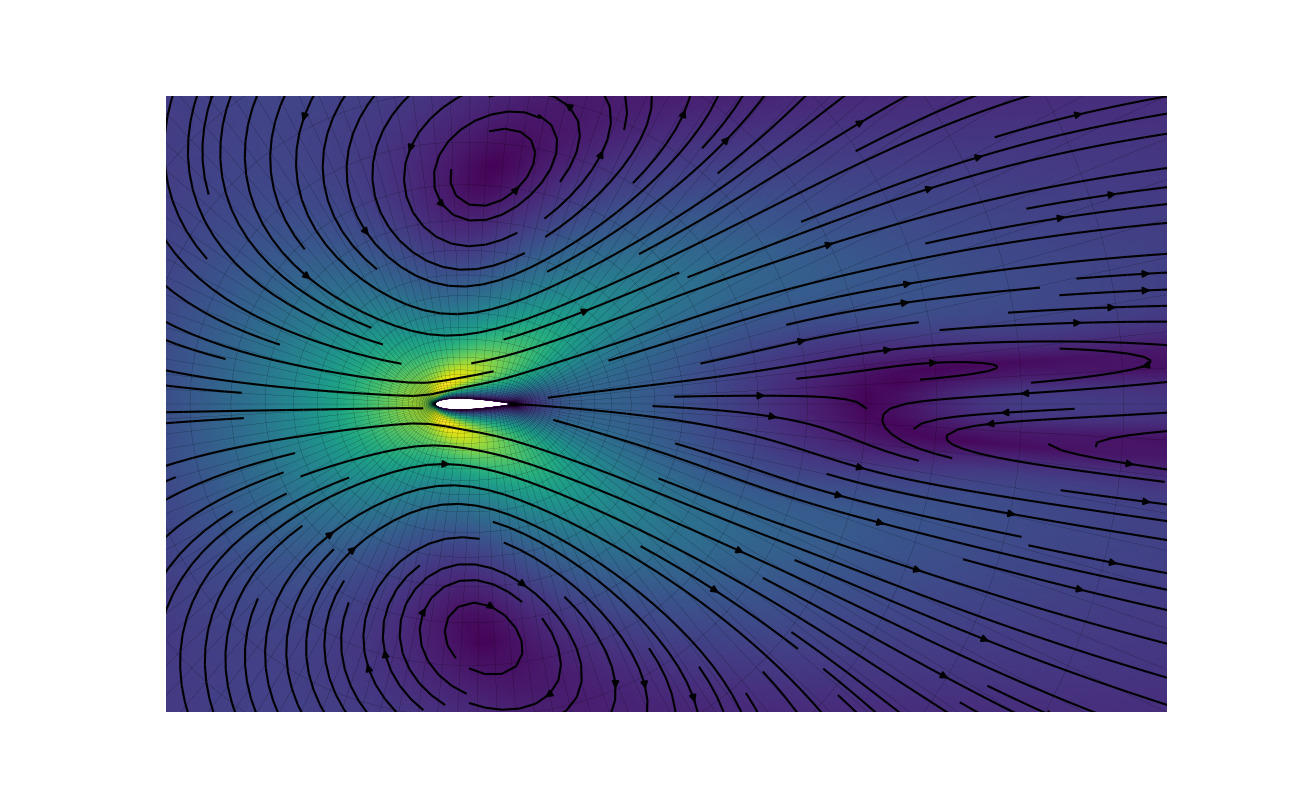}
    }
    \adjustbox{trim={0.13\width} {0.12\height} {0.13\width} {0.12\height}, clip}{
      \includegraphics[width=0.4\textwidth]{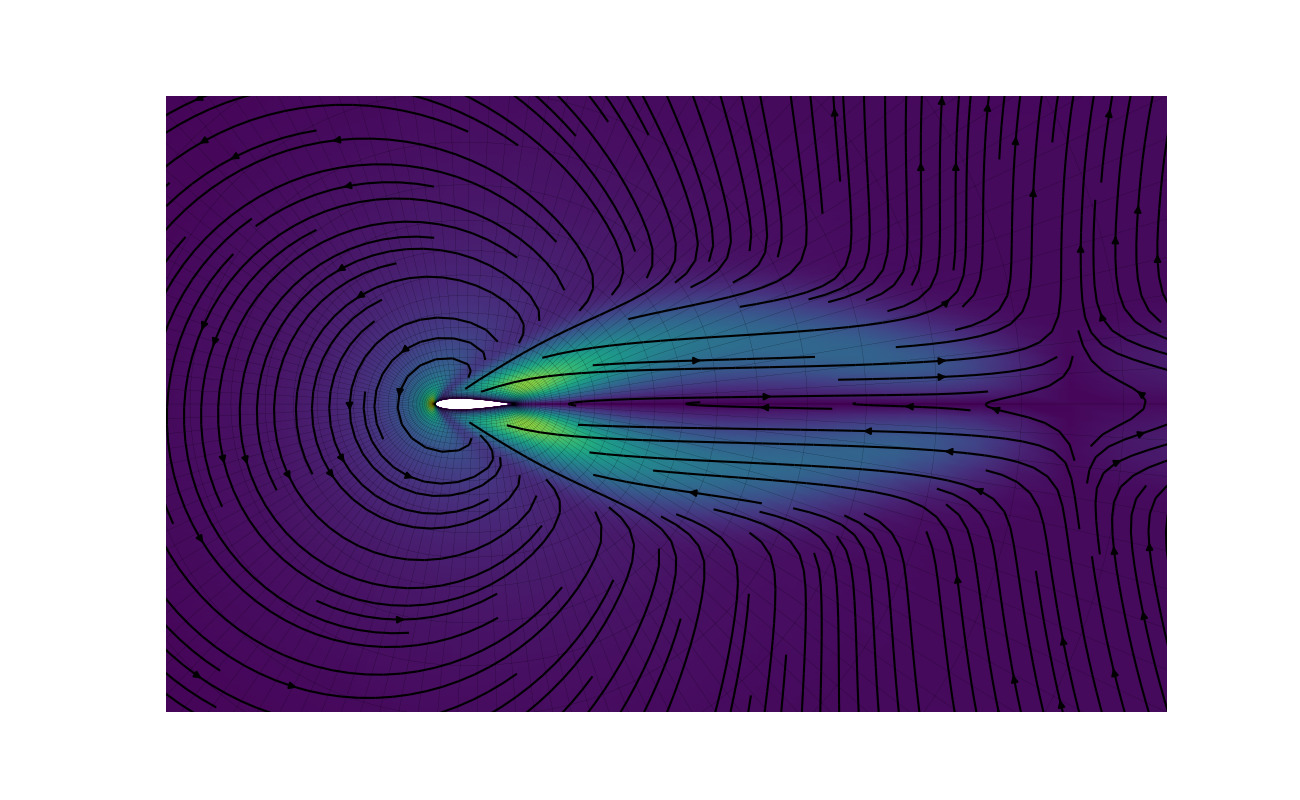}
    }
    \adjustbox{trim={0.13\width} {0.12\height} {0.13\width} {0.12\height}, clip}{
      \includegraphics[width=0.4\textwidth]{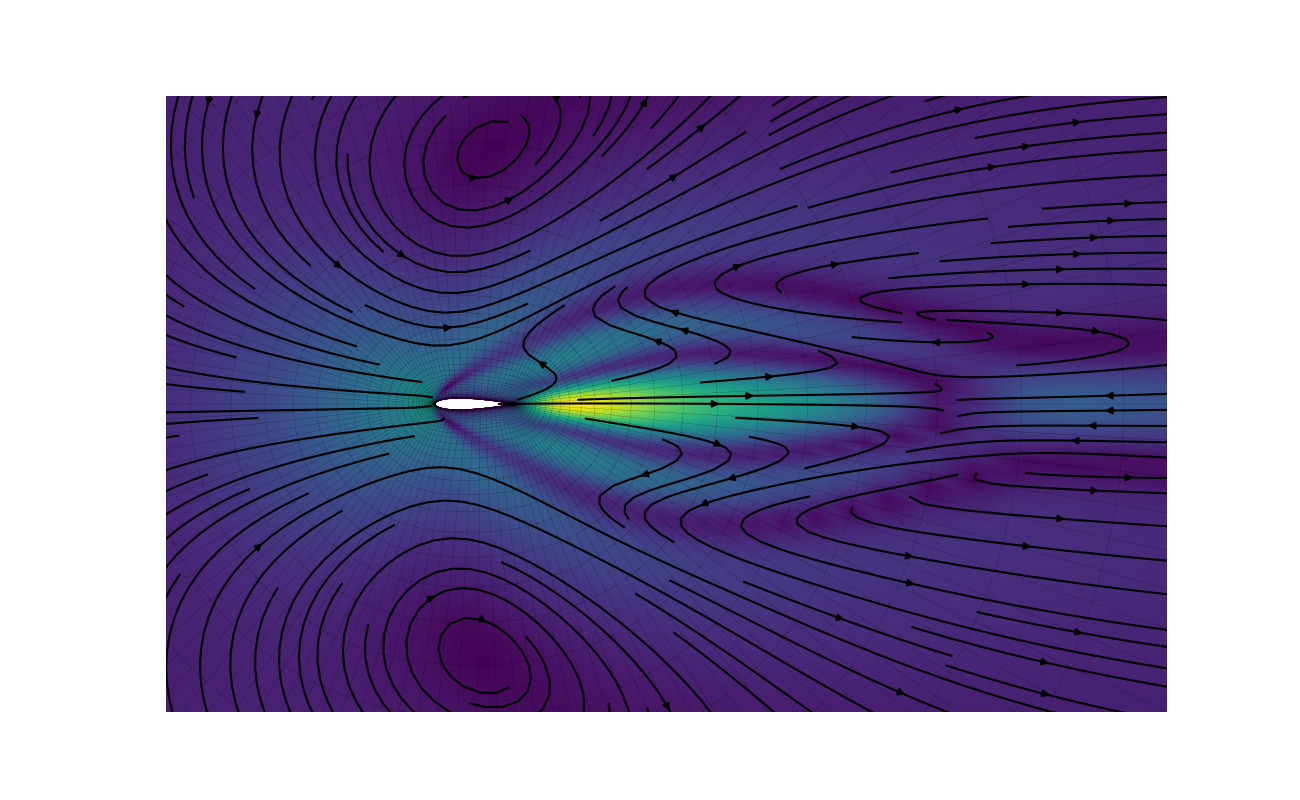}
    } \\
    \adjustbox{trim={0.13\width} {0.12\height} {0.13\width} {0.12\height}, clip}{
      \includegraphics[width=0.4\textwidth]{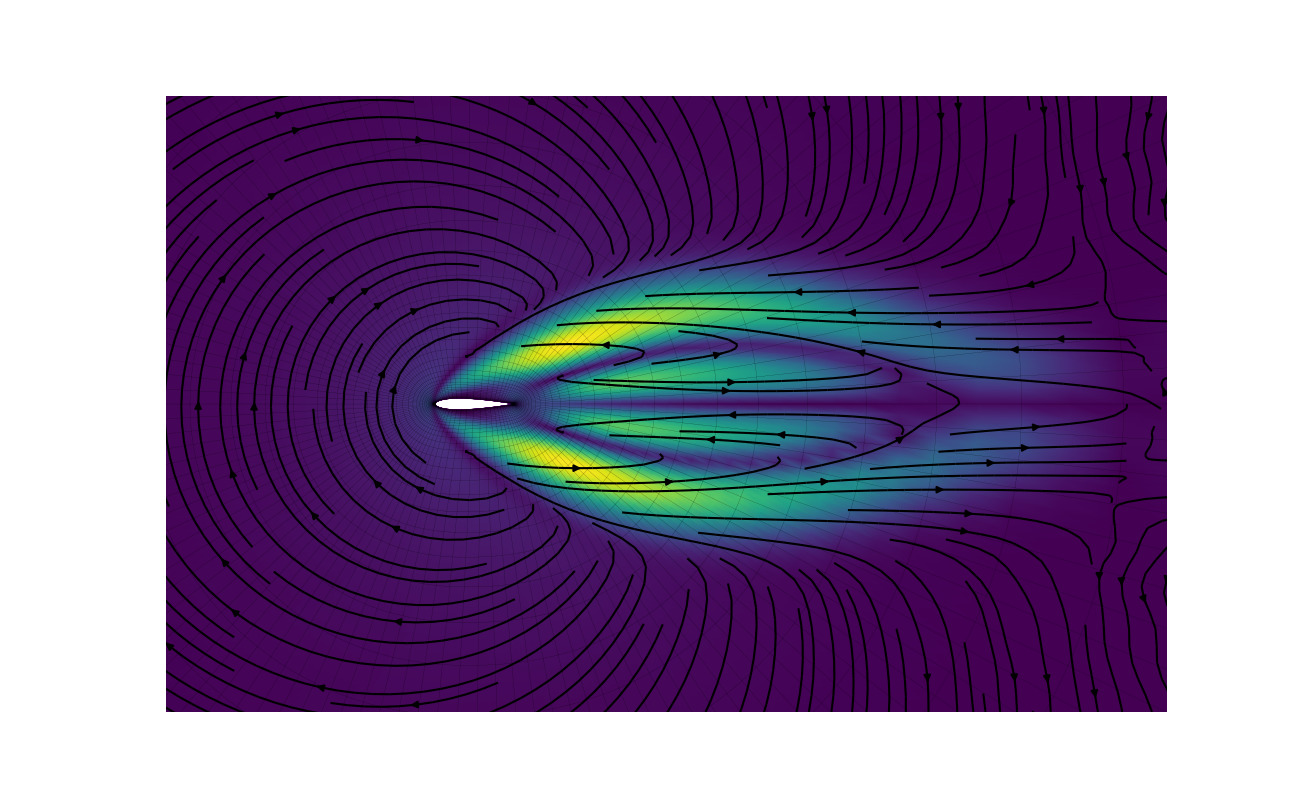}
    }
    \adjustbox{trim={0.13\width} {0.12\height} {0.13\width} {0.12\height}, clip}{
      \includegraphics[width=0.4\textwidth]{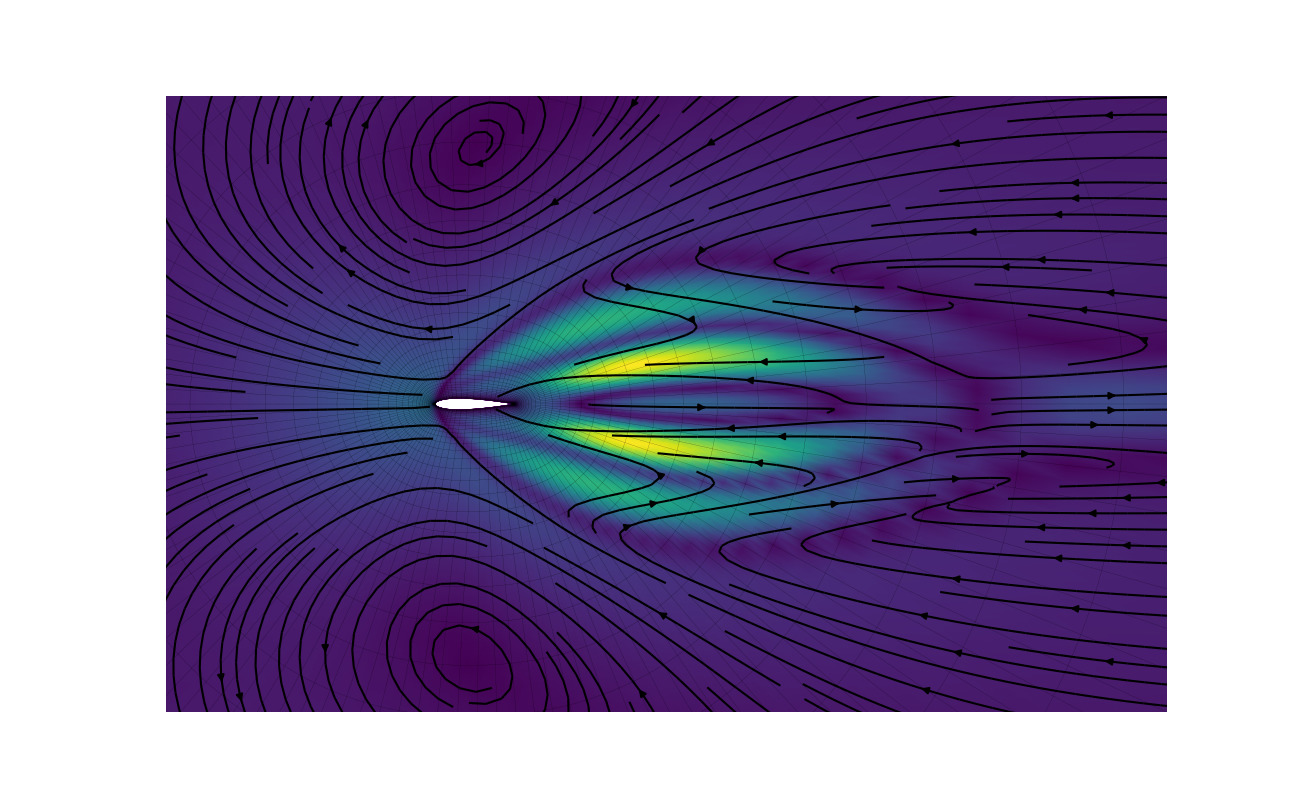}
    }
    \adjustbox{trim={0.13\width} {0.12\height} {0.13\width} {0.12\height}, clip}{
      \includegraphics[width=0.4\textwidth]{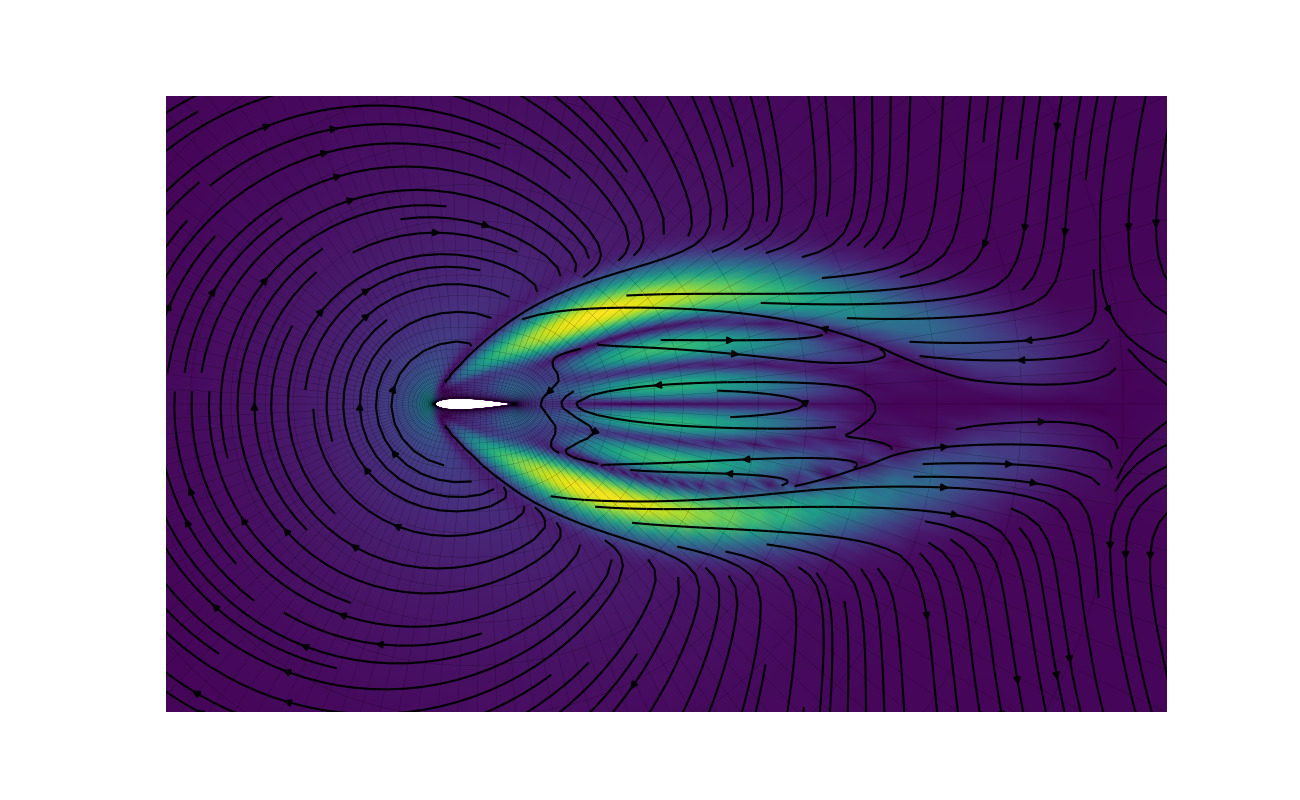}
    }
    \caption{First six velocity modes (Taylor-Hood method).}
    \label{fig:modes-vel-noconf}
    \vspace{\floatsep}
    \adjustbox{trim={0.3\width} {0.3\height} {0.3\width} {0.3\height}, clip}{
      \includegraphics[width=0.4\textwidth]{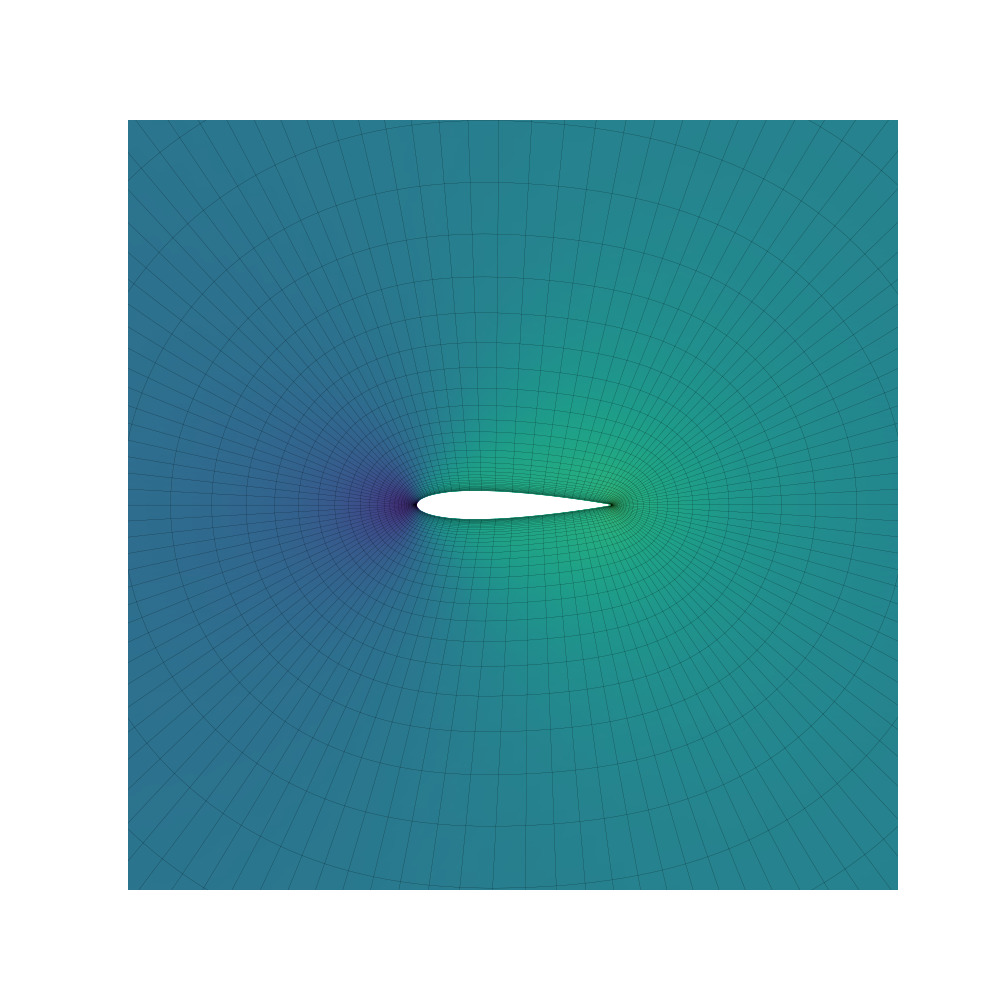}
    }
    \adjustbox{trim={0.3\width} {0.3\height} {0.3\width} {0.3\height}, clip}{
      \includegraphics[width=0.4\textwidth]{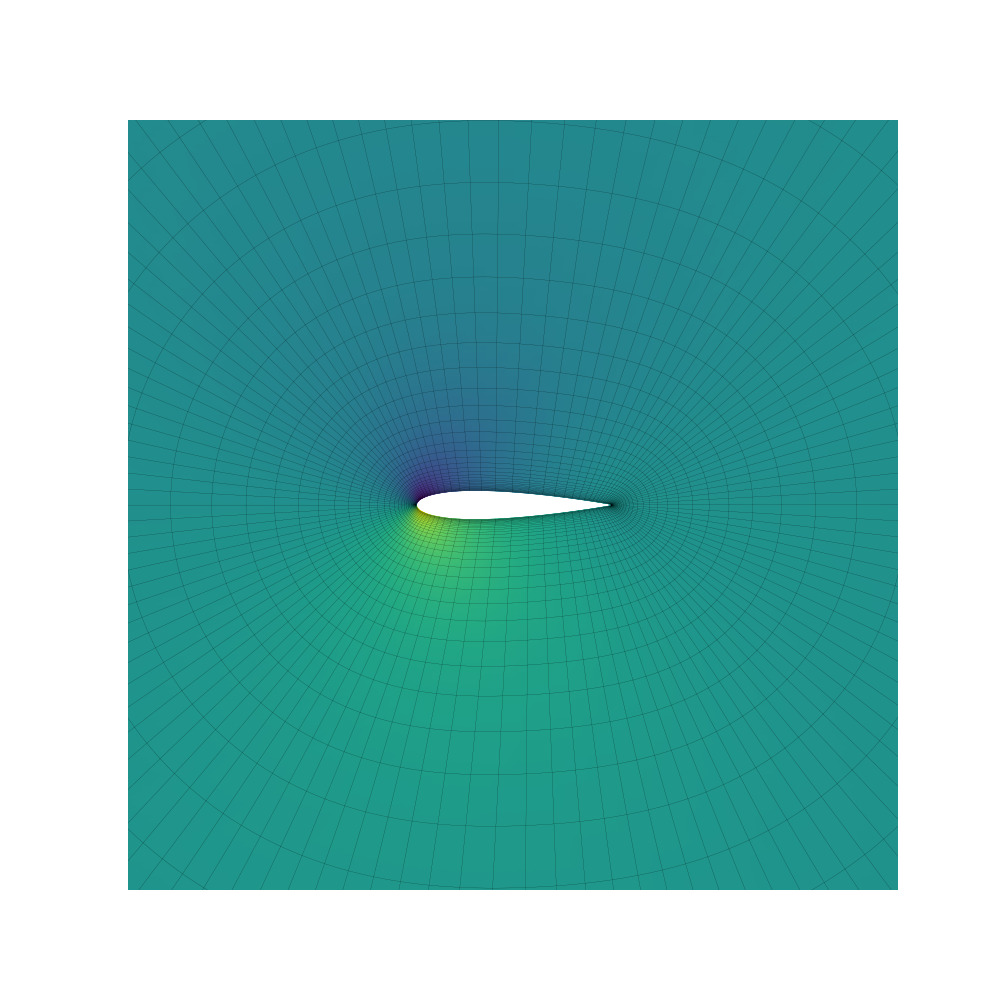}
    }
    \adjustbox{trim={0.3\width} {0.3\height} {0.3\width} {0.3\height}, clip}{
      \includegraphics[width=0.4\textwidth]{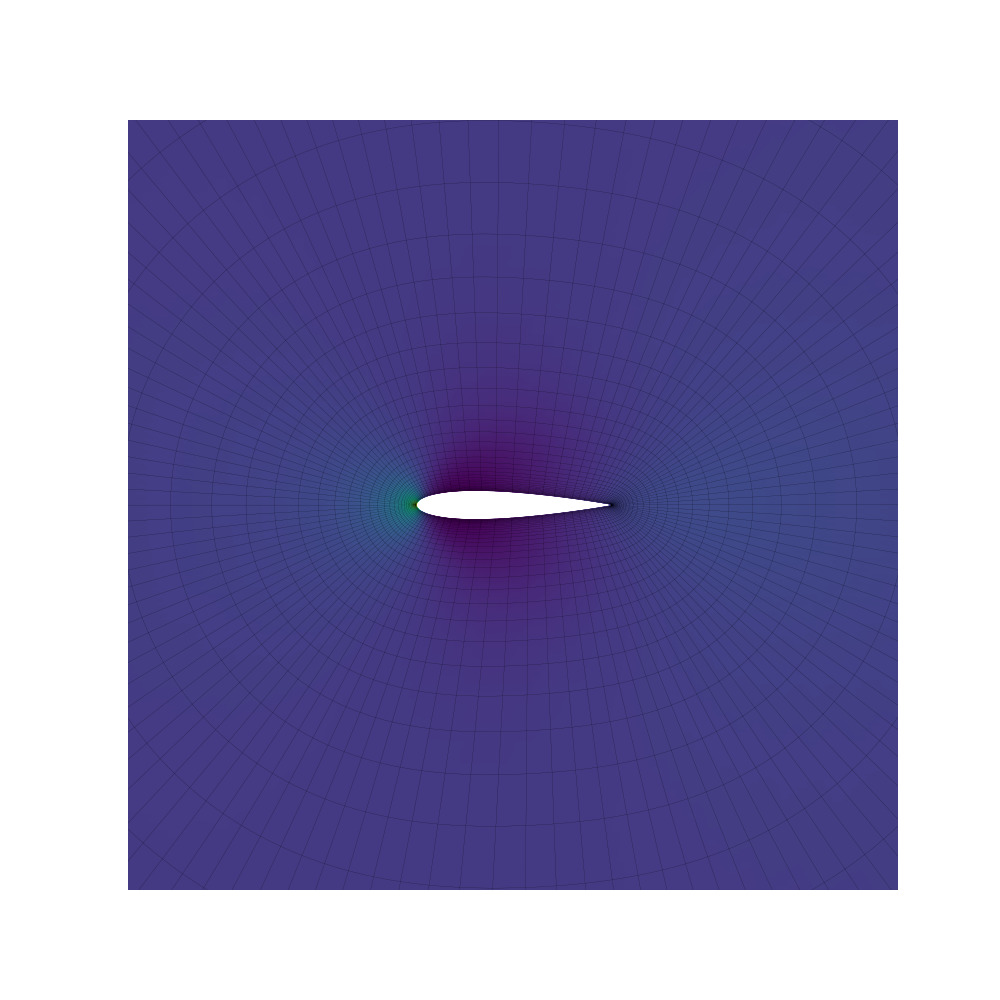}
    } \\
    \adjustbox{trim={0.3\width} {0.3\height} {0.3\width} {0.3\height}, clip}{
      \includegraphics[width=0.4\textwidth]{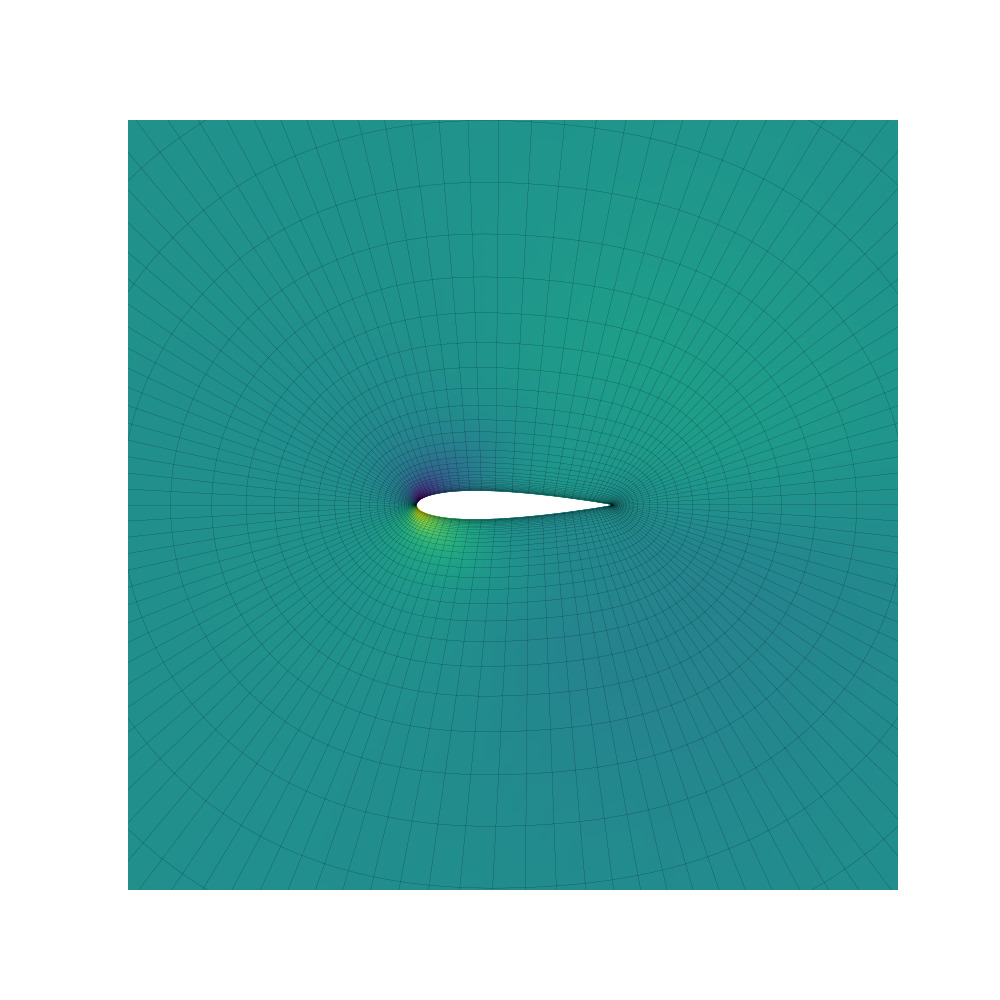}
    }
    \adjustbox{trim={0.3\width} {0.3\height} {0.3\width} {0.3\height}, clip}{
      \includegraphics[width=0.4\textwidth]{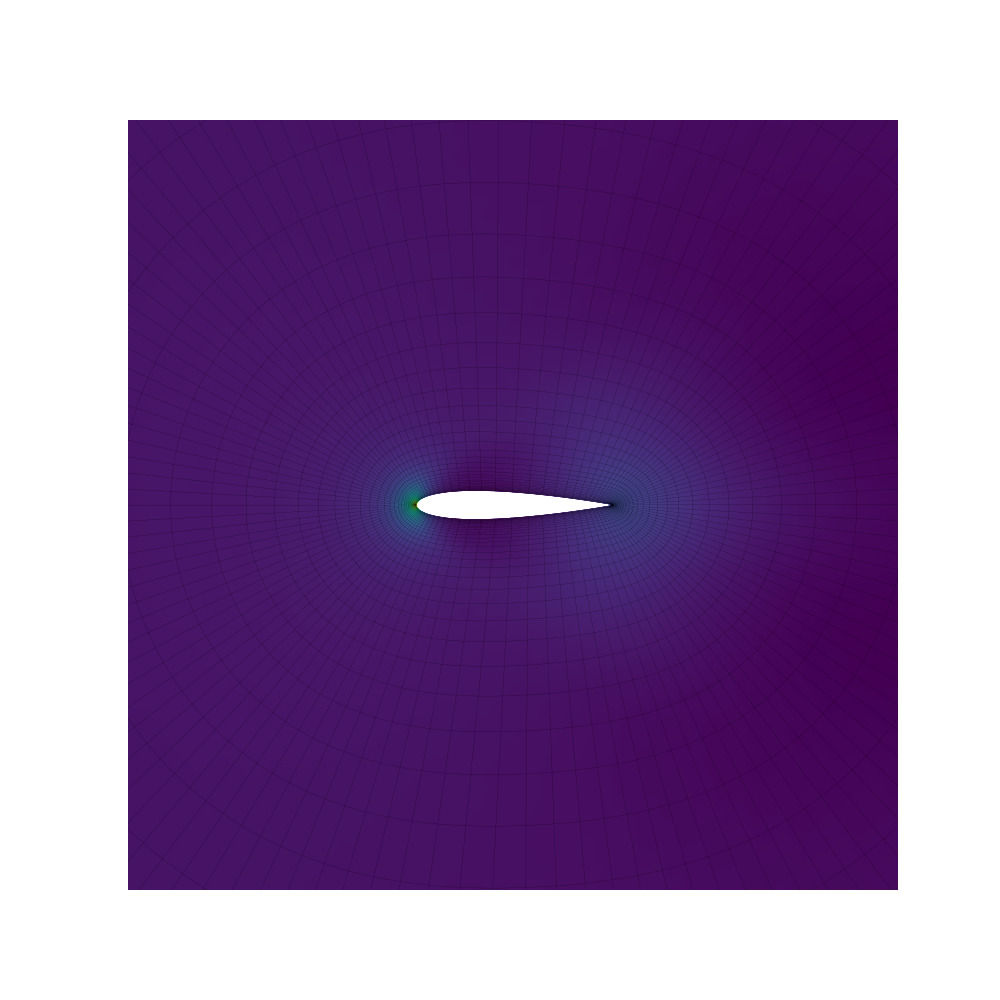}
    }
    \adjustbox{trim={0.3\width} {0.3\height} {0.3\width} {0.3\height}, clip}{
      \includegraphics[width=0.4\textwidth]{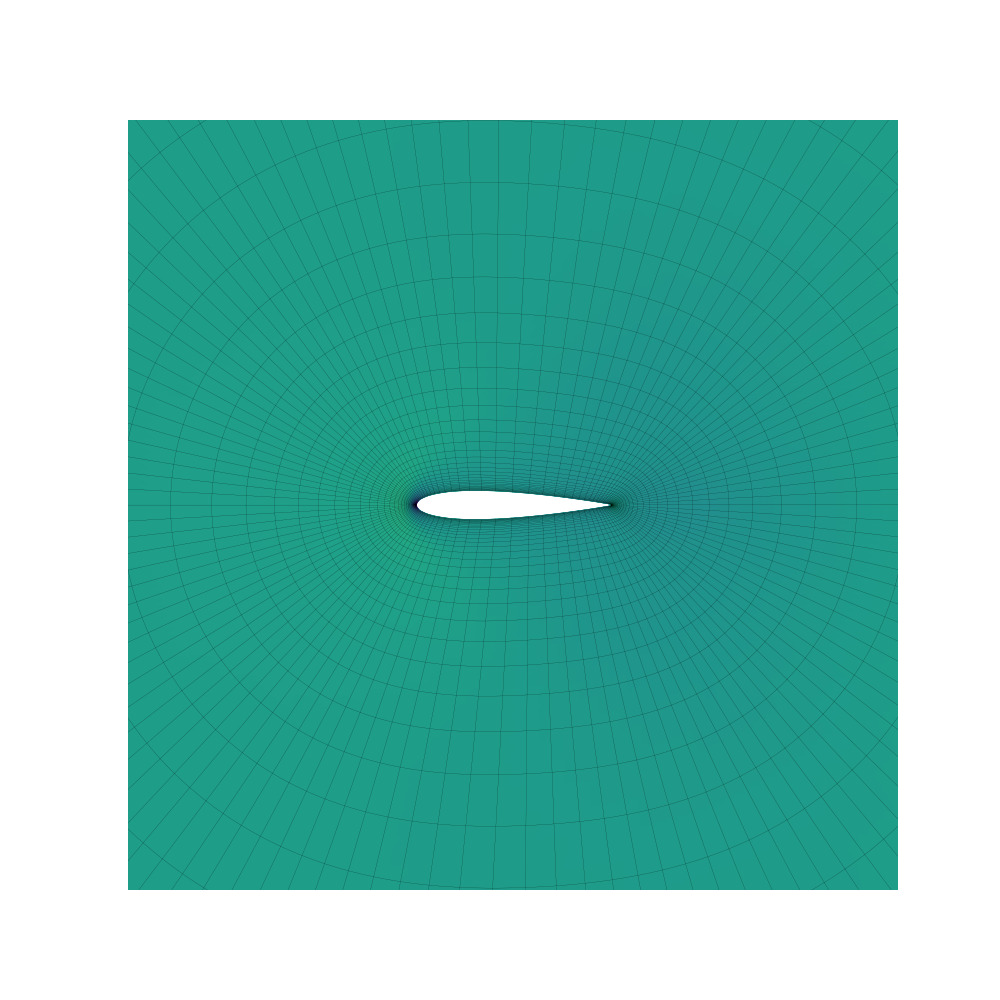}
    }
    \caption{First six pressure modes (Taylor-Hood method).}
    \label{fig:modes-press-noconf}
    \vspace{\floatsep}
    \adjustbox{trim={0.24\width} {0.23\height} {0.22\width} {0.24\height}, clip}{
      \includegraphics[width=0.4\textwidth]{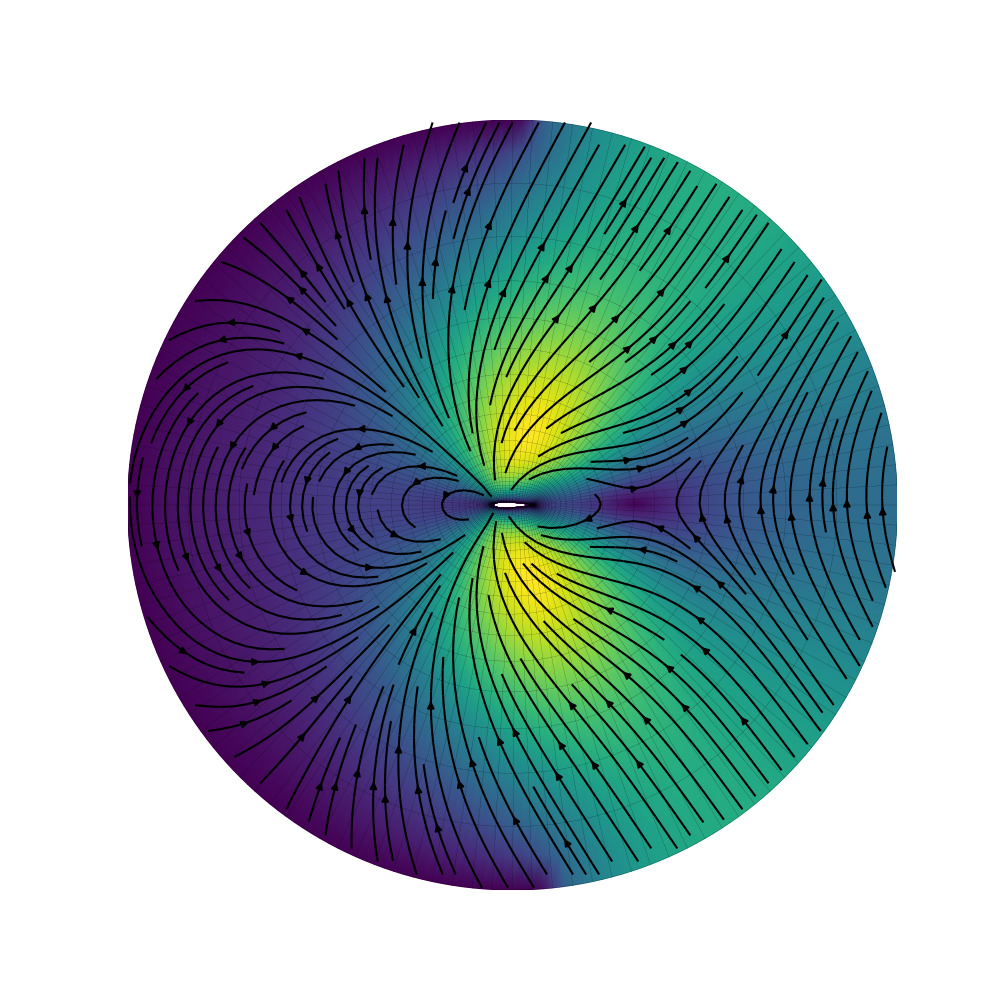}
    }
    \adjustbox{trim={0.24\width} {0.23\height} {0.22\width} {0.24\height}, clip}{
      \includegraphics[width=0.4\textwidth]{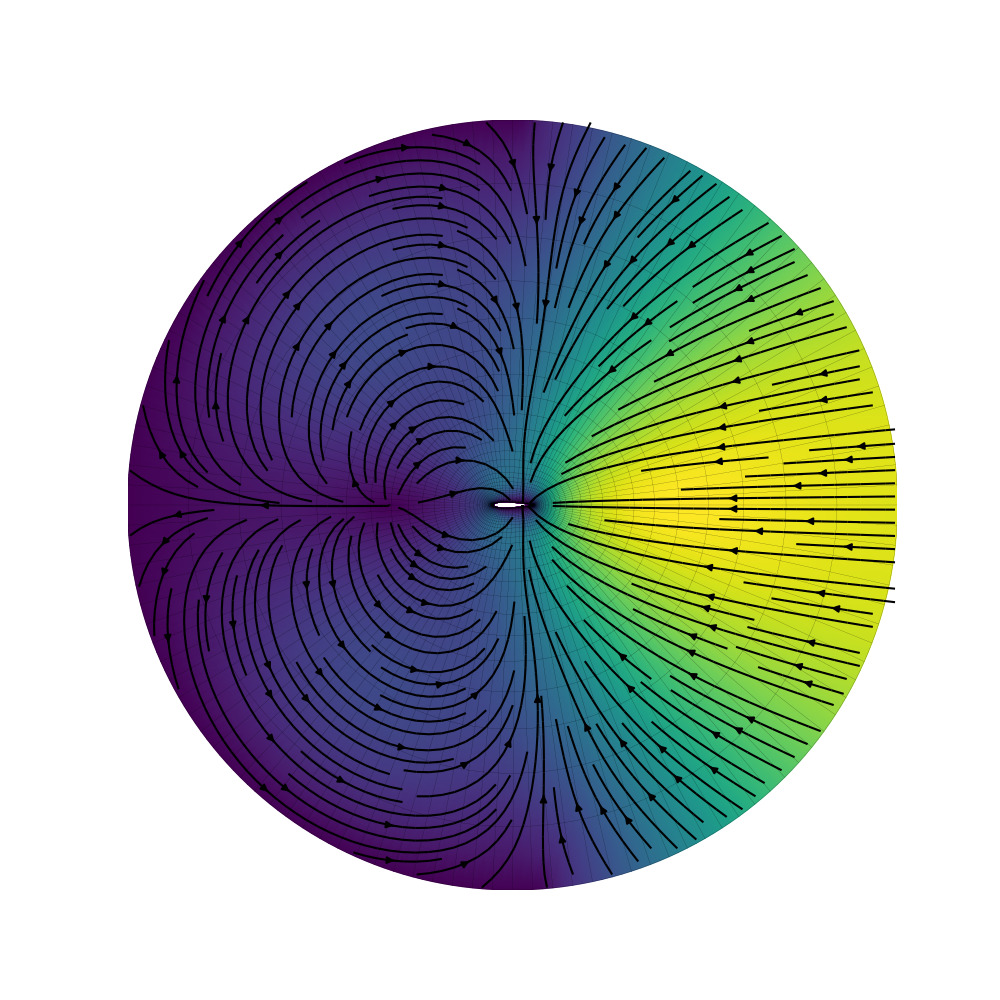}
    }
    \adjustbox{trim={0.24\width} {0.23\height} {0.22\width} {0.24\height}, clip}{
      \includegraphics[width=0.4\textwidth]{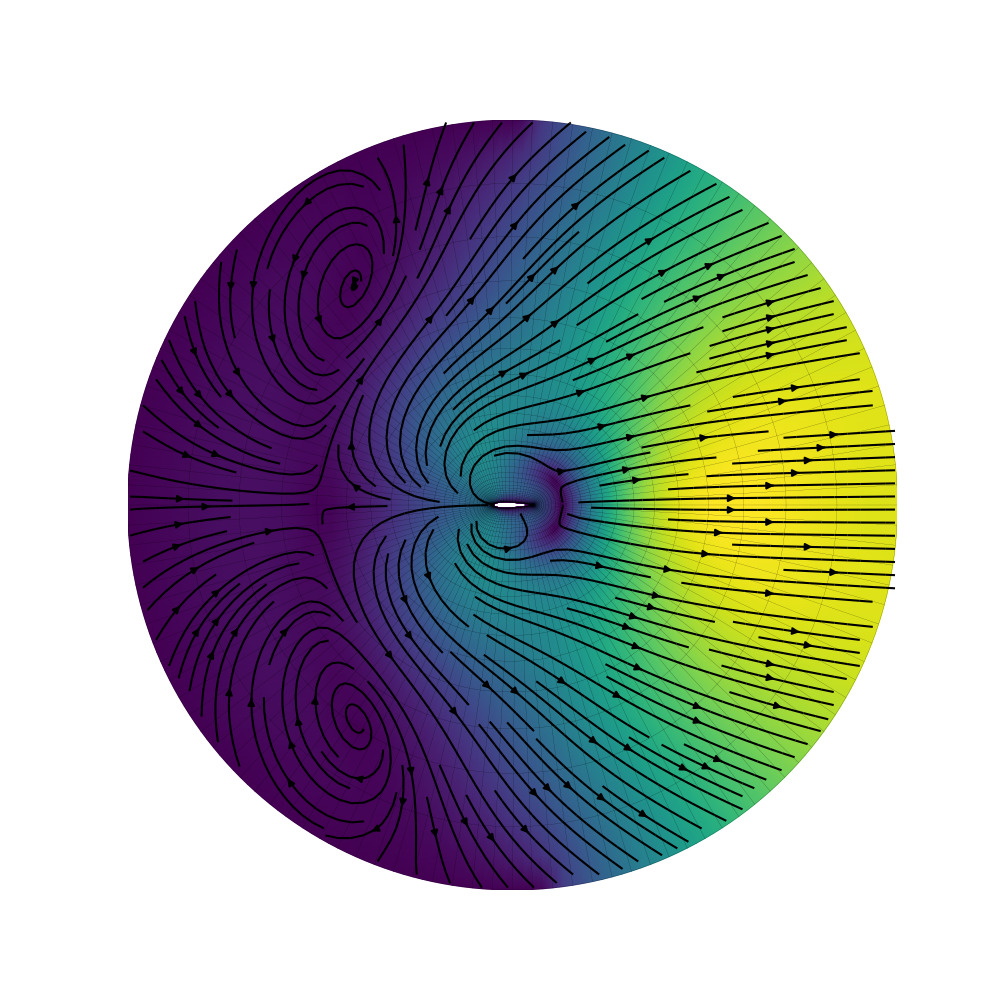}
    } \\
    \adjustbox{trim={0.24\width} {0.23\height} {0.22\width} {0.24\height}, clip}{
      \includegraphics[width=0.4\textwidth]{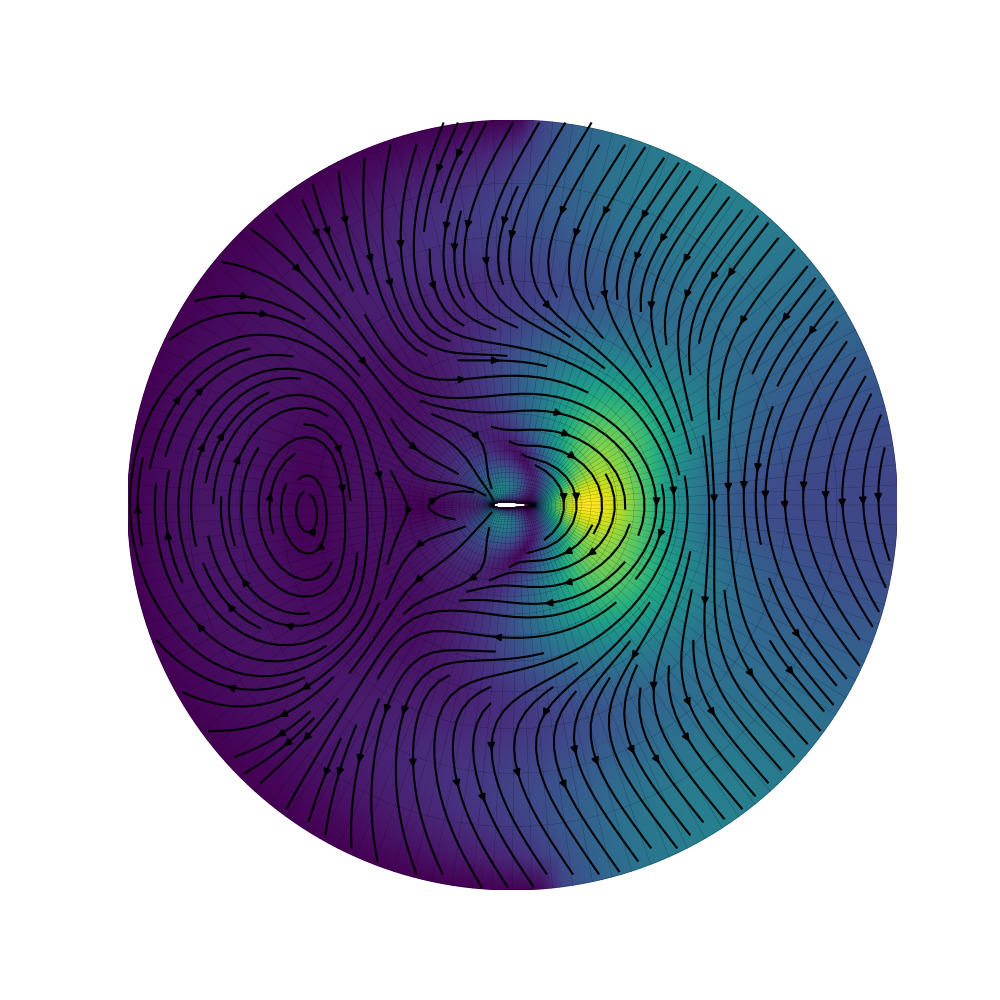}
    }
    \adjustbox{trim={0.24\width} {0.23\height} {0.22\width} {0.24\height}, clip}{
      \includegraphics[width=0.4\textwidth]{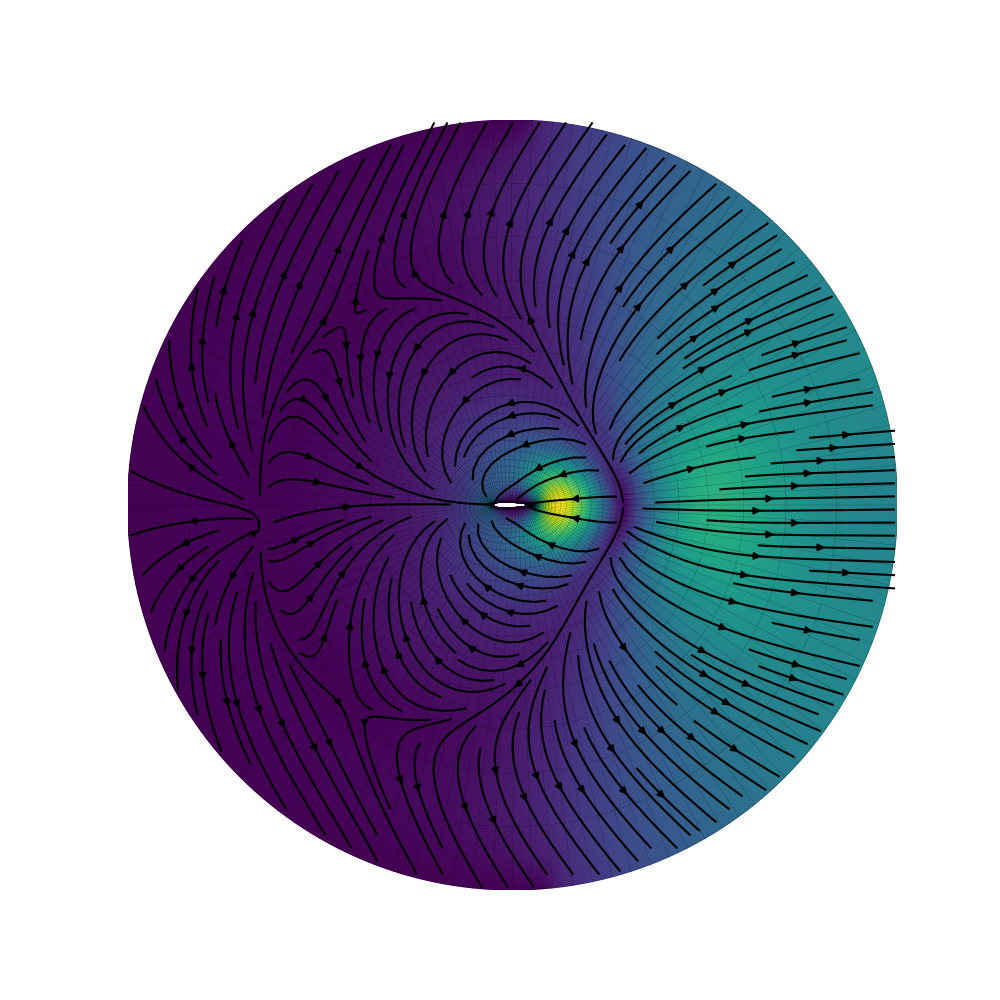}
    }
    \adjustbox{trim={0.24\width} {0.23\height} {0.22\width} {0.24\height}, clip}{
      \includegraphics[width=0.4\textwidth]{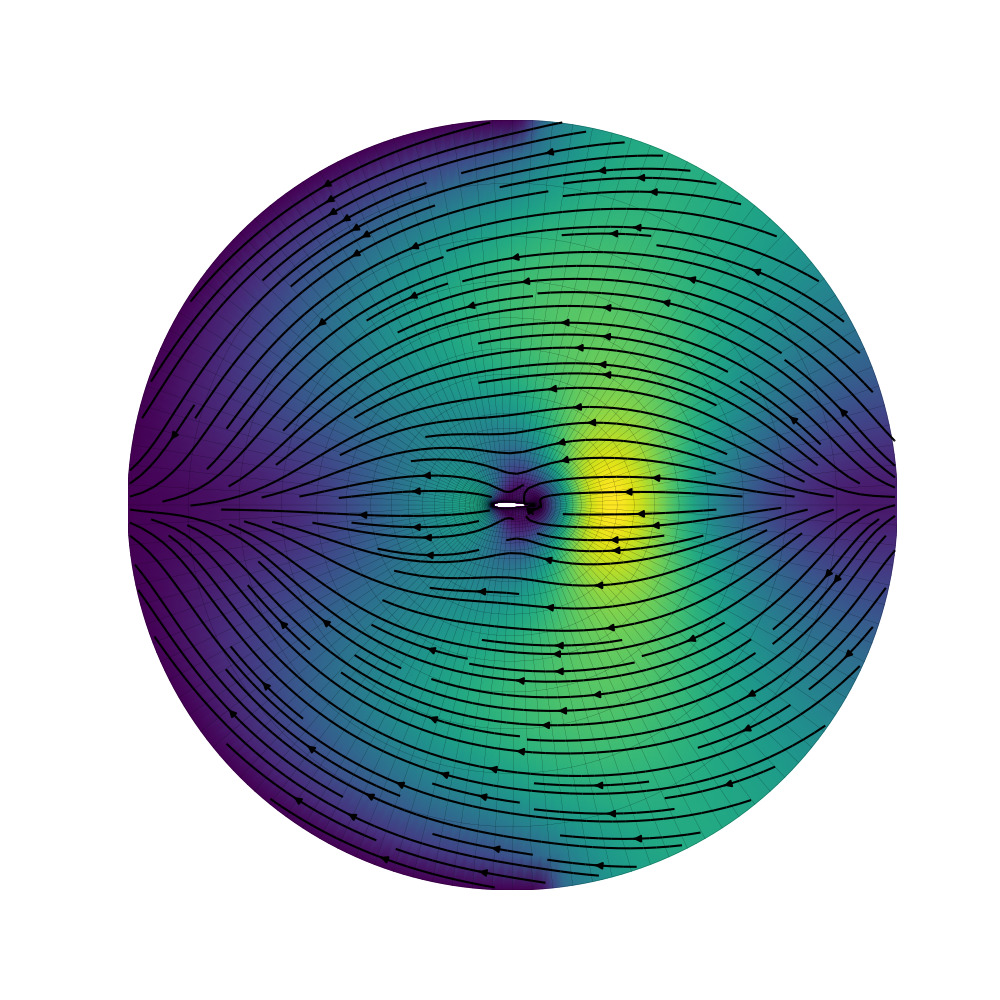}
    }
    \caption{First six supremizer modes (Taylor-Hood method).}
    \label{fig:modes-sup-noconf}
  \end{center}
\end{figure}

\begin{figure}
  \begin{center}
    \adjustbox{trim={0.13\width} {0.12\height} {0.13\width} {0.12\height}, clip}{
      \includegraphics[width=0.4\textwidth]{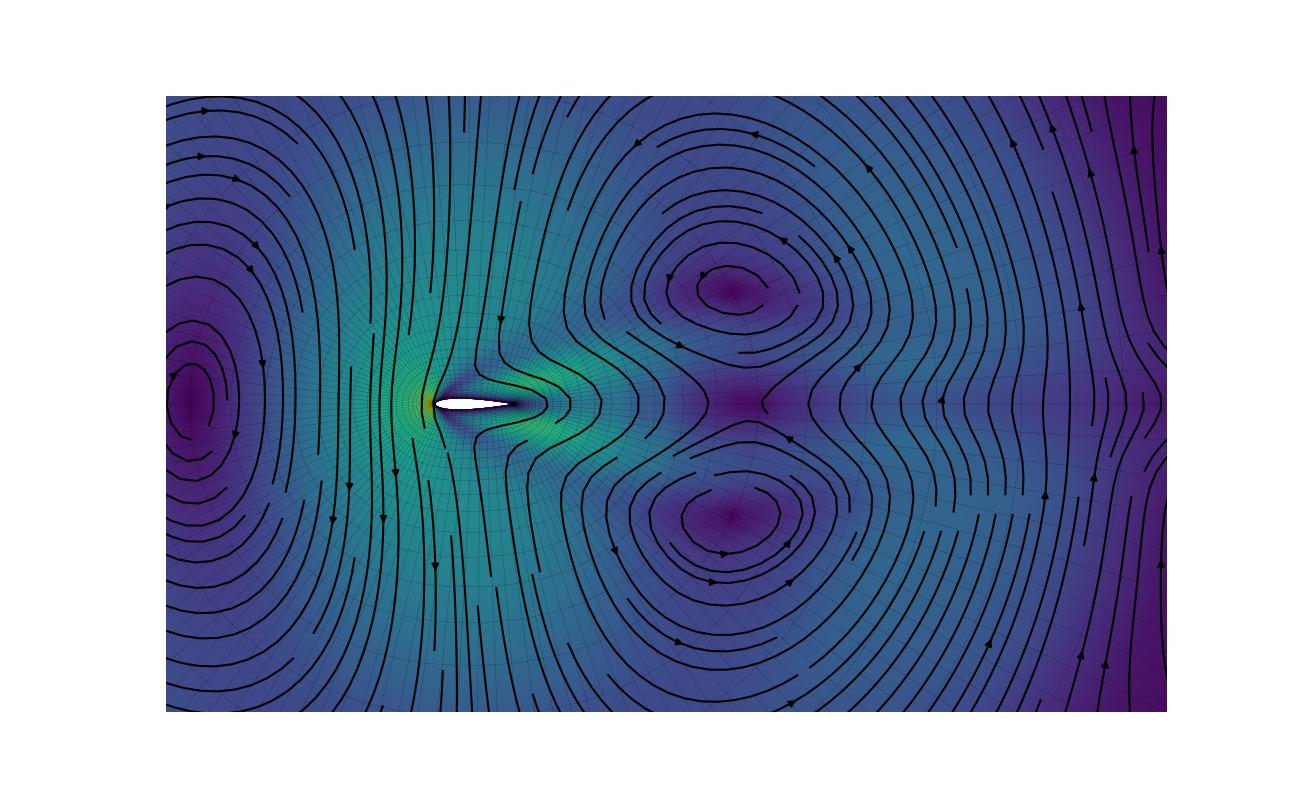}
    }
    \adjustbox{trim={0.13\width} {0.12\height} {0.13\width} {0.12\height}, clip}{
      \includegraphics[width=0.4\textwidth]{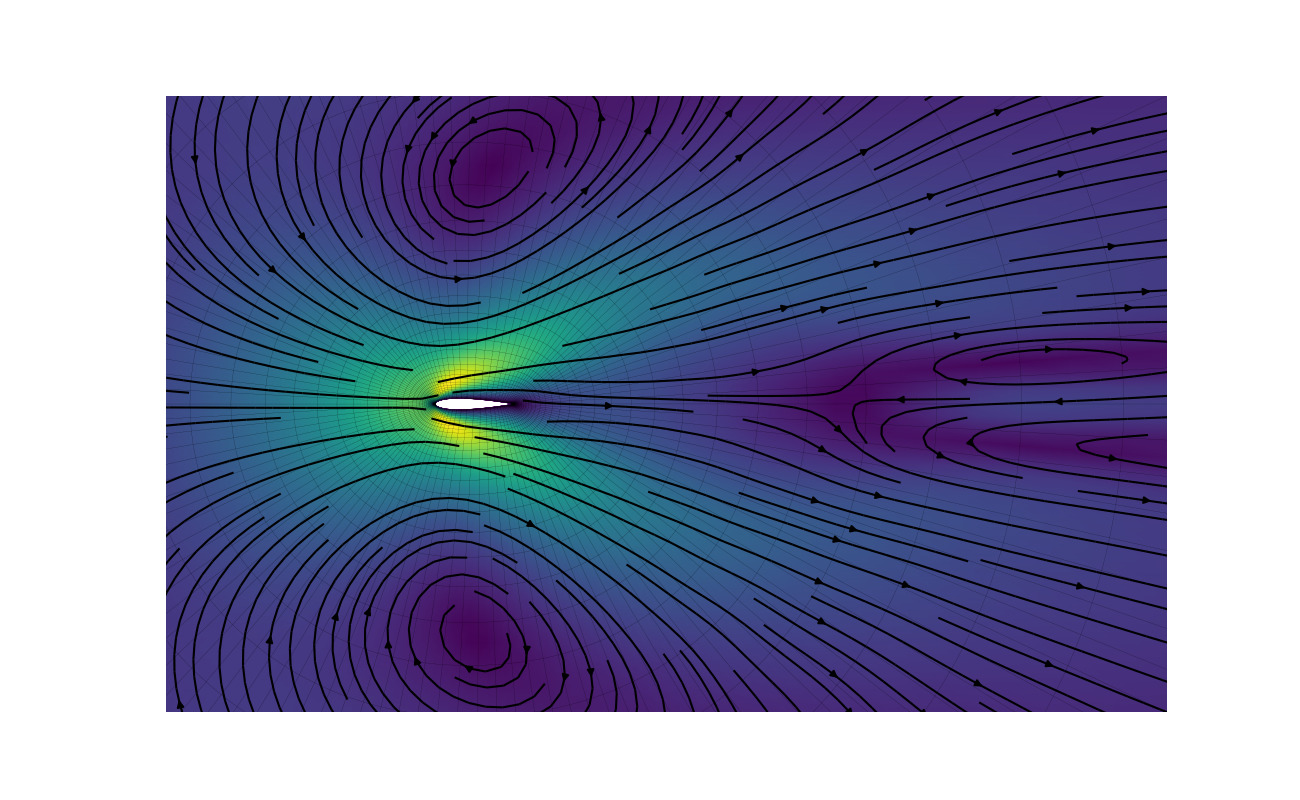}
    }
    \adjustbox{trim={0.13\width} {0.12\height} {0.13\width} {0.12\height}, clip}{
      \includegraphics[width=0.4\textwidth]{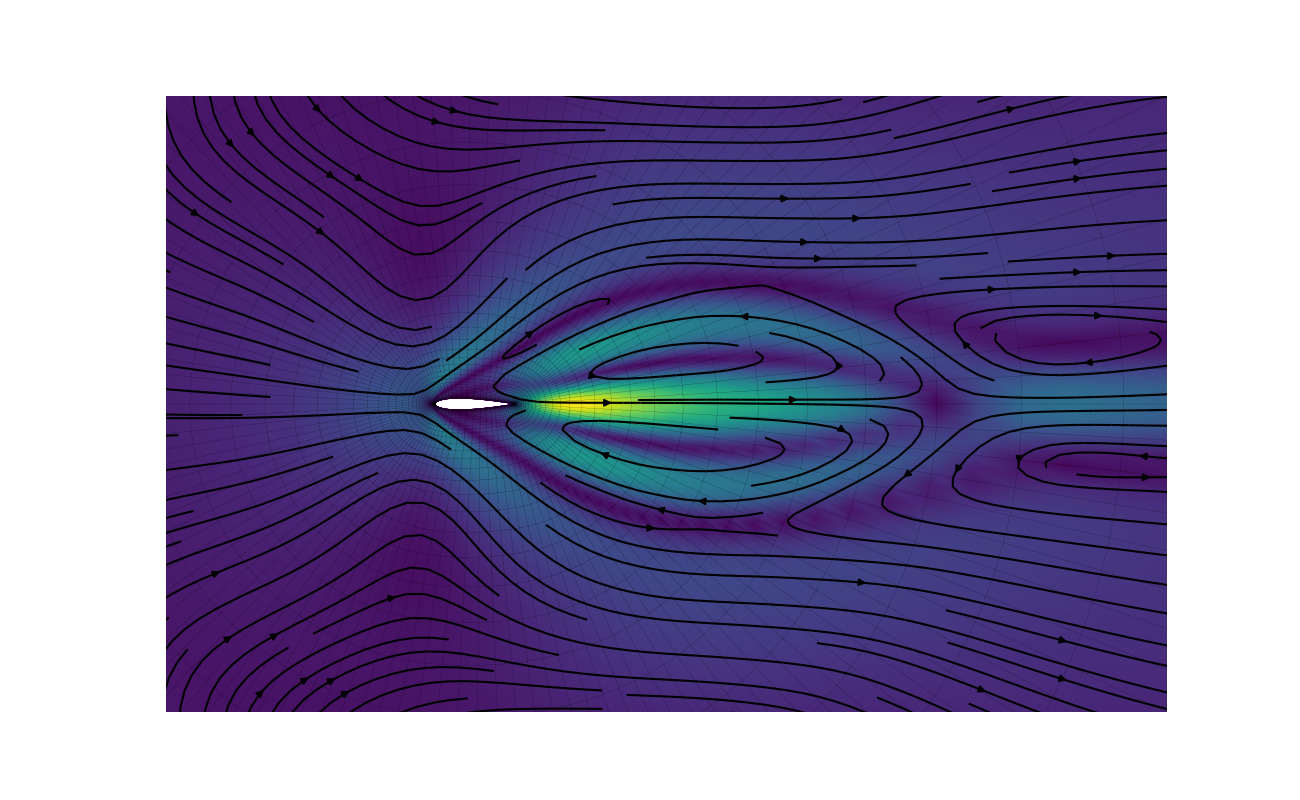}
    } \\
    \adjustbox{trim={0.13\width} {0.12\height} {0.13\width} {0.12\height}, clip}{
      \includegraphics[width=0.4\textwidth]{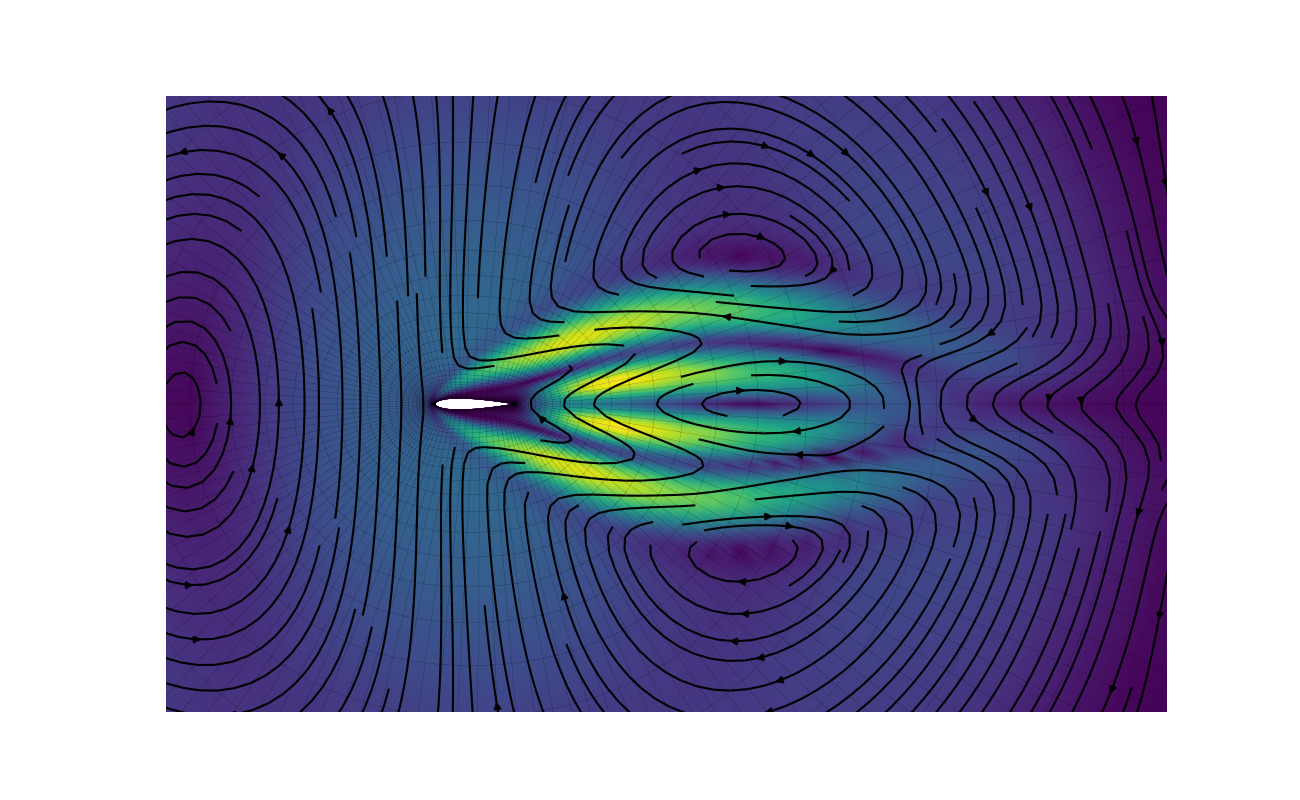}
    }
    \adjustbox{trim={0.13\width} {0.12\height} {0.13\width} {0.12\height}, clip}{
      \includegraphics[width=0.4\textwidth]{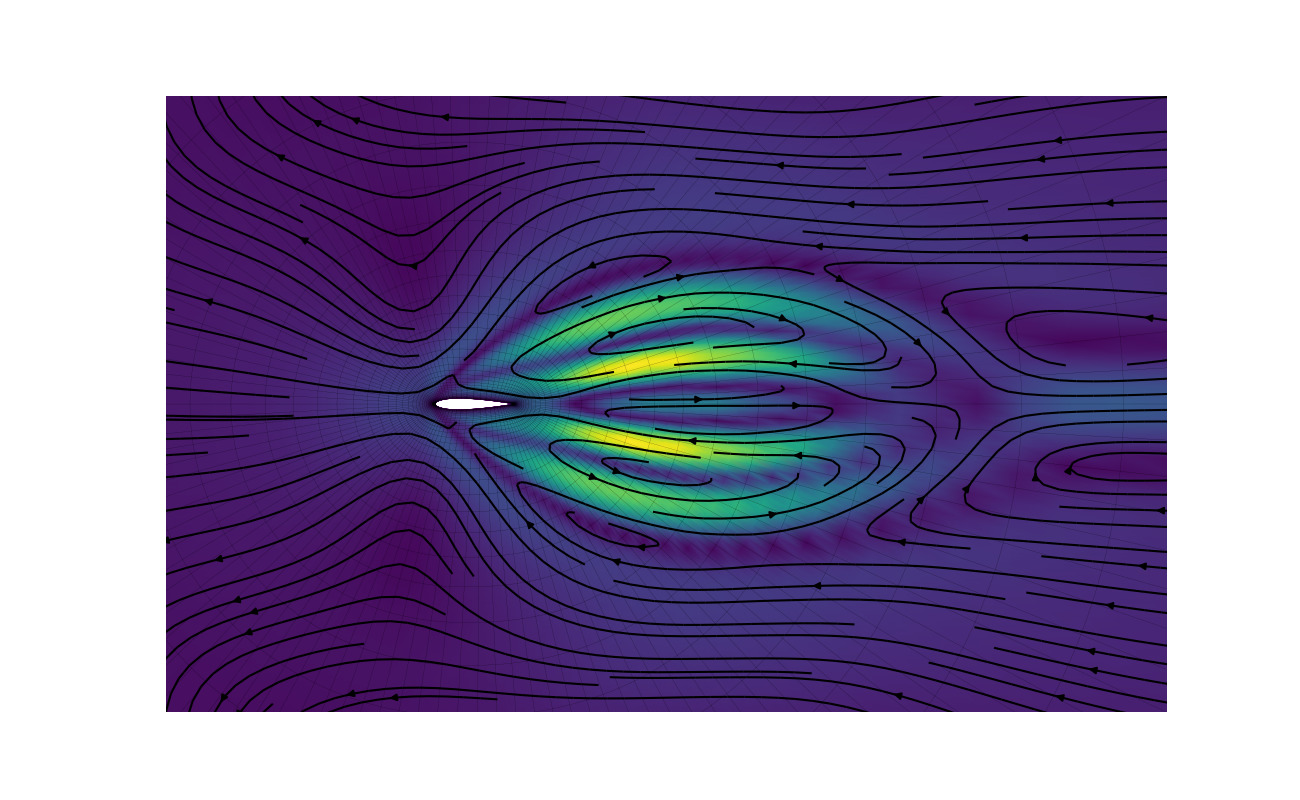}
    }
    \adjustbox{trim={0.13\width} {0.12\height} {0.13\width} {0.12\height}, clip}{
      \includegraphics[width=0.4\textwidth]{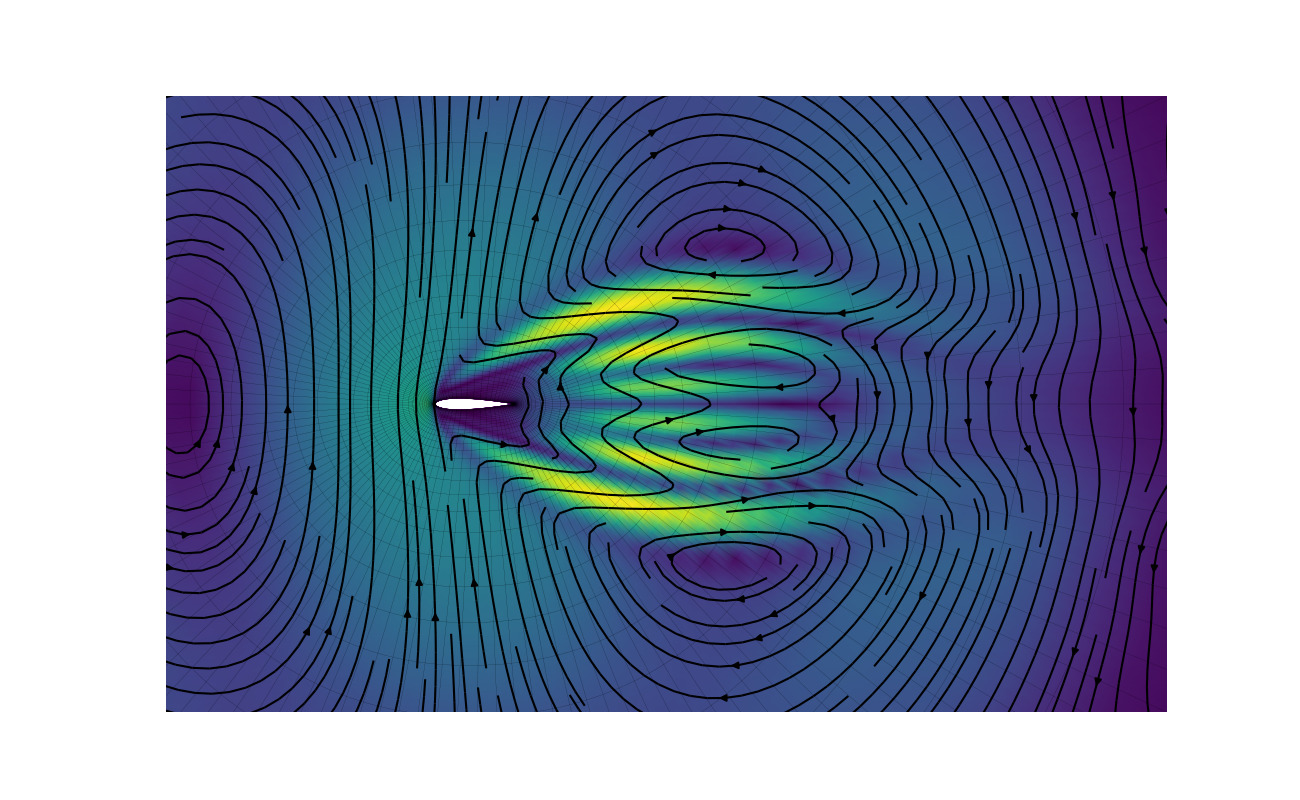}
    }
    \caption{First six velocity modes (conforming method).}
    \label{fig:modes-vel-conf}
    \vspace{\floatsep}
    \adjustbox{trim={0.3\width} {0.3\height} {0.3\width} {0.3\height}, clip}{
      \includegraphics[width=0.4\textwidth]{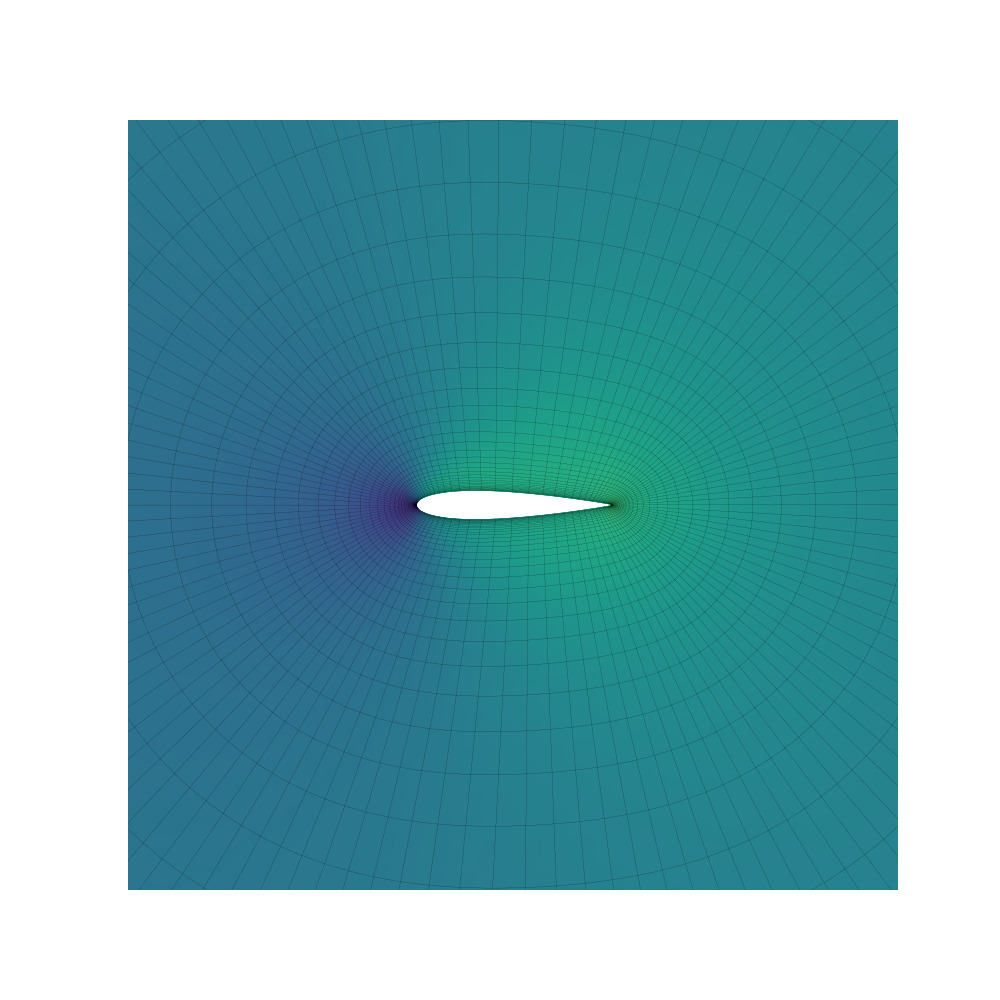}
    }
    \adjustbox{trim={0.3\width} {0.3\height} {0.3\width} {0.3\height}, clip}{
      \includegraphics[width=0.4\textwidth]{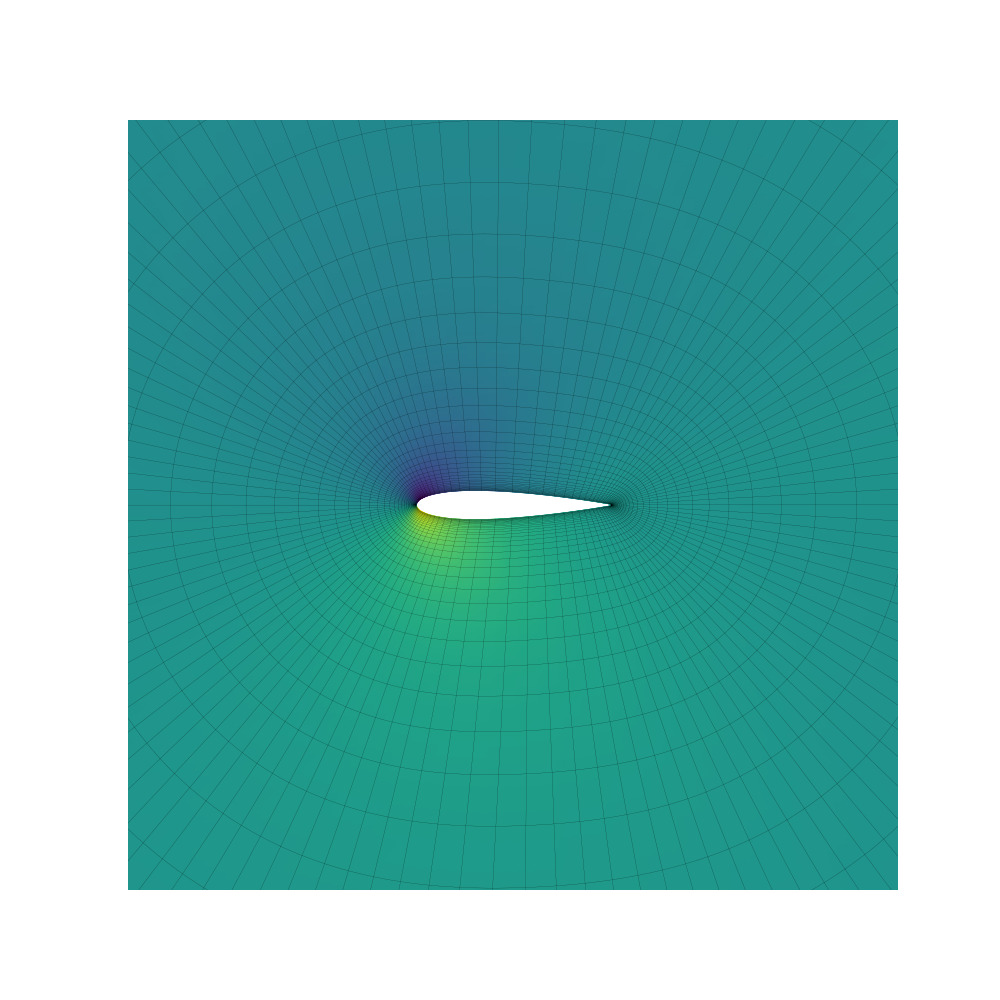}
    }
    \adjustbox{trim={0.3\width} {0.3\height} {0.3\width} {0.3\height}, clip}{
      \includegraphics[width=0.4\textwidth]{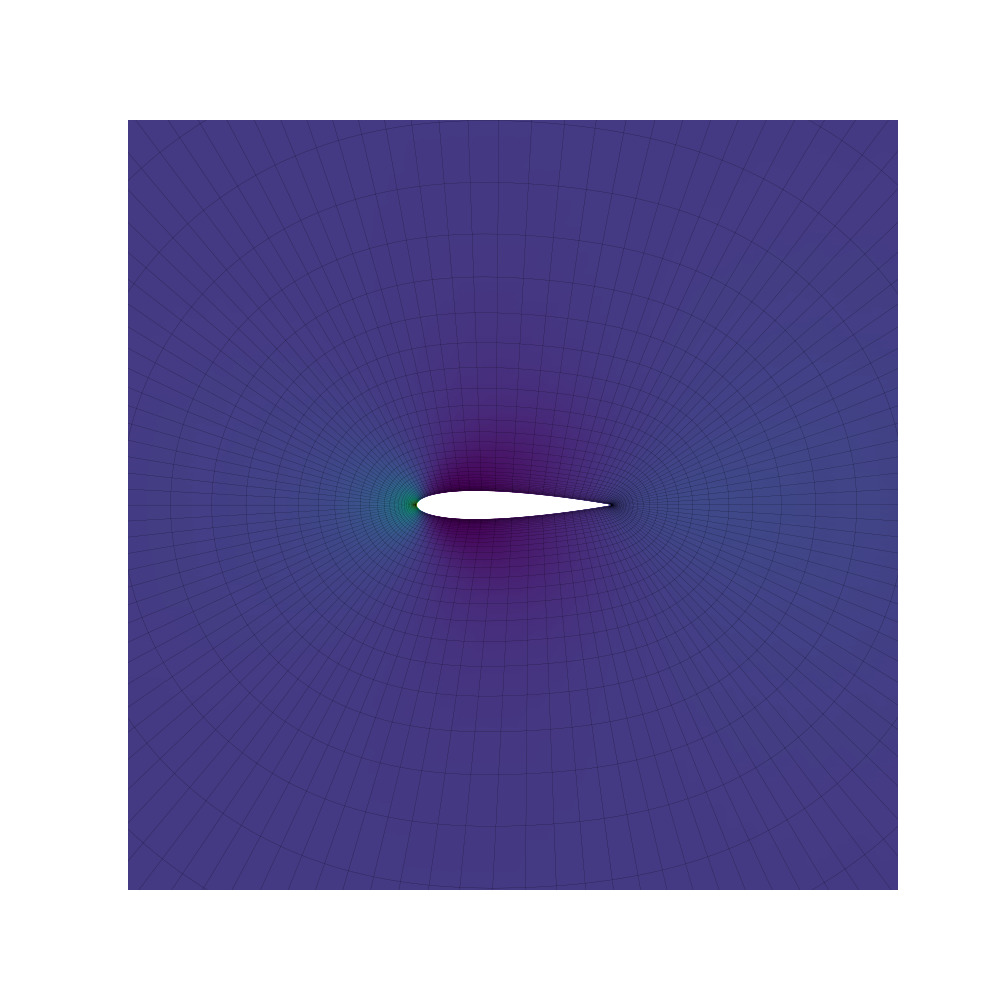}
    } \\
    \adjustbox{trim={0.3\width} {0.3\height} {0.3\width} {0.3\height}, clip}{
      \includegraphics[width=0.4\textwidth]{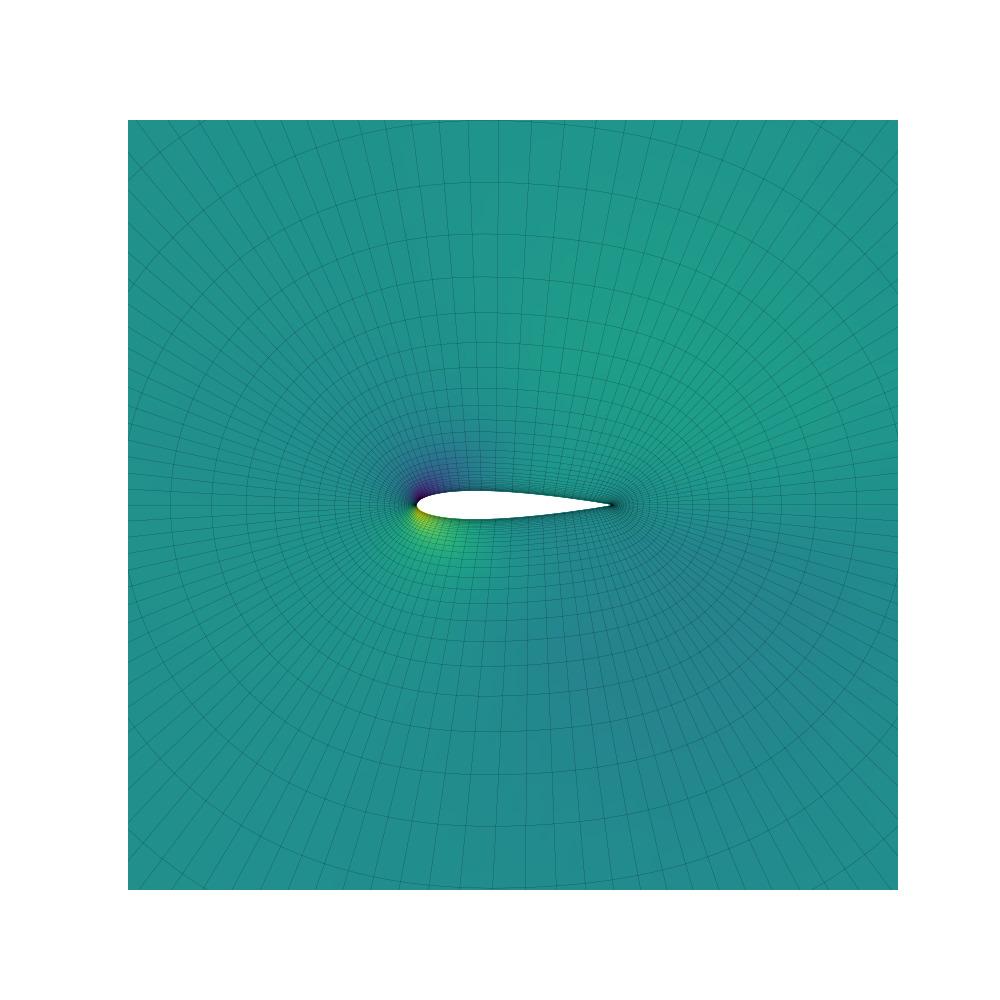}
    }
    \adjustbox{trim={0.3\width} {0.3\height} {0.3\width} {0.3\height}, clip}{
      \includegraphics[width=0.4\textwidth]{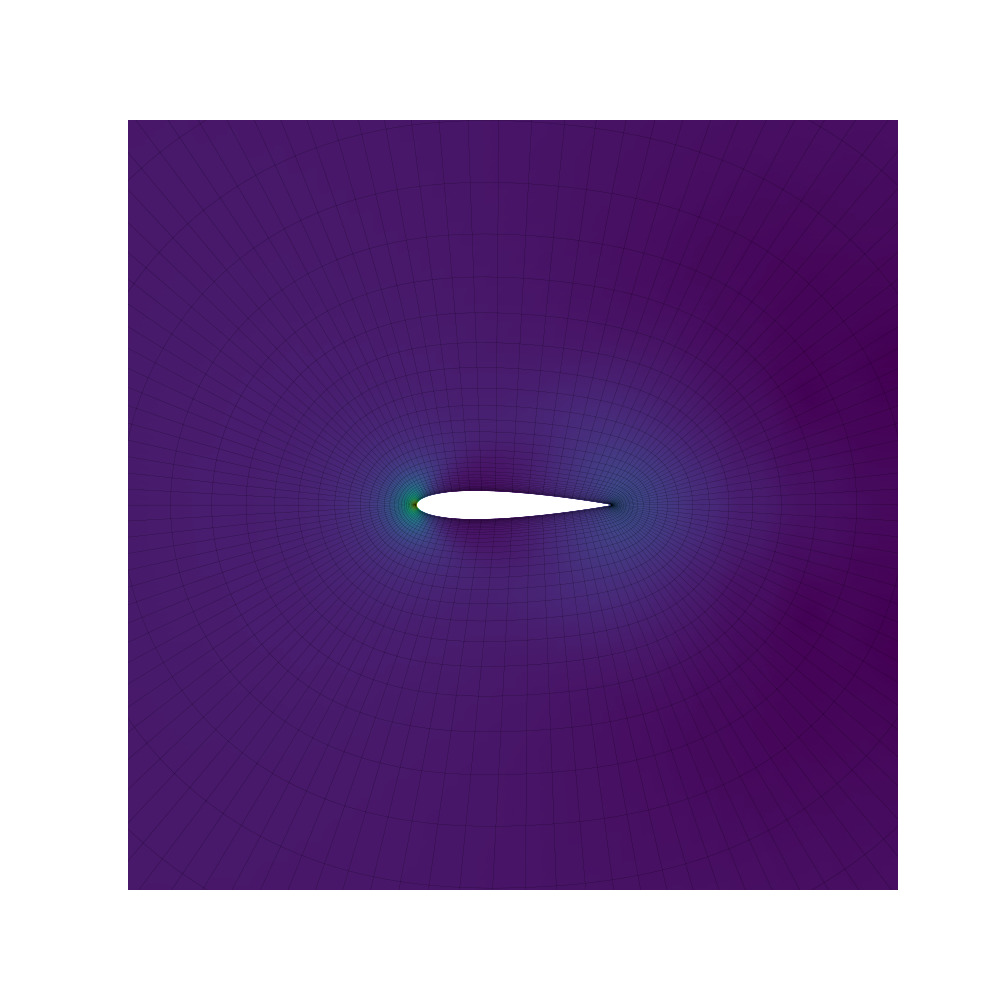}
    }
    \adjustbox{trim={0.3\width} {0.3\height} {0.3\width} {0.3\height}, clip}{
      \includegraphics[width=0.4\textwidth]{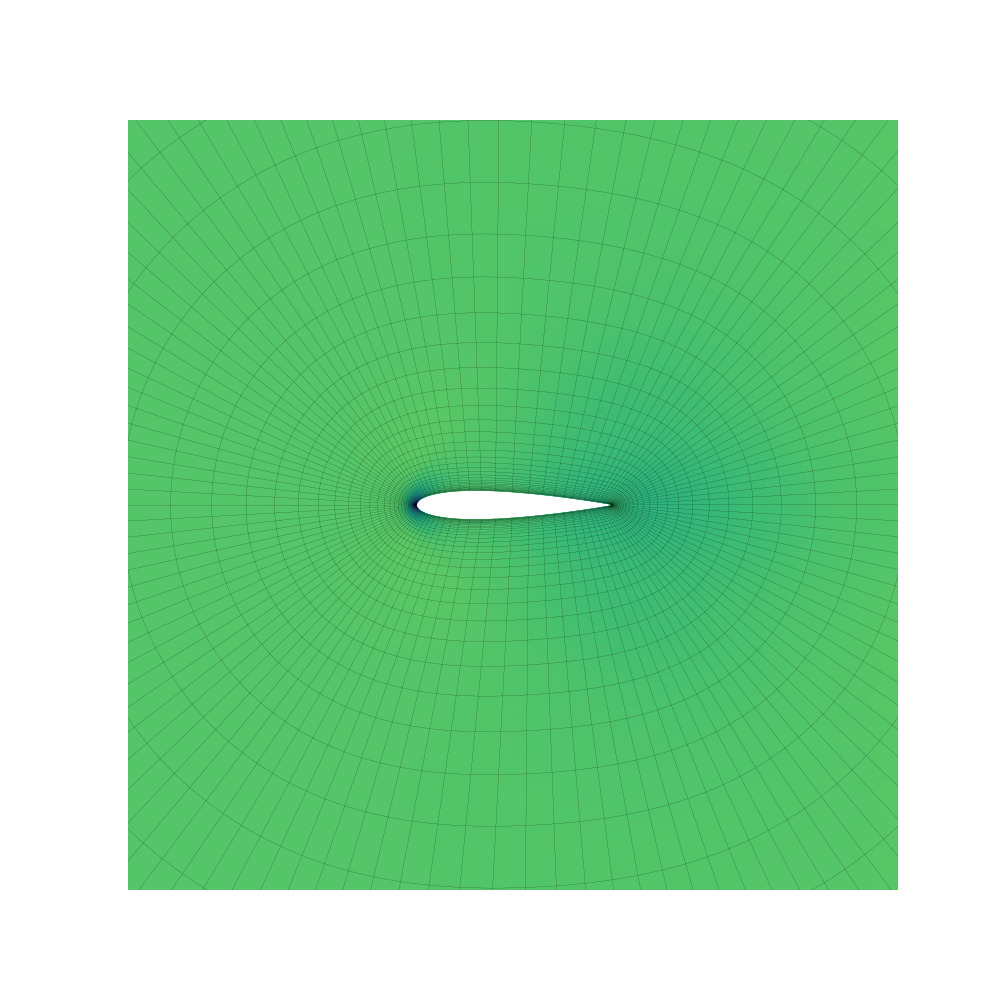}
    }
    \caption{First six pressure modes (conforming method).}
    \label{fig:modes-press-conf}
    \vspace{\floatsep}
    \adjustbox{trim={0.24\width} {0.23\height} {0.22\width} {0.24\height}, clip}{
      \includegraphics[width=0.4\textwidth]{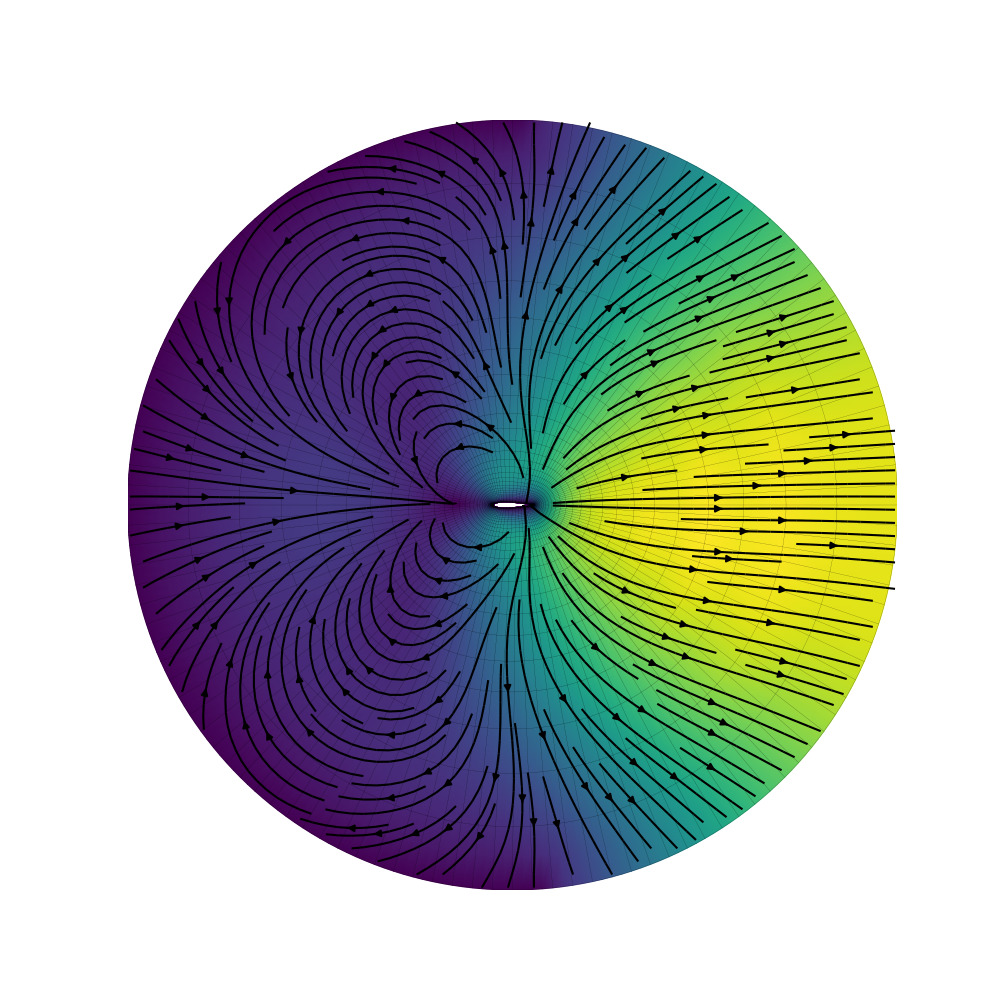}
    }
    \adjustbox{trim={0.24\width} {0.23\height} {0.22\width} {0.24\height}, clip}{
      \includegraphics[width=0.4\textwidth]{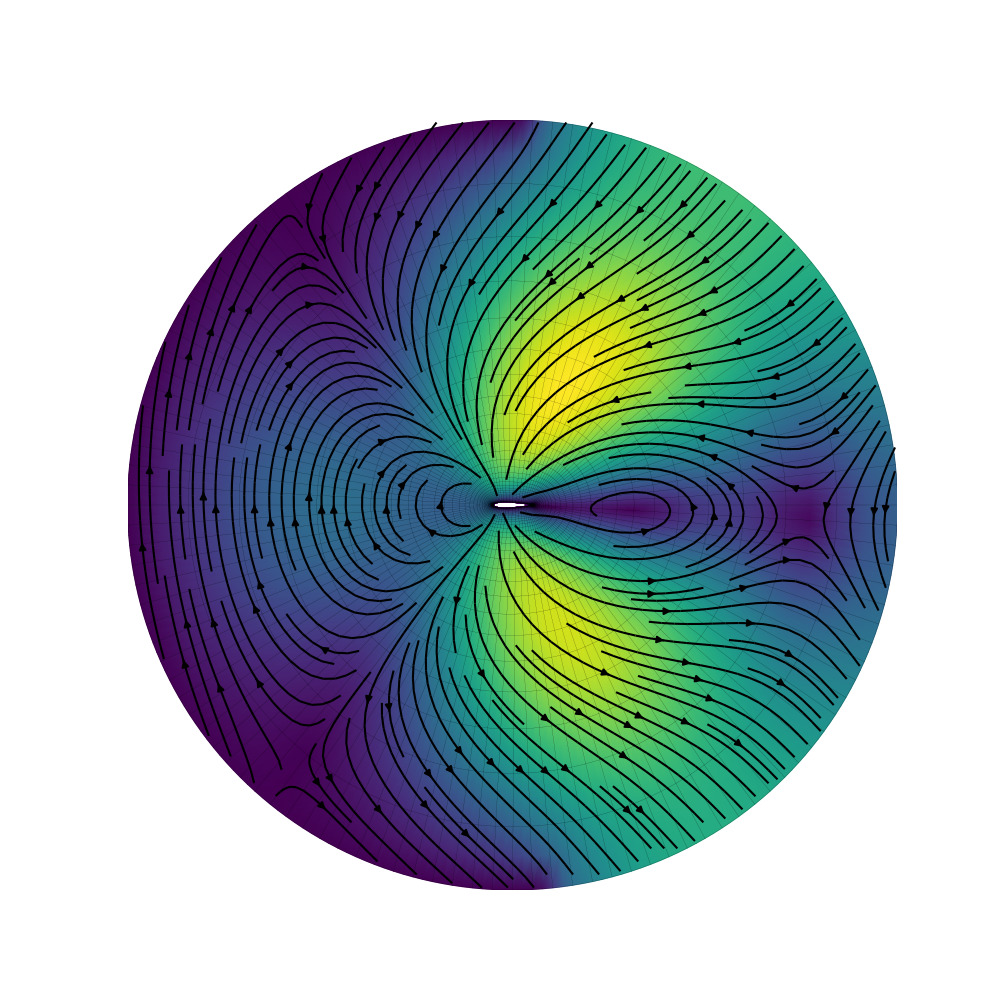}
    }
    \adjustbox{trim={0.24\width} {0.23\height} {0.22\width} {0.24\height}, clip}{
      \includegraphics[width=0.4\textwidth]{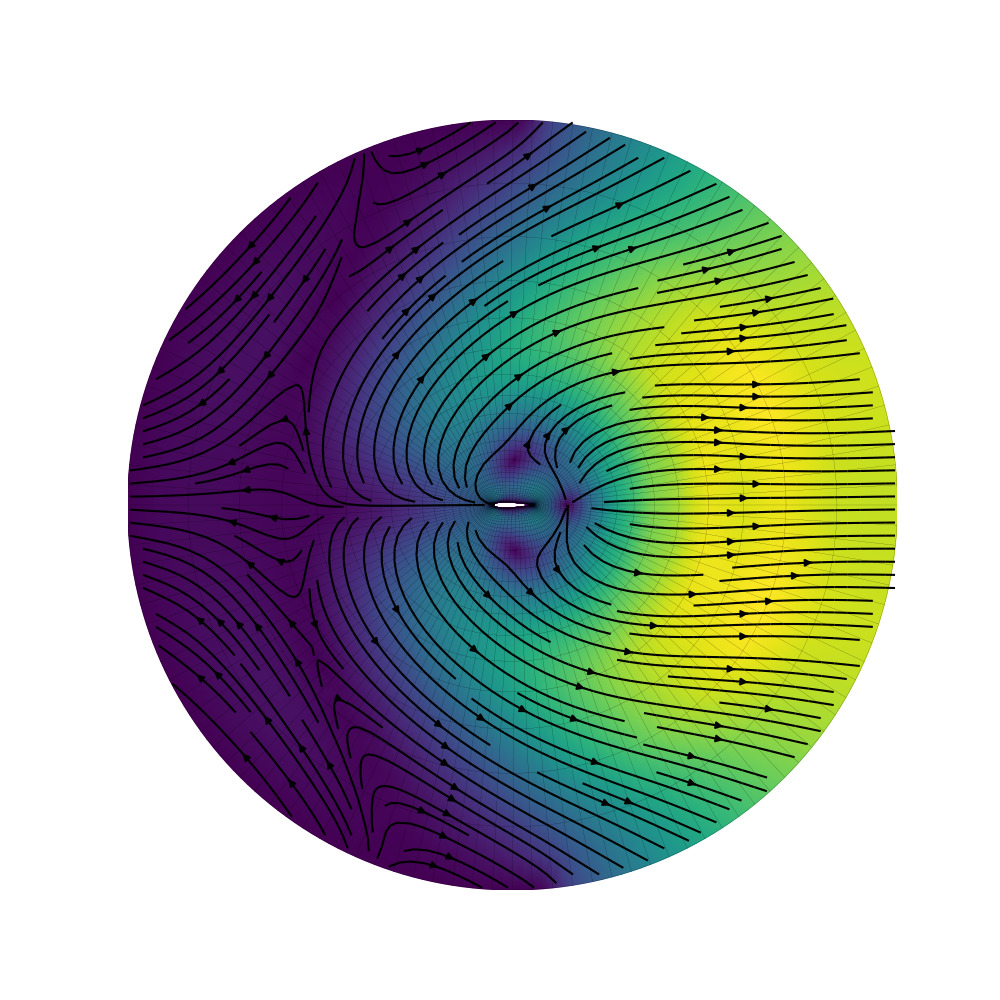}
    } \\
    \adjustbox{trim={0.24\width} {0.23\height} {0.22\width} {0.24\height}, clip}{
      \includegraphics[width=0.4\textwidth]{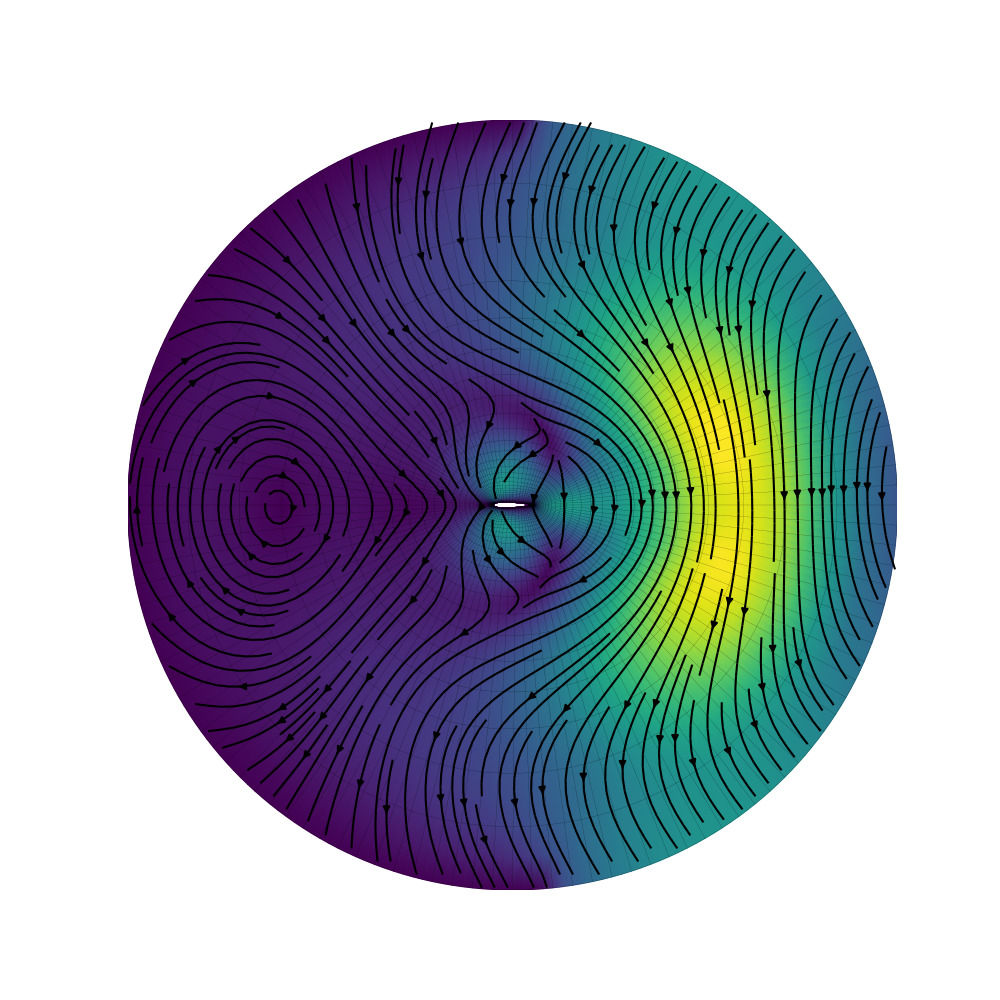}
    }
    \adjustbox{trim={0.24\width} {0.23\height} {0.22\width} {0.24\height}, clip}{
      \includegraphics[width=0.4\textwidth]{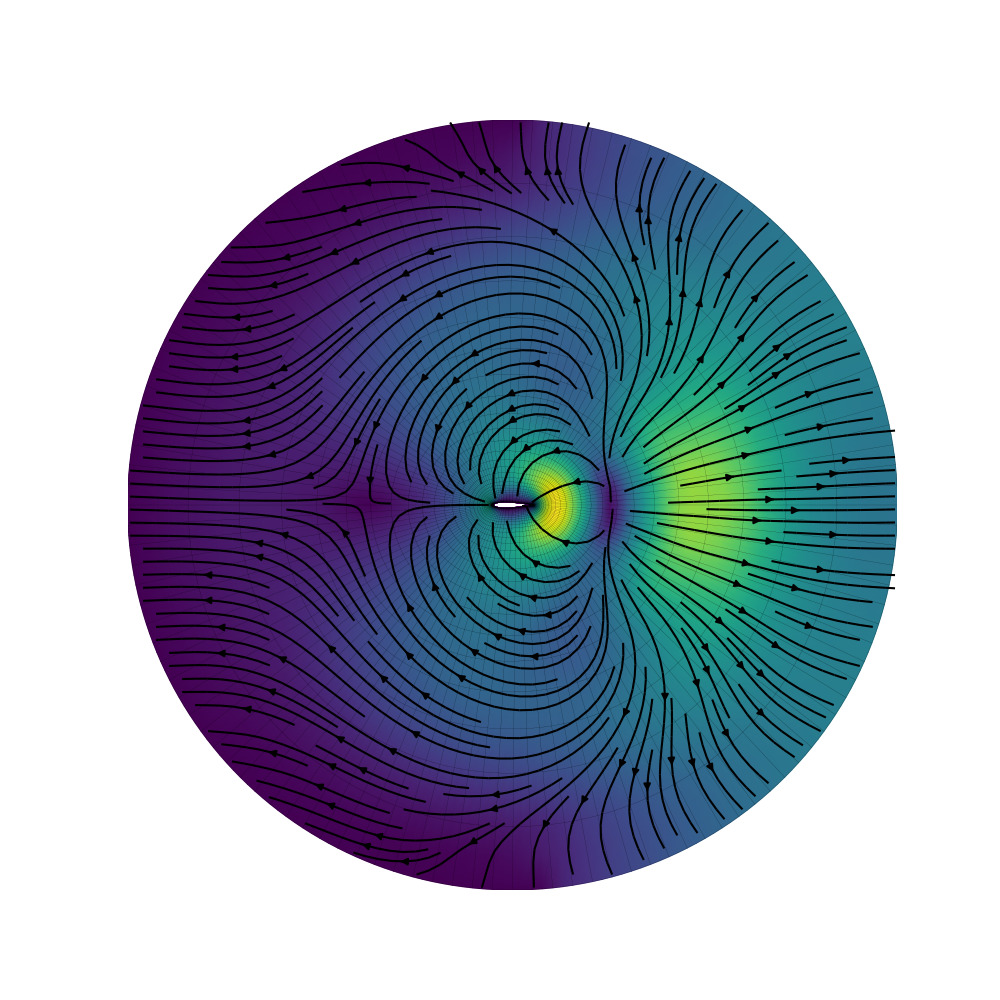}
    }
    \adjustbox{trim={0.24\width} {0.23\height} {0.22\width} {0.24\height}, clip}{
      \includegraphics[width=0.4\textwidth]{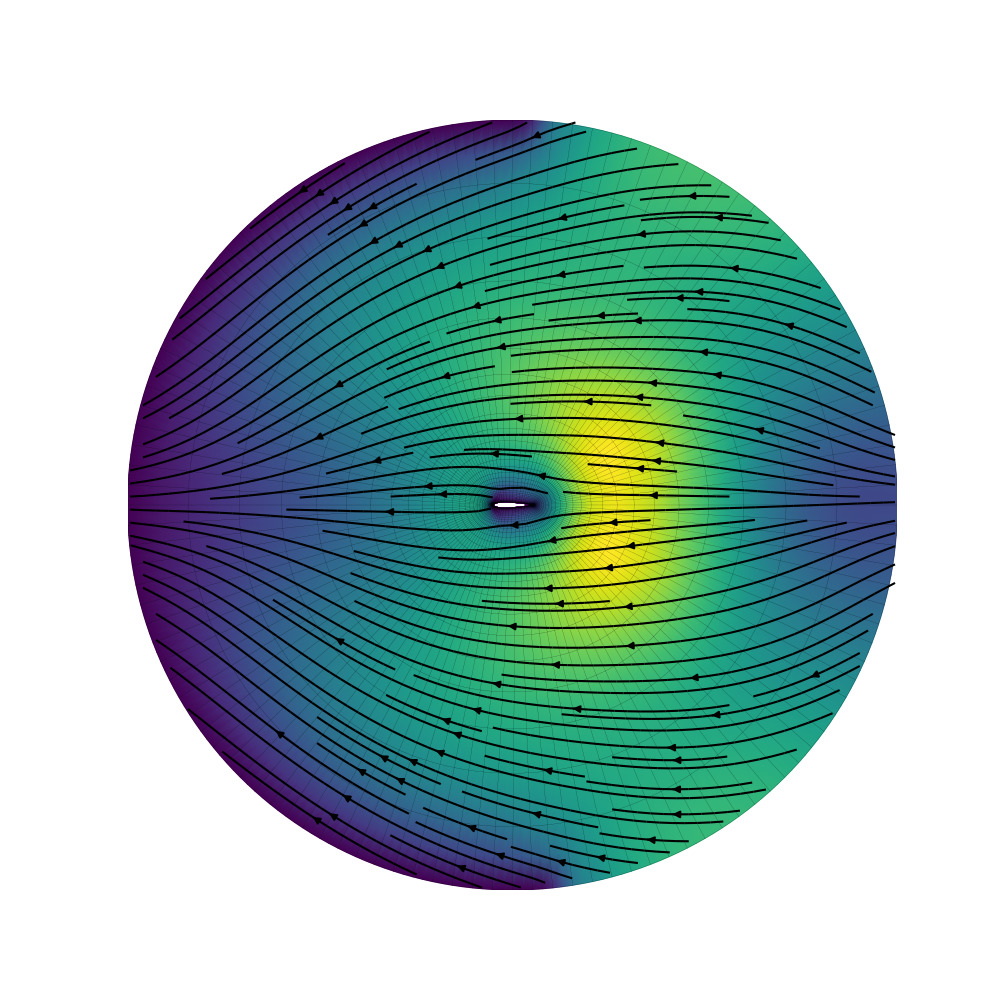}
    }
    \caption{First six supremizer modes (conforming method).}
    \label{fig:modes-sup-conf}
  \end{center}
\end{figure}

The six dominant modes for velocity, pressure and supremizers are shown. See
Figures~\ref{fig:modes-vel-noconf}--\ref{fig:modes-sup-noconf} for the regular method and
Figures~\ref{fig:modes-vel-conf}--\ref{fig:modes-sup-conf} for the conforming method. As is common
with reduced basis methods, we can see that the first modes represent the most dominant flow
patterns, and that it is left for the higher order modes to represent finer details, in this case
wake effects. It is interesting that the first velocity modes of the two methods do not agree (in
fact, the second mode of the conforming method agrees with the first mode of the regular
method). The first two supremizer modes of the two methods also appear to be ``switched''. As for
the pressure modes, aside from irrelevant sign flips, they are in good agreement between the two
methods, which is natural considering that the function space mapping for the pressure is identical.

\begin{figure}
  \begin{tikzpicture}
    \begin{axis}[
      xlabel={Angle of attack ($\varphi$, degrees)},
      ylabel={Divergence ($L^2$-norm)},
      ymode=log,
      width=0.9\textwidth,
      height=0.4\textwidth,
      grid=both,
      axis lines=left,
      legend style={
        at={(0.5, -0.25)},
        anchor=north,
        draw=none,
      },
      legend cell align=left,
      legend columns=2,
      ]
      \addplot[blue, thick, mark=*, mark options={solid}]
      table[x index={0}, y index={2}]{data/airfoil-divs-no-piola.csv};
      \addplot[blue, thick, densely dashed, mark=o, mark options={solid}]
      table[x index={0}, y index={1}]{data/airfoil-divs-no-piola.csv};
      \addplot[red, thick, mark=*, mark options={solid}]
      table[x index={0}, y index={2}]{data/airfoil-divs-piola.csv};
      \addplot[red, thick, densely dashed, mark=o, mark options={solid}]
      table[x index={0}, y index={1}]{data/airfoil-divs-piola.csv};
      \legend{
        Regular mean,
        Regular max,
        Conforming mean,
        Conforming max,
      }
    \end{axis}
  \end{tikzpicture}
  \caption{
    Mean and maximal divergences of the ten first basis functions for the
    regular and conforming methods, as measured in $L^2$-norms, for various
    angles of attack.
  }
  \label{fig:divs}
\end{figure}
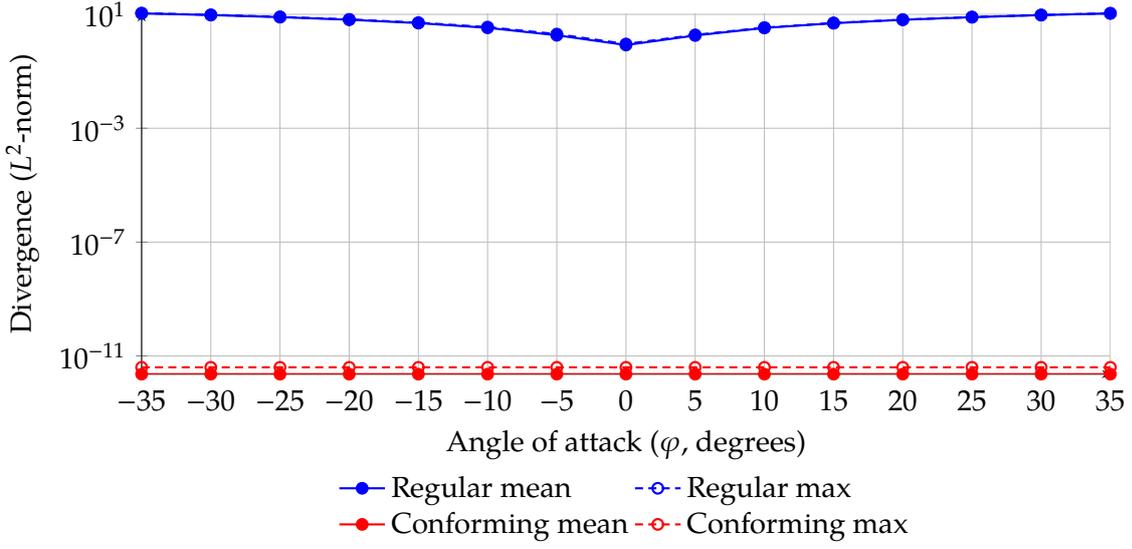

The divergence of the velocity basis functions is shown in Figure~\ref{fig:divs}. Here, the \emph{mean} and \emph{maximal}
divergence of the ten first basis functions for each method (measured in the $L^2$-norm) are shown
for various angles of attack. The figure shows that the divergences of the basis functions for the
conventional method range between approximately~1 and~10. The divergences of the reduced-basis functions
for the conforming method are
consistently zero to $11$ or $12$ digits of accuracy (well within the $10$ digits chosen as a
baseline when choosing the number of terms $n$.) More importantly, it reveals that the solenoidality
property is consistently satisfied throughout the parameter domain, with no significant variability.

\subsection{Stability}
\label{sec:stability}

To understand the stability properties of the four methods, viz.~stabilized and unstabilized Taylor-Hood and
conforming, we computed the LBB constants $\beta_h$ from \eqref{eqn:lbb}. The results are presented
in Table~\ref{tbl:lbb}. The results clearly demonstrate the ability of the supremizer enrichment
procedure to stabilize the method. In addition, the results highlight that no significant degradation
of stability is to be expected for divergence-conforming reduced methods.

Table~\ref{tbl:subspaceangles} shows the minimal principal angles between the velocity and supremizer
spaces, computed according to \cite{Knyazev2002pab}. Since linear independence between the two is
not enforced \emph{a priori}, these angles are relevant to the question whether the sum velocity
space may possibly be degenerate. This shows that for all bases considered in this work, the two
spaces are suitably linearly independent. It is noteworthy that the divergence-conforming bases show
remarkable near-orthogonality for larger $M$ that the Taylor-Hood method does not match.

\begin{table}
  \begin{center}
    \bgroup\def\arraystretch{1.2}
    \begin{tabular}{crrrr}
      & \multicolumn{2}{c}{\bf Taylor-Hood} & \multicolumn{2}{c}{\bf Conforming} \\
      \hline
      $\sharp$ DoFs ($M$)
      & \multicolumn{1}{c}{Unstabilized} & \multicolumn{1}{c}{Stabilized}
      & \multicolumn{1}{c}{Unstabilized} & \multicolumn{1}{c}{Stabilized} \\
      \hline $10$ & $\SI{8.06e-6}{}$ & $\SI{3.93e-1}{}$ & $\SI{2.68e-18}{}$ & $\SI{4.11e-1}{}$ \\
      $20$ & $\SI{5.68e-6}{}$ & $\SI{2.54e-1}{}$ & $\SI{4.05e-18}{}$ & $\SI{3.37e-1}{}$ \\
      $30$ & $\SI{1.03e-6}{}$ & $\SI{3.15e-1}{}$ & $\SI{1.96e-19}{}$ & $\SI{3.13e-1}{}$ \\
      $40$ & $\SI{2.65e-8}{}$ & $\SI{2.67e-1}{}$ & $\SI{1.69e-18}{}$ & $\SI{2.90e-1}{}$ \\
      $50$ & $\SI{1.77e-9}{}$ & $\SI{2.75e-1}{}$ & $\SI{3.17e-18}{}$ & $\SI{2.72e-1}{}$ \\
      \hline
    \end{tabular}
    \egroup
  \end{center}
  \caption{
    LBB constants $\beta_h$, as defined by \eqref{eqn:lbb}, for various basis
    sizes and methods. The reported values are \emph{minima} over a sampling of
    $15 \times 15$ parameter values. The LBB constants for the unstabilized conforming
    method were all indistinguishable from zero to machine precision.
  }
  \label{tbl:lbb}
\end{table}

\begin{table}
  \begin{center}
    \bgroup\def\arraystretch{1.2}
    \begin{tabular}{crr}
      $\sharp$ DoFs ($M$) & {\bf Taylor-Hood} & {\bf Conforming} \\
      \hline $10$ & $\SI{82.2}{\degree}$ & $\SI{89.9}{\degree}$ \\
      \hline $20$ & $\SI{75.4}{\degree}$ & $\SI{88.9}{\degree}$ \\
      \hline $30$ & $\SI{65.7}{\degree}$ & $\SI{88.3}{\degree}$ \\
      \hline $40$ & $\SI{62.7}{\degree}$ & $\SI{86.9}{\degree}$ \\
      \hline $50$ & $\SI{59.4}{\degree}$ & $\SI{85.7}{\degree}$ \\
      \hline
    \end{tabular}
    \egroup
  \end{center}
  \caption{
    Minimal principal angles between the velocity and supremizer spaces for the two methods for various
    reduced basis sizes.
  }
  \label{tbl:subspaceangles}
\end{table}

\subsection{Performance}

Performance metrics are presented in terms of error, speedup and degrees of
freedom.
\begin{itemize}
  \item \emph{Error}, when reported, is always given as
    \emph{mean relative error} between the result of the high-fidelity and the
    reduced basis method. The mean is taken over $225 = 15 \times 15$ uniformly
    spaced points in the parameter space. The error norms used are $H^1$
    seminorm for velocities and $L^2$ norm for pressure. For the
    unstabilized reduced methods, we found that they consistently failed to
    converge, thus the errors reported for those methods are based on fewer
    parameter values (those where convergence could be achieved).
  \item \emph{Expected error} is, given $M$, the smallest $\epsilon$ satisfying
    \eqref{eqn:error}. We used prescribed values of $M=10,20,30,40,50$,
    computing $\epsilon$ as a function of $M$ rather than the other way around.
  \item \emph{Time usage} is the mean time spent solving one instance of the problem. This includes
    time spent in both the assembly and the linear solver phase, and notably will also be affected
    by the number of \emph{iterations} the nonlinear solver needs to converge.
  \item \emph{Degrees of freedom} is the number $M$ of basis functions in a reduced space for a
    single field. For the regular method, this implies a total of $2M$ degrees of freedom in
    velocity (although only $M$ of these can be expected to have approximative power, the supremizer
    functions are indispensable and \emph{do} contribute to the final solution), and $M$ degrees of
    freedom in pressure.  For the conforming method, there are $M$ degrees of freedom in both
    spaces.
\end{itemize}

Figure~\ref{fig:perf1-unstab} first shows measured error as a function of expected error, for the
stabilized and unstabilized Taylor-Hood methods. While the stabilized method can be seen to
converge, the velocity solution of the un-stabilized methods does not, and the pressure solution
diverges. In addition to this, the Taylor-Hood un-stabilized method has serious convergence
problems in the nonlinear solver, especially for larger $M$, and the conforming un-stabilized
methods essentially never converges for any parameter values. In the following we will therefore
ignore the unstabilized methods.

Figure~\ref{fig:perf1} shows the same results, comparing the regular and the conforming methods
(both stabilized). We see a clear linear relationship between the expected and measured errors. The
conforming method obtains somewhat better velocity results and comparable pressure results to the
Taylor-Hood method.

Figure~\ref{fig:perf2} shows error as a function of degrees of freedom. Here, it is clear that the
pressure results for the conventional method are slightly better \emph{per degree of freedom}. The
difference in the pressure errors is more pronounced than in Figure~\ref{fig:perf1} because part
of the discrepancy can be attributed to corresponding differences in the spectrum; cf. Figure~\ref{fig:spectra}.
That is to say, due to the slightly slower decay
of eigenvalues in the conforming pressure snapshots as compared to the Taylor-Hood modes, the
expected error for reduced models with the same number of degrees of freedom will be higher for the
latter method.
In Figure~\ref{fig:perf2} we can also see that the combined pressure-velocity basis
(Section~\ref{sec:combined}) has reconstructed pressure errors that are only slightly worse than the
recovered pressures from Section~\ref{sec:conforming}. It is to be noted, however, that it is not
evident that these good approximation properties of the reconstructed pressure will be retained for
more strongly nonlinear cases.

\begin{figure}
  \begin{tikzpicture}
    \begin{axis}[
      xlabel={Time (seconds)},
      ylabel={Mean relative error},
      ymode=log,
      xmode=log,
      width=0.9\textwidth,
      height=0.5\textwidth,
      grid=both,
      axis lines=left,
      legend style={
        at={(0.5, -0.2)},
        anchor=north,
        draw=none,
      },
      legend cell align=left,
      legend columns=2,
      ]
      \addplot[blue, thick, mark=*, mark options={solid}]
      table[x index={15}, y index={4}]{data/airfoil-results-no-piola-sups-no-block.csv};
      \addplot[blue, thick, densely dashed, mark=o, mark options={solid}]
      table[x index={15}, y index={8}]{data/airfoil-results-no-piola-sups-no-block.csv};
      \addplot[red, thick, mark=*, mark options={solid}]
      table[x index={15}, y index={4}]{data/airfoil-results-piola-sups-no-block.csv};
      \addplot[red, thick, densely dashed, mark=o, mark options={solid}]
      table[x index={15}, y index={8}]{data/airfoil-results-piola-sups-no-block.csv};
      \addplot[green, thick, mark=*, mark options={solid}]
      table[x index={11}, y index={4}]{data/airfoil-results-piola-sups-block.csv};
      \addplot[green, thick, densely dashed, mark=o, mark options={solid}]
      table[x index={11}, y index={8}]{data/airfoil-results-piola-sups-block.csv};
      \addplot[cyan, thick, mark=*, mark options={solid}]
      table[x index={11}, y index={4}]{data/airfoil-results-combined.csv};
      \addplot[cyan, thick, densely dashed, mark=o, mark options={solid}]
      table[x index={11}, y index={8}]{data/airfoil-results-combined.csv};
      \legend{
        Regular stabilized ($v$),
        Regular stabilized ($p$),
        Conforming stabilized naive ($v$),
        Conforming stabilized naive ($p$),
        Conforming stabilized block ($v$),
        Conforming stabilized block ($p$),
        Conforming combined ($v$),
        Conforming combined ($p$),
      }
    \end{axis}
  \end{tikzpicture}
  \caption{
    Measured mean relative error as a function of mean time usage (in seconds), for velocity
    ($H^1$-seminorm) and pressure ($L^2$-norm). The time reported includes both the velocity
    solution and the pressure recovery or reconstruction, as appropriate. Both regular and
    conforming are shown, the latter with two different solvers (the naive solver and the block
    solver), as well as with the joint velocity-pressure reconstruction scheme. Bottom left is
    better.
  }
  \label{fig:perf3}
\end{figure}
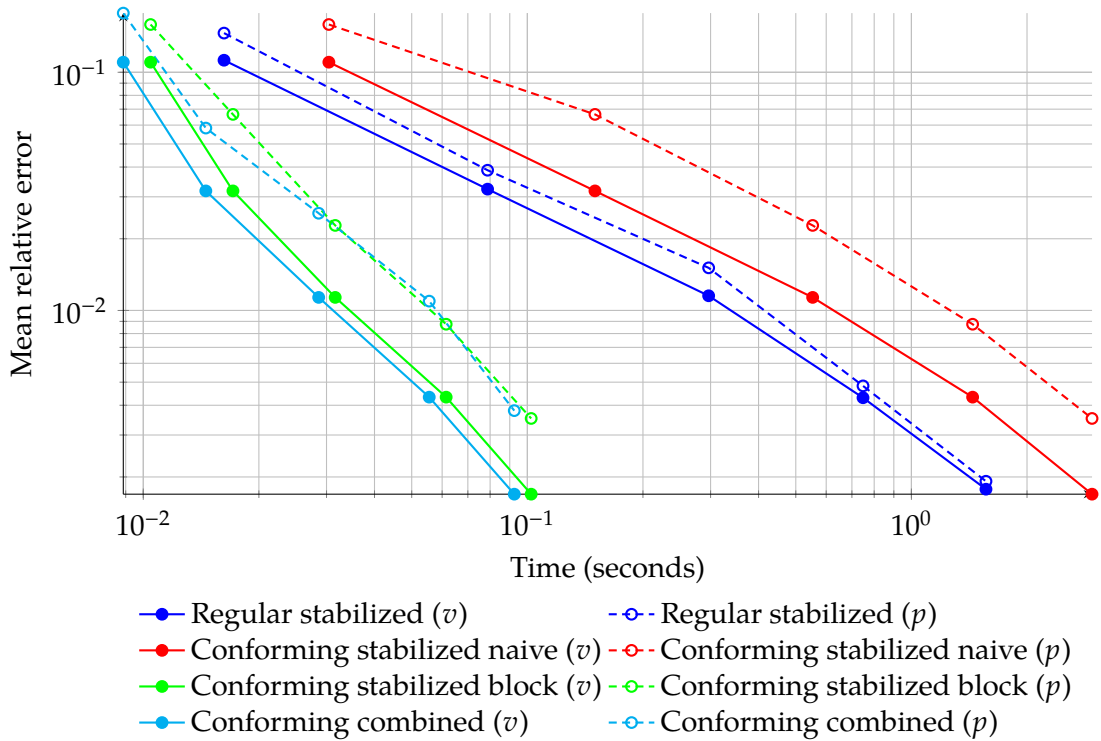

Figure~\ref{fig:perf3} shows the measured errors as a function of mean time usage, measured in
seconds. The time usage reported is the total time taken for both velocity solution and pressure
recovery, as appropriate.
This clearly shows the ability of the conforming method to achieve the same results significantly
faster if the block solver algorithm is employed. It also shows how the combined velocity-pressure
basis is able to achieve the same velocity results slightly faster, on account of not requiring an
additional linear system to solve for the pressure.

\begin{figure}
  \begin{tikzpicture}
    \begin{axis}[
      xlabel={Expected mean relative error},
      ylabel={Mean relative error},
      ymode=log,
      xmode=log,
      width=0.9\textwidth,
      height=0.5\textwidth,
      grid=both,
      axis lines=left,
      legend style={
        at={(0.5, -0.2)},
        anchor=north,
        draw=none,
      },
      legend cell align=left,
      legend columns=2,
      ]
      \addplot[blue, thick, mark=*, mark options={solid}]
      table[x index={1}, y index={4}]{data/airfoil-results-no-piola-sups-no-block.csv};
      \addplot[blue, thick, densely dashed, mark=o, mark options={solid}]
      table[x index={2}, y index={8}]{data/airfoil-results-no-piola-sups-no-block.csv};
      \addplot[magenta, thick, mark=*, mark options={solid}]
      table[x index={1}, y index={4}]{data/airfoil-results-no-piola-no-sups-no-block.csv};
      \addplot[magenta, thick, densely dashed, mark=o, mark options={solid}]
      table[x index={2}, y index={6}]{data/airfoil-results-no-piola-no-sups-no-block.csv};
      \legend{
        Regular stabilized ($v$),
        Regular stabilized ($p$),
        Regular un-stabilized ($v$),
        Regular un-stabilized ($p$),
      }
    \end{axis}
  \end{tikzpicture}
  \caption{
    Measured mean relative error as a function of expected mean relative error,
    for velocity ($H^1$-seminorm) and pressure ($L^2$-norm). The regular
    stabilized and un-stabilized methods are shown. Stabilization is clearly
    necessary to obtain any sort of useful solution.
  }
  \label{fig:perf1-unstab}
\end{figure}
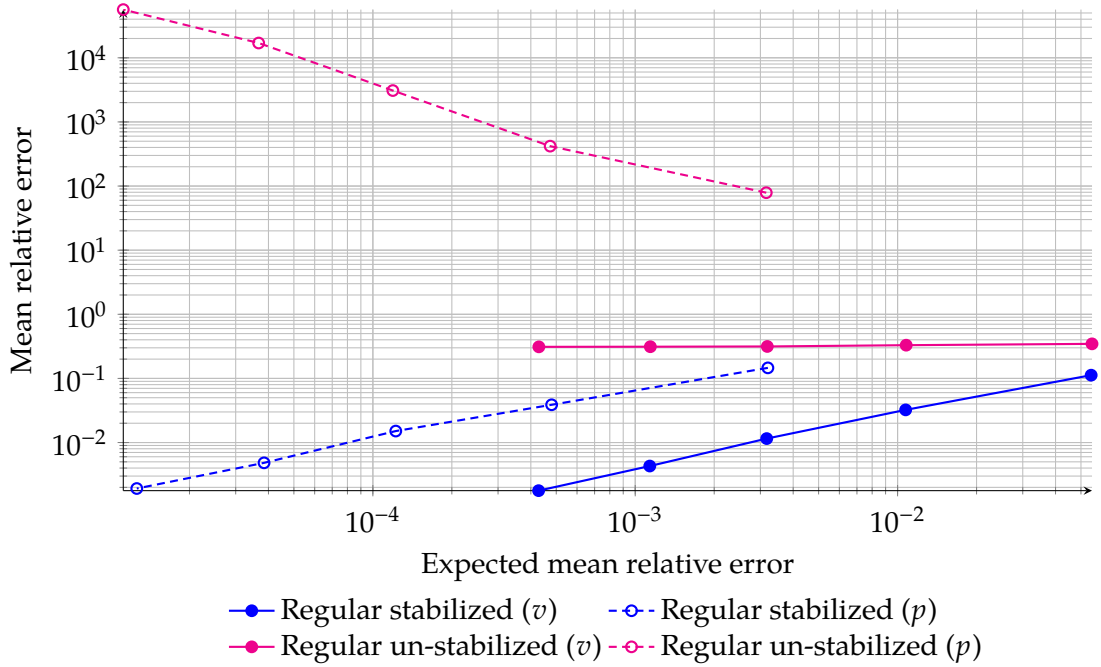

\begin{figure}
  \begin{tikzpicture}
    \begin{axis}[
      xlabel={Expected mean relative error},
      ylabel={Mean relative error},
      ymode=log,
      xmode=log,
      width=0.9\textwidth,
      height=0.5\textwidth,
      grid=both,
      axis lines=left,
      legend style={
        at={(0.5, -0.2)},
        anchor=north,
        draw=none,
      },
      legend cell align=left,
      legend columns=2,
      ]
      \addplot[blue, thick, mark=*, mark options={solid}]
      table[x index={1}, y index={4}]{data/airfoil-results-no-piola-sups-no-block.csv};
      \addplot[blue, thick, densely dashed, mark=o, mark options={solid}]
      table[x index={2}, y index={8}]{data/airfoil-results-no-piola-sups-no-block.csv};
      \addplot[red, thick, mark=*, mark options={solid}]
      table[x index={1}, y index={4}]{data/airfoil-results-piola-sups-no-block.csv};
      \addplot[red, thick, densely dashed, mark=o, mark options={solid}]
      table[x index={2}, y index={8}]{data/airfoil-results-piola-sups-no-block.csv};
      \legend{
        Regular stabilized ($v$),
        Regular stabilized ($p$),
        Conforming stabilized ($v$),
        Conforming stabilized ($p$),
      }
    \end{axis}
  \end{tikzpicture}
  \caption{
    Measured mean relative error as a function of expected mean relative error,
    for velocity ($H^1$-seminorm) and pressure ($L^2$-norm). Both regular and
    conforming are shown.
    The conforming block solver is indistinguishable from the naive conforming solver,
    and the velocity errors of the Taylor-Hood and the conforming methods are also virtually
    indistinguishable from each other.
    Bottom right is better.
  }
  \label{fig:perf1}
\end{figure}
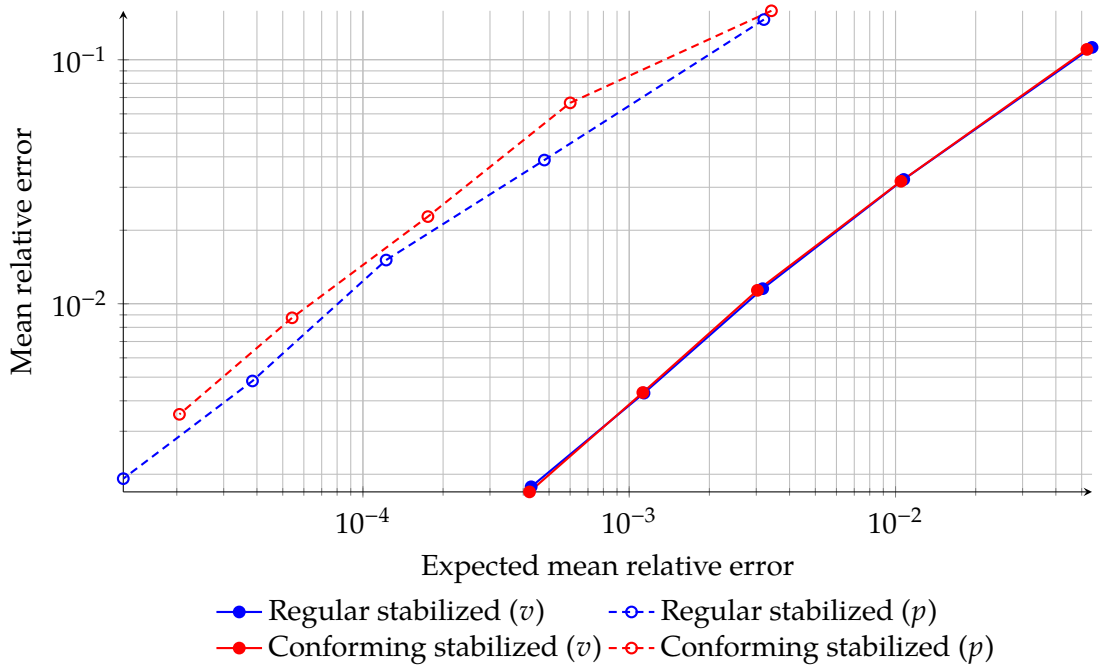

\begin{figure}
  \begin{tikzpicture}
    \begin{axis}[
      xlabel={Degrees of freedom ($M$)},
      ylabel={Mean relative error},
      ymode=log,
      width=0.9\textwidth,
      height=0.5\textwidth,
      grid=both,
      axis lines=left,
      legend style={
        at={(0.5, -0.2)},
        anchor=north,
        draw=none,
      },
      legend cell align=left,
      legend columns=2,
      ]
      \addplot[blue, thick, mark=*, mark options={solid}]
      table[x index={0}, y index={4}]{data/airfoil-results-no-piola-sups-no-block.csv};
      \addplot[blue, thick, densely dashed, mark=o, mark options={solid}]
      table[x index={0}, y index={8}]{data/airfoil-results-no-piola-sups-no-block.csv};
      \addplot[red, thick, mark=*, mark options={solid}]
      table[x index={0}, y index={4}]{data/airfoil-results-piola-sups-no-block.csv};
      \addplot[red, thick, densely dashed, mark=o, mark options={solid}]
      table[x index={0}, y index={8}]{data/airfoil-results-piola-sups-no-block.csv};
      \addplot[cyan, thick, mark=*, mark options={solid}]
      table[x index={0}, y index={4}]{data/airfoil-results-combined.csv};
      \addplot[cyan, thick, densely dashed, mark=o, mark options={solid}]
      table[x index={0}, y index={8}]{data/airfoil-results-combined.csv};
      \legend{
        Regular stabilized ($v$),
        Regular stabilized ($p$),
        Conforming stabilized ($v$),
        Conforming stabilized ($p$),
        Conforming combined ($v$),
        Conforming combined ($p$),
      }
    \end{axis}
  \end{tikzpicture}
  \caption{
    Measured mean relative error as a function of $M$, for velocity
    ($H^1$-seminorm) and pressure ($L^2$-norm). Both regular and conforming are
    shown. The velocity errors for the two conforming methods overlap.
    Bottom left is better.
  }
  \label{fig:perf2}
\end{figure}
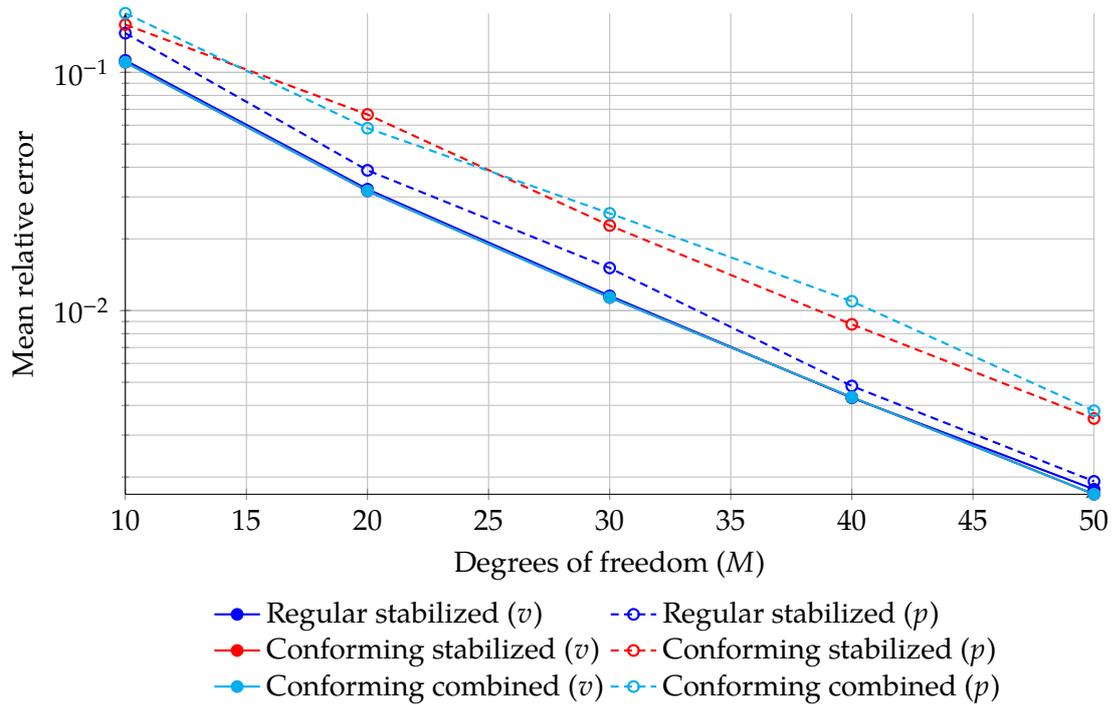

\clearpage

\section{Conclusions}
\label{sec:conc}

We have formulated the fundamentals for Reduced Basis Methods (RBMs) for stationary Navier--Stokes
problems. Using a divergence-conforming isogeometric high-fidelity method, we were able to generate
divergence-free reduced bases, yielding a velocity-only reduced formulation. To solve for pressure,
we employed supremizers as test functions, a technique used to stabilize conventional reduced basis
methods. This allowed us to achieve accurate reduced solutions for both velocity and pressure with
significant speed benefits over a coupled conventional formulation. In addition, we considered the
option of direct pressure reconstruction, using the coefficients of the RB velocity solution.

Investigations were performed for the parametrized problem of flow around a NACA0015 airfoil with
varying angle-of-attack and inflow velocity. The divergence-conforming method was compared with a
supremizer stabilized reduced method based on the Taylor-Hood high-fidelity scheme with the same
mesh and comparable polynomial degree. For the divergence-conforming method, a velocity-only
formulation was applied, with two choices for the pressure approximation, viz.~recovery and
reconstruction.

While the accuracy of the reduced methods were virtually identical, the divergence-conforming models
were able to achieve a significantly faster online stage thanks to the velocity-pressure decoupling.

The pressure recovery and reconstruction methods were also able to achieve similar performance for the
considered problem. Since the pressure reconstruction method relies on a linearity
assumption that is not valid in general, we must conclude that pressure
recovery is to be preferred in most cases.

The faster online stage of the conforming RB method is offset by a more complicated offline stage.
In particular, the affine
representations \eqref{eqn:split-a}--\eqref{eqn:split-d1} of the divergence-conforming method
required approximately three times as many terms (${\sim}12n$ versus ${\sim}4n$). In spite of this,
the conforming method yielded significantly faster online stages. We aim to investigate in future work
in which manner automatic ``black box'' methods such as Empirical Interpolation
(EIM, see~\cite{Barrault2004eim,Grepl2007erb,Maday2007gmi} and~\cite[Chapter 10]{Quarteroni2016rbm})
can alleviate the additional work needed to produce such affine representations.

We will also investigate the performance of the present approach when using adaptive high-fidelity
models based on \emph{a posteriori} error estimates developed in~\cite{Kumar2015sap,Kumar2017spr}.

\section{Acknowledgements}

The authors acknowledge the financial support from the Norwegian Research Council and the industrial
partners of OPWIND: Operational Control for Wind Power Plants (Grant No.: 268044/E20).

\bibliography{refs}
\bibliographystyle{plain}

\end{document}